\documentclass[12pt]{article}
\usepackage{amsmath}
\usepackage{amssymb}
\usepackage{amsthm}
\usepackage[totalwidth=17cm, totalheight=25cm]{geometry}
\usepackage{graphicx}
\usepackage{mathabx}
\usepackage{tabularx}
	\newcolumntype{C}[1]{>{\centering\arraybackslash}m{#1}} 
	\newcolumntype{R}[1]{>{\raggedleft\arraybackslash}m{#1}} 
\newtheoremstyle{boldplain}% name
{9pt}%      Space above
{9pt}%      Space below
{\itshape}%         Body font
{}%         Indent amount (empty = no indent, \parindent = para indent)
{\bfseries}% Thm head font
{.}%        Punctuation after thm head
{.5em}%     Space after thm head: " " = normal interword space;
{\thmname{#1}\thmnumber{ #2}\thmnote{ (#3)}}%

\newtheoremstyle{bolddefinition}% name
{9pt}%      Space above
{9pt}%      Space below
{}%         Body font
{}%         Indent amount (empty = no indent, \parindent = para indent)
{\bfseries}% Thm head font
{.}%        Punctuation after thm head
{.5em}%     Space after thm head: " " = normal interword space;
{\thmname{#1}\thmnumber{ #2}\thmnote{ (#3)}}%

\theoremstyle{boldplain}

\newtheorem{cor}[equation]{Corollary}

\newtheorem{lem}[equation]{Lemma}
\newtheorem{lemma}[equation]{Lemma}
\newtheorem{prop}[equation]{Proposition}

\newtheorem{thm}[equation]{Theorem}

\newtheorem{notation}[equation]{Notation}

\theoremstyle{bolddefinition}
\newtheorem{dfn}[equation]{Definition}
\newtheorem{definition}[equation]{Definition}

\newtheorem{example}[equation]{Example}

\newtheorem{ques}[equation]{Question}

\newtheorem{rem}[equation]{Remark}

\newfont{\bigbf}{cmbx10 scaled\magstep1}
\normalbaselines
\numberwithin{equation}{section}

\def\B{\operatorname{B}}

\def\R{{\mathbb R}}

\def\Z{{\mathbb Z}}

\def\al{\alpha}
\def\ga{\gamma}
\def\Ga{\Gamma}
\def\de{\delta}
\def\De{\Delta}
\def\eps{\epsilon}
\def\la{\lambda}
\def\La{\Lambda}
\def\si{\sigma}
\def\Si{\Sigma}

\def\om{\omega}
\def\Om{\Omega}

\def\3{\ss}

\def\Acc{\operatorname{Acc}}
\def\acts{\curvearrowright}
\def\amod{a_{mod}}

\def\card{\operatorname{card}}
\def\CF{C_{F\ddot u}}

\def\D{\partial}
\def\DF{\partial_{F\ddot u}}

\def\diam{\mathop{\hbox{diam}}}

\def\diamot{\diamondsuit_{\tau_{mod}}}

\def\dist{\operatorname{dist}}
\def\embed{\hookrightarrow}
\def\Fix{\operatorname{Fix}}
\def\Flag{\operatorname{Flag}}
\def\Flagn{\operatorname{Flag_{\nu_{mod}}}}
\def\Flagin{\Flag_{\iota\nu_{mod}}}
\def\Flags{\operatorname{Flag_{\si_{mod}}}}

\def\Flagt{\operatorname{Flag_{\tau_{mod}}}}
\def\Flagit{\Flag_{\iota\tau_{mod}}}

\def\Flagpmt{\Flag_{\pm\tau_{mod}}}

\def\Fmod{F_{mod}}

\def\geo{\partial_{\infty}}

\def\id{\mathop{\hbox{id}}}
\def\im{\mathop{\hbox{im}}}
\def\inte{\operatorname{int}}
\def\interior{\operatorname{int}}

\def\Lan{\La_{\nu_{mod}}}

\def\Lat{\La_{\tau_{mod}}}
\def\Latm{\La_{\tau_{mod}}^-}
\def\Latp{\La_{\tau_{mod}}^+}

\def\Lapmt{\La_{\pm\tau_{mod}}}

\def\lra{\longrightarrow}

\def\numod{\nu_{mod}}

\def\oa{\overrightarrow}
\def\ol{\overline}
\def\OmF{\operatorname{\Om}_{F\ddot u}}

\def\pihalf{\frac{\pi}{2}}

\def\2pithird{\frac{2\pi}{3}}

\def\pos{\mathop{\hbox{pos}}\nolimits}
\def\cpos{\mathop{\hbox{c-pos}}}
\def\ccpos{\mathop{\hbox{c-c-pos}}}
\def\prect{\prec_{\taumod}}

\def\rank{\mathop{\hbox{rank}}}

\def\Ra{\Rightarrow}

\def\simod{\si_{mod}}

\def\st{\operatorname{st}}
\def\stF{\operatorname{st}_{F\ddot u}}

\def\ost{\operatorname{ost}}

\def\Stab{\operatorname{Stab}}

\def\tangle{\angle_{Tits}}
\def\taumod{\tau_{mod}}

\def\tits{\partial_{Tits}}
\def\Th{\mathop{\hbox{Th}}\nolimits}
\def\ThF{\operatorname{Th}_{F\ddot u}}
\def\cTh{\mathop{\hbox{Th}^c}\nolimits}

\def\Wn{W_{\nu_{mod}}}

\def\Wt{W_{\tau_{mod}}}

\def\8{\infty}
\def\<{\langle}
\def\>{\rangle}

\def\BI{\begin{itemize}}
\def\EI{\end{itemize}}

\long\def\comment#1\endcomment{}

\input diagrams.tex

\title{Finsler bordifications of symmetric and certain locally symmetric spaces}
\author{Michael Kapovich, Bernhard Leeb}
\date{March 14, 2016}

\begin{document}

\maketitle

\begin{abstract}
\noindent
We give a geometric interpretation of the maximal Satake compactification 
of symmetric spaces $X=G/K$ of noncompact type, 
showing that it arises by attaching the horofunction boundary for a suitable $G$-invariant Finsler metric on $X$. 
As an application, we establish the existence of natural bordifications, as orbifolds-with-corners, 
of locally symmetric spaces $X/\Gamma$ for arbitrary discrete subgroups $\Gamma< G$. 
These bordifications result from attaching $\Gamma$-quotients of suitable domains of proper discontinuity at infinity. We further prove that such bordifications are compactifications in the case of Anosov subgroups. 
We show, conversely, that Anosov subgroups are characterized by the existence of such compactifications among uniformly regular subgroups. Along the way, we give a positive answer, in the torsion free case, 
to a question of Ha\"issinsky and Tukia on convergence groups
regarding the cocompactness of their actions on the domains of discontinuity.
\end{abstract}

\tableofcontents

\section{Introduction}

The goal of this paper is four-fold:

1.\ We give a geometric interpretation of the maximal Satake compactification 
of a symmetric space $X=G/K$ of noncompact type
by obtaining it as the horoclosure with respect to a suitable $G$-invariant Finsler metric.

2.\ This compactification turns out to have good dynamical properties,
better, for our purposes, than the usual visual compactification as a CAT(0) space.
In it we find natural domains of proper discontinuity for discrete subgroups $\Ga<G$.
For Anosov subgroups we show that the actions on these domains are also cocompact,
thereby providing natural orbifold-with-corner compactifications of the corresponding locally symmetric spaces. 

3.\ We use these dynamical results to establish new characterizations of Anosov subgroups.

4.\ We apply our techniques for proving cocompactness to the theory of abstract convergence groups
and verify the cocompactness on the domain of discontinuity for a certain class of actions. 

We now describe our main results in more detail.

1.\ 
We prove that the {\em maximal Satake compactification} $\ol{X}^S_{max}$  
(see \cite[Chapter 2]{Borel-Ji}) 
is $G$-equivariantly homeomorphic, as a manifold-with-corners, 
to a {\em regular Finsler compactification} $ \ol X^{Fins}=\ol X^{\bar\theta}$ 
obtained by adding to $X$ points at infinity represented by {\em Finsler horofunctions}. 
These horofunctions arise as limits, modulo additive constants, of distance functions
$$
d^{\bar\theta}_x = d^{\bar\theta}(\cdot,x)
$$
where $d^{\bar\theta}$ is a certain $G$-invariant Finsler distance on $X$ associated 
with an interior point $\bar\theta$ of the model spherical Weyl  chamber $\simod$ of $X$. 
This {\em horoclosure} construction is a special case of a well-known general construction 
of compactifications for metric spaces.
For instance, applying it to CAT(0) spaces yields their visual compactification. 
The novelty here is finding the right metric on the symmetric space $X$ which yields $\ol{X}^S_{max}$. 
Our first main result, 
proven in sections~\ref{sec:fico}, \ref{sec:cox} and~\ref{sec:mfcr}, 
describes geometric and dynamical properties of the Finsler compactification:

\begin{thm}
\label{thm:comp}
For every regular type $\bar\theta\in\interior(\simod)$,
$$\ol X^{\bar\theta}=X\sqcup\geo^{\bar\theta}X$$ 
is a compactification of $X$ as a $G$-space 
which satisfies the following properties:

(i) There are finitely many $G$-orbits $S_{\taumod}$ indexed by the faces $\taumod$ of $\simod$.
($X=S_{\emptyset}$.)

(ii) The stratification of $\ol X^{\bar\theta}$ by $G$-orbits is a $G$-invariant manifold--with--corners structure. 

(iii) 
$\ol X^{\bar\theta}$ is homeomorphic to the closed ball,
with $X$ corresponding to the open ball.

(iv) The compactification $\ol X^{\bar\theta}$ is independent of the regular type $\bar\theta$ 
in the sense that the identity map $\id_X$ extends to a natural homeomorphism of any two such compactifications. 

(v) There exists a $G$-equivariant homeomorphism of manifolds with corners 
between $\ol X^{\bar\theta}$ and the 
maximal Satake compactification $\ol{X}^S_{max}$ which yields a natural correspondence of strata.
\end{thm}
In view of (iv) we will denote the Finsler compactification from now on by $\ol X^{Fins}$. 

\begin{rem}
(i) We also give a geometric interpretation of the points in $\geo^{Fins}X$ as strong asymptote classes of Weyl sectors,
see Lemma~\ref{lem:distlifct}. 

(ii) The strata $S_{\taumod}\subset \geo^{Fins} X$ at infinity 
naturally fiber over the partial flag manifolds $\Flagt\cong G/P_{\taumod}$.
The fibers $X_{\tau}$ for $\tau\in\Flagt$, called {\em small strata},
are naturally identified with symmetric subspaces of $X$,
namely with cross sections of parallel sets. 
In the case $\taumod=\simod$ the fibration is a homeomorphism,
$S_{\simod}\cong\Flags=\DF X\cong G/B$,
i.e.\ the Furstenberg boundary $\DF X$ embeds at infinity as the unique closed stratum.

(iii) The Finsler view point had emerged in several instances during our earlier study 
\cite{coco15,morse,mlem}
of asymptotic and coarse properties 
of regular discrete isometry groups acting on symmetric spaces and euclidean buildings.
For instance, the notion of {\em flag convergence} 
(defined earlier in \cite[\S 7.2]{coco15})  
is a special case of the Finsler convergence at infinity considered in this paper, see Proposition \ref{prop:relconvflfi}. 
Furthermore, the {\em Morse Lemma} proven in \cite{mlem} can be rephrased to the effect 
that regular quasigeodesics in symmetric spaces and euclidean buildings 
are uniformly close to Finsler geodesics,  see section \ref{sec:geodesics}. 
In the same vein, 
{\em Morse subgroups} $\Ga<G$ can be characterized as Finsler quasiconvex, 
see section~\ref{sec:fqc}. 

(iv) The maximal Satake compactification is known to carry a $G$-invariant real-analytic structure, see \cite{Borel-Ji}.
\end{rem}

\begin{rem}
After finishing this work we learnt about work of Anne Parreau \cite{Parreau}
where she studies the geometry of CAT(0) model spaces, 
i.e.\ of symmetric spaces of noncompact type and euclidean buildings, 
from a very natural perspective, regarding them as metric spaces 
with a vector valued distance function with values in the euclidean Weyl chamber $\De$
(called $\De$-distance in our paper). 
Among other things, she shows that basic properties of CAT(0) spaces persist in this setting, 
notably the convexity of the distance, 
and develops a comparison geometry for the $\De$-distance function.  
Furthermore, she proves that the resulting  $\De$-valued horofunction
compactifications of model spaces 
are naturally homeomorphic to 
their maximal Satake compactifications.
\end{rem}

2. Our main application of Theorem \ref{thm:comp} concerns discrete subgroups $\Ga<G$. Recall that 
if $X$ is a negatively curved symmetric space, then the locally symmetric space $X/\Ga$ (actually, an orbifold) admits the standard bordification
$$
X/\Ga\embed (X\sqcup \Om(\Ga))/\Ga
$$
where $\Om(\Ga)\subset \geo X$ is the domain of discontinuity of $\Ga$ at infinity. 
The quotient $(X\sqcup \Om(\Ga))/\Ga$ is an orbifold with boundary $\Om(\Ga)/\Ga$. 
Furthermore, a subgroup $\Ga$ is {\em convex cocompact} if and only if $(X\sqcup \Om(\Ga))/\Ga$ is compact. 
The main purpose of this paper is to generalize these bordifications and compactifications 
to suitable classes of discrete subgroups of higher rank Lie groups. 

In our earlier papers \cite{coco15, morse, mlem}, 
we introduced several conditions for discrete subgroups $\Ga$ of semisimple Lie groups $G$, 
generalizing the notions of discreteness and convex cocompactness in rank one. 
These properties are defined relative to faces $\taumod$ of the spherical model Weyl chamber $\simod$,
equivalently, with respect to conjugacy classes of parabolic subgroups of $G$. 
The most important properties for the purposes of this paper are {\em regularity} 
and {\em asymptotic embeddedness}.

The regularity conditions capture the asymptotics of the orbits in $X$
and are reflected by the location of their accumulation sets in $\geo^{Fins}X$.
A discrete subgroup $\Ga<G$ is {\em $\taumod$-regular}
if its orbits $\Ga x\subset X$ accumulate in $\geo^{Fins}X$ at the closure of the stratum $S_{\taumod}$.
The {\em $\taumod$-limit set} $\Lat$ of $\Ga$ is then defined as the compact set of simplices $\tau\in\Flagt$
such that the small stratum closure $\ol X_{\tau}$ contains accumulation points. 
The subgroup $\Ga$ is {\em $\taumod$-antipodal} 
if the simplices in $\Lat$ are pairwise antipodal.\footnote{Here we require 
the face $\taumod$ to be invariant under the opposition involution of $\simod$.
For the correponding parabolic subgroups this means that they are conjugate to their opposite parabolic subgroups.}
In rank one regularity is equivalent to discreteness. 

Asymptotically embedded subgroups form a certain subclass of regular subgroups,
which turns out to coincide with the class of Anosov subgroups, 
cf.\ \cite{morse}. 
A discrete subgroup $\Ga<G$ is {\em $\taumod$-asymptotically embedded}
if it is $\taumod$-regular, $\taumod$-antipodal,
intrinsically word hyperbolic 
and its Gromov boundary $\geo\Ga$ is equivariantly homeomorphic to $\Lat$.\footnote{cf.\ Definition~\ref{def:asyemb}}

In order to obtain bordifications of locally symmetric spaces $X/\Ga$,
we construct domains of proper discontinuity for $\Ga$ in $\ol X^{Fins}$. 
These domains will depend on an auxiliary combinatorial datum,
namely a subset $\Th\subset W$ of the Weyl group, called a {\em thickening}. 
It can be thought of as a set of ``sufficiently special'' relative positions 
of pairs of chambers (full flags) in the Tits boundary. 
This datum is used to construct 
the {\em Finsler thickening} $\Th^{Fins}(\Lat)\subset\geo^{Fins}X$ of $\Lat$
a certain $\Ga$-invariant saturated\footnote{A 
subset $S\subset \geo^{Fins}X$ is called {\em saturated} if it is a union of small strata. }
compact subset, 
see section~\ref{sec:thickenings},
where the reader also finds the definitions of fat and balanced thickenings. 
(Balanced implies fat.)

The following result establishes the existence of natural bordifications and compactifications 
(as orbifolds-with-corners) 
for locally symmetric spaces $X/\Ga$ by attaching 
$\Ga$-quotients of suitably chosen saturated domains in the Finsler boundary of $X$. 
It is a combination of Theorems~\ref{thm:pdwconv} and~\ref{thm:main-coco}.
\begin{thm}
\label{thm:grp2}
Let $\Ga<G$ be a $\taumod$-regular subgroup. Then: 

(i) For each $W_{\taumod}$-left invariant fat thickening $\Th\subset W$, the action 
$$
\Ga \acts X\sqcup\Om_{Th}^{Fins}:= \ol{X}^{Fins} - \Th^{Fins}(\Lat) 
$$
is properly discontinuous. The quotient
\begin{equation}
\label{eq:bordf}
\left(X\sqcup\Om_{Th}^{Fins}\right)/\Ga 
\end{equation}
provides a real-analytic bordification of the orbifold $X/\Ga$
as an orbifold-with-corners.

(ii) If $\Ga$ is $\taumod$-asymptotically embedded
and $\Th$ is a $\Wt$-left invariant balanced thickening,
then $\left(X\sqcup\Om^{Fins}_{Th}\right)/\Ga$ is compact.  In
particular, the bordification is a compactification.
\end{thm}

\begin{rem}
(i)
We note that if $\Ga$ is $\taumod$-antipodal, $\Th$ balanced and $\rank(X)\geq2$,
then the domains $\Om^{Fins}_{Th}$ at infinity are non-empty,
cf.\ Proposition~\ref{prop:connected}.

(ii)
The construction of the domains of proper discontinuity extends in a straightforward way
to all discrete subgroups $\Ga<G$, 
see Theorems~\ref{thm:pdwdsc} and~\ref{thm:pdwdsc2}.

(iii)
We note that 
the existence of an orbifold-with-boundary compactification of locally symmetric quotients by 
Anosov subgroups of some special classes of simple Lie groups (namely, $Sp(2n,\R), SU(n,n), SO(n,n)$) 
appeared in \cite{GW}.

(iv)
The main results of this e-print are already contained in its second version. 
After this work had been completed, the e-print \cite{shameless} was posted,
also addressing the compactification of locally symmetric spaces. 
Theorem~1.2 there provides orbifolds-with-corners compactifications of $X/\Ga$ 
via the maximal Satake compactification of $X$
for $\taumod$-Anosov subgroups $\Ga<G$ of special face types $\taumod$.
However, 
Theorem~1.1 of \cite{shameless} dealing with the general case
still lacks a complete proof. 
It remains unclear 
whether the compactifications constructed there
are orbifolds-with-corners.
Namely,
the approach uses ``generalized'' Satake compactifications
and the proof of Lemma A9, establishing that the latter are manifolds-with-corners,
lacks details. 
Note also that in the first version of \cite{shameless}
there was a mistake in the compactness argument.
It was corrected in the third version 
using methods from the first version of \cite{coco15}.
\end{rem}

3.\ 
As already mentioned,
the class of asymptotically embedded subgroups coincides with the class of Anosov subgroups.
We also prove a converse of part (ii) of the previous theorem,
thereby providing a new characterization of Anosov subgroups among uniformly regular subgroups 
in terms of the existence of certain compactifications of the locally symmetric spaces.
To this end,
we say that a discrete subgroup $\Ga< G$ is {\em S-cocompact} if there exists a
$\Ga$-invariant saturated open subset $\Om\subset \geo^{Fins}X$ such that $\Ga$ acts properly discontinuously and cocompactly on $X\sqcup \Om$. 
Theorem \ref{thm:grp2} shows that $\taumod$-asymptotically embedded subgroups 
are S-cocompact with $\Om= \Om^{Fins}_{Th}$.
Conversely, we prove for $\iota$-invariant face types $\taumod$,
see section~\ref{sec:sccrt}:

\begin{thm}\label{thm:SCC->RCA}
Uniformly $\taumod$-regular S-cocompact subgroups $\Ga<G$ are $\taumod$-Anosov.
\end{thm}
Combining the last two theorems, we obtain the characterization:
\begin{cor}
\label{cor:rcachar}
A uniformly $\taumod$-regular subgroup $\Ga< G$ is $\taumod$-Anosov
iff it is S-cocompact. 
\end{cor}
Our cocompactness results thus provide a precise higher-rank analogue of the characterization of convex cocompact subgroups of rank 1 Lie groups in terms of compactifications of the corresponding locally symmetric spaces.

While proving Theorem \ref{thm:SCC->RCA}, we establish yet another coarse-geometric characterization 
of $\taumod$-Anosov subgroups,
namely as uniformly $\taumod$-regular subgroups 
which are {\em coarse retracts}, see sections \ref{sec:coarse_geometry} and \ref{sec:ecRCA} for the details. This theorem is a higher-rank analogue of the characterization  of quasiconvex subgroups of Gromov-hyperbolic groups as coarse retracts. 

\medskip
4.\ 
In section~\ref{sec:haiss},  
as an intermediate step in the proof of Theorem \ref{thm:grp2}, 
we verify a conjecture by Ha\"issinsky and Tukia 
regarding the cocompactness of convergence group actions on their domains of discontinuity
under mild extra assumptions:

\begin{thm}
Let $\Ga\acts \Si$ be a convergence group action 
of a virtually torsion-free hyperbolic group on a metrizable compact space $\Si$, 
and suppose that $\La\subset \Si$ is an invariant compact subset 
which is equivariantly homeomorphic to $\geo\Ga$. 
Then the action $\Ga\acts \Si - \La$ is cocompact 
provided that $\Si - \La$ has finitely many path connected components.
\end{thm}

{\bf Acknowledgements.} 
We thank the MSRI program ``Dynamics on Moduli Spaces of Geometric Structures'' for support
during part of the work on this project.
We are grateful to the Korea Institute for Advanced Study
for support and excellent working conditions.
The first author was also supported by NSF grant  DMS-12-05312. We are grateful to Lizhen Ji for helpful discussions 
and to Anne Parreau for informing us about her work \cite{Parreau}.

\section{Preliminaries}

\subsection{Notations and definitions}

We note that for Hausdorff paracompact topological spaces (and in this paper we will be dealing only with such topological spaces), Alexander-Spanier and \v{C}ech cohomology theories are naturally isomorphic, see \cite[Ch. 6.9]{Spanier}. 
Therefore, in our paper, all cohomology is Alexander-Spanier-\v{C}ech with field coefficients (the reader can assume that the field of coefficients is $\Z_2$). For manifolds and CW complexes, singular and cellular cohomology is naturally isomorphic to the \v{C}ech cohomology. We will use the notation $H^*_c$ for cohomology with compact support. 
As for homology, we will use it again with field coefficients and only for locally-finite CW complexes, where we will be using singular homology and singular locally finite homology, denoted $H^{lf}_*$. By Kronecker duality, for each locally-finite CW complex $X$,
$$
(H^{lf}_k(X))^* \cong H^k_c(X), \quad k\ge 0. 
$$

\medskip 
We refer the reader to  \cite{Joyce} 
for the definitions of manifolds and orbifolds with corners. The only examples of orbifolds with corners which appear in this paper are the {\em good} ones, i.e., quotients of manifolds with corners by properly discontinuous group actions.

\medskip
Throughout the paper, $\angle$ will denote the angle between vectors in a euclidean vector space,
respectively, the angular metric on spherical simplices.

\subsection{Some point set topology}

Let $Z$ and $Z'$ be first countable Hausdorff spaces,
and let $O\subset Z$ and $O'\subset Z'$ be dense open subsets.
Let $f:Z\to Z'$ be a map such that $f(O)\subseteq O'$, 
and suppose that $f$ has the following partial continuity property:
If $(y_n)$ is a sequence in $O$ which converges to $z\in Z$,
then $f(y_n)\to f(z)$ in $Z'$. 
In particular, $f|_O$ is continuous.

\begin{lem}\label{lem:topologylemma}
Under these assumptions, the map $f$ is continuous.
\end{lem}
\proof The lemma follows from a standard diagonal subsequence argument.
\qed

\medskip
Let $(A_n)$ be a sequence of subsets of a metrizable topological space $Z$.
We denote by $\Acc((A_n))$ the closed subset consisting of the accumulation points of all 
sequences $(a_n)$ of points $a_n\in A_n$.

We say that the sequence of subsets $(A_n)$ 
{\em accumulates at} a subset $S\subset Z$
if $\Acc((A_n))\subseteq S$. 

If $Z$ is compact and $C\subset Z$ is a closed subset,
then the sequence $(A_n)$ accumulates at $S$ 
if and only if every neighborhood $U$ of $C$ contains all but finitely many of the subsets $A_n$.

\subsection{Properness and dynamical relation}\label{sec:discdynrl}

In the paper we will use the notion of {\em dynamical relation} between points of a topological space $Z$,
which is an open subset of a compact metrizable space,
with respect to a topological action $\Ga\acts Z$ of a discrete group. The reader will find this definition in \cite{Frances}, see also \cite{coco15}. 
We write 
$$\xi\stackrel{\Ga}{\sim}\xi'$$
if the points $\xi$ and $\xi'$ are dynamically related with respect to the action of the group $\Ga$,
and 
$$\xi\stackrel{(\ga_n)}{\sim}\xi'$$
if $\xi$ is dynamically related to $\xi'$ with respect to a sequence $\ga_n\to\infty$ in $\Ga$,
cf \cite[\S 2.1]{coco15}.
An action is properly discontinuous if and only if no points of $Z$ are dynamically related to each other, see \cite{Frances}.

\subsection{A transformation group lemma}

Let $K$ be a compact Hausdorff topological group,
and let $K\acts Y$
be a continuous action on a compact Hausdorff space $Y$. 
We suppose that there exists a {\em cross section} for the action,
i.e. a compact subset $C\subset Y$ 
which contains precisely one point of every orbit.

Consider the natural surjective map
\begin{equation*} 
K\times C\stackrel{\al}{\lra} Y 
\end{equation*}
given by the action, 
$\al(k,y)=ky$. 
We observe that $Y$ carries the quotient topology with respect to $\al$,
because $K\times C$ is compact and $Y$ is Hausdorff. 
The identifications by $\al$ are determined by the stabilizers of the points in $C$,
namely $\al(k,y)=\al(k',y')$ iff $y=y'$ and $k^{-1}k'\in\Stab_K(y)$.

Consider now two such actions 
$K\acts Y_1$ and $K\acts Y_2$ by the same group 
with cross sections $C_i\subset Y_i$,
and suppose that 
$$ C_1\stackrel{\phi}{\lra} C_2 $$
is a homeomorphism. 
\begin{lem}\label{lem:transformation lemma}
If $\phi$ respects point stabilizers,
i.e. $\Stab_K(y_1)=\Stab_K(\phi(y_1))$ for all $y_1\in C_1$,
then $\phi$ extends to a $K$-equivariant homeomorphism $\Phi:Y_1\to Y_2$.
\end{lem}
\proof
According to the discussion above,
the stabilizer condition implies that there exists a bijection $\Phi:Y_1\to Y_2$
for which the diagram 
$$
\begin{diagram}
K\times C_1 & \rTo^{id_K \times\phi}           &    K\times C_2  \\
   \dTo^{\al_1}         &                 &            \dTo^{\al_2} \\
 Y_1 & \rTo^{\Phi}           &    Y_2 \\
\end{diagram}
$$
commutes.
Since the $\al_i$ are quotient projections,
it follows that $\Phi$ is a homeomorphism.
\qed

\subsection{Thom class} 

In this section $H^{lf}_*$ denotes locally finite homology with $\Z_2$-coefficients. 

\begin{lem}[Thom class] 
\label{lem:thom}
Let $F\stackrel{\iota}{\to} E\to B$ be a fiber bundle 
whose base $B$ is a compact CW-complex and whose fiber $F$ 
is a connected $m$-manifold (without boundary). 
Suppose that there exists a section $s:B\to E$.
Then the map
$$
\underbrace{H^{lf}_m(F)}_{\cong\Z_2}\stackrel{\iota_*}{\lra} H^{lf}_m(E)
$$
induced by an inclusion of the fiber 
is nonzero. 
\end{lem}
\proof By thickening the section, one obtains a closed disk subbundle $D\to B$.
Then we have the commutative diagram:
\begin{diagram}
 H^{lf}_m(F) & \rTo^{\iota_*}           &    H^{lf}_m(E)  \\
   \dTo^{j}         &                 &            \dTo \\
 H_m(D_F,\D D_F) & \rTo^{\iota'_*}           &    H_m(D, \D D)\\
\end{diagram}
The map $j$ is an isomorphism.
By Thom's isomorphism theorem, the map $\iota'_*$ is injective. 
It follows that the map $\iota_*$ is injective.
\qed

\subsection{The horoboundary of metric spaces}
\label{sec:horob}

We refer the reader to \cite{Gromov}, \cite[ch. II.1]{Ballmann} for the definition and basic properties of horofunction compactification of metric spaces. In this section we describe these notions  
in the context of {\em nonsymmetric metrics},
compare \cite{Walsh14}.

Let $(Y,d)$ be a metric space.
We allow the distance $d$ to be {\em non-symmetric},
i.e. we only require that it is positive,
$$ d(y,y')\geq0 \hbox{ with equality iff } y=y' ,$$
and satisfies the triangle inequality
$$ d(y,y')+d(y',y'')\geq d(y,y'') .$$
The symmetrized distance 
$$ d^{sym}(y,y') := d(y,y')+d(y',y) $$
is a metric in the standard sense
and induces a {\em topology} on $Y$.
One observes that $d$ is continuous,
and the distance functions 
$$ d_y := d(\cdot,y) $$
are 1-Lipschitz with respect to $d^{sym}$. These functions satisfy  the inequality 
\begin{equation}
\label{ineq:osc}
-d(y,y') \leq d_y -d_{y'} \leq d(y',y) . 
\end{equation}
Let ${\mathcal C}(Y)$ denote the space of continuous real valued functions,
equipped with the topology of uniform convergence on bounded subsets.
Moreover, let 
$$ \ol{\mathcal C}(Y) := {\mathcal C}(Y)/\R $$
be the quotient space of continuous functions modulo additive constants.
We will denote by $[f]\in\ol{\mathcal C}(Y)$ 
the equivalence class represented by a function $f\in{\mathcal C}(Y)$,
and our notation $f\equiv g$ means that the difference $f-g$ is constant. 

We consider the natural map 
\begin{equation}
\label{eq:embfct}
Y\lra \ol{\mathcal C}(Y), \quad y\mapsto [d_y] .
\end{equation}
It is continuous as a consequence of the triangle inequality. This map is a topological embedding provided that $Y$ is a geodesic  space; 
 see \cite[Ch. II.1]{Ballmann}, where this is proven for symmetric metrics, but the same proof goes through for nonsymmetric metrics as well. 
We assume from now on that the space $Y$ is {\em geodesic}. 

We identify $Y$ with its image in $\ol{\mathcal C}(Y)$ 
and call the closure $\ol Y$ the {\em horoclosure} of $Y$,
and $\geo Y:=\ol Y-Y$ the {\em horoboundary} or {\em boundary at infinity},
i.e. we have the decomposition 
$$ \ol Y = Y \sqcup\geo Y .$$
We note that the horoclosure $\ol{Y}$ is Hausdorff and 1st countable since the space $\ol{\mathcal C}(Y)$ is. 

The functions representing points in $\geo Y$ are called {\em horofunctions}.
We write 
$$ y_n \to [h] $$
for a divergent sequence of points $y_n\to\infty$ in $Y$ 
which converges to a point $[h]\in\geo Y$ represented by a horofunction $h$,
i.e. $d_{y_n}\to h$ modulo additive constants, 
and say that $(y_n)$ {\em converges at infinity}. Each horofunction is $1$-Lipschitz with respect to the symmetrized metric. 

If the metric space $(Y,d^{sym})$ is proper (which will be the case in this paper since we are interested in symmetric spaces), then the Arzel\`a-Ascoli theorem implies that the closure $\ol Y$ and the boundary $\geo Y$ at infinity are {\em compact}. 
In this case, $\ol Y$ is called the {\em horofunction compactification} of $Y$. 

Suppose that 
$$G\acts Y $$ 
is a $d$-isometric group action.
Then the embedding (\ref{eq:embfct}) is equivariant with respect to the induced action on functions
by $g\cdot f=f\circ g^{-1}$. 
For every 
$L>0$,
the subspace of $L$-Lipschitz functions 
$\ol{ Lip}_L(Y,d^{sym})\subset\ol{\mathcal C}(Y)$ 
is preserved by the action 
and contains, for  $L\geq1$, the horoclosure $\ol Y$.
We equip $G$ with the topology of uniform convergence on bounded subsets,
using the symmetrized metric $d^{sym}$ for both. 
Then the action $G\acts\ol{ Lip}_L(Y,d^{sym})$ is continuous.  
In particular,
the action 
$$G\acts\ol Y$$ 
is continuous. 
We will use this fact in the situation when $G$ is the isometry group of a Riemannian symmetric space of noncompact type.
In this case the topology of uniform convergence on bounded subsets 
coincides with the Lie group topology. 

An {\em oriented geodesic} in $(Y,d)$ is a ``forward" isometric embedding $c:I\to Y$,
i.e. for any parameters $t_1\leq t_2$ in $I$ it holds that 
$$ d(c(t_1),c(t_2)) = t_2-t_1 .$$
In particular, $c$ is continuous with respect to the symmetrized metric $d^{sym}$.
The metric space $(Y,d)$ is called a {\em geodesic} space,
if any pair of points $(y,y')$ can be connected by an oriented geodesic from $y$ to $y'$. 

If $(Y,d)$ is a geodesic space, 
then the horofunctions arising as limits of sequences along geodesic rays
are called {\em Busemann functions}, 
and their sublevel and level sets are called {\em horoballs} and {\em horospheres}.
We will denote by $Hb_b$ a horoball for the Busemann function $b$,
and more specifically, by $Hb_{b,y}$ the horoball of $b$ which contains the point $y$ in its boundary horosphere. 

In the situations studied in this paper,
all horofunctions will turn out to be Busemann functions,
cf. section~\ref{sec:idpts}.

Suppose that $Z\subset Y$ is a closed convex subset.
Then $Z$ is a geodesic space with respect to the induced metric,
and proper if $Y$ is proper.
There is a natural map 
\begin{equation}
\label{eq:extrintr}
\ol Z^Y\to\ol Z
\end{equation}
from the extrinsic closure $\ol Z^Y\subset\ol Y$ of $Z$ in $Y$
to the intrinsic horoclosure $\ol Z$ of $Z$.
It extends $\id_Z$, and at infinity is the map $\D^YZ\to\geo Z$ 
from the boundary of $Z$ in $\ol Y$ into the horoboundary of $Z$
given by the restriction of horofunctions to $Z$. 
If the latter map is injective, 
then (\ref{eq:extrintr}) is a homeomorphism
and there is a natural {\em embedding} of horoclosures 
$$ \ol Z \to\ol Y $$
given by ``unique extension of horofunctions''.

\subsection{Some notions of coarse geometry} \label{sec:coarse_geometry}

\begin{definition}
A correspondence $f: (X,d)\to (X',d')$ between metric spaces is {\em coarse Lipschitz} if there exist constants $L, A$ such that 
for all $x, y\in X$ and $x'\in f(x), y'\in f(y)$, we have
$$
d'(x',y')\le L d(x,y) + A. 
$$ 
\end{definition}

Note that if $(X,d)$ is a geodesic metric space, then in order to show that $f$ is coarse Lipschitz it suffices to verify that there exists a constant $C$ such that 
$$
d'(x', y')\le C
$$
for all $x, y\in X$ with $d(x,y)\le 1$ and all $x'\in f(x), y'\in f(y)$. 

\medskip
Two correspondences $f_1, f_2: (X,d)\to (X',d')$ are said to be {\em within distance $\le D$ from each other},  $\dist(f,g)\le D$, if 
for all $x\in X, y_i\in f_i(x)$, we have
$$
d'(y_1, y_2)\le D. 
$$
Two correspondences $f_1, f_2$ are said to be within finite distance from each other if $\dist(f_1, f_2)\le D$ for some $D$. 

A correspondence $(X,d)\to (X,d)$ is said to have {\em bounded displacement} if it is within finite distance from the identity map.

\begin{definition}
A coarse Lipschitz correspondence $f: (X,d)\to (X',d')$ is said to have a {\em  coarse left inverse} if there exists a coarse Lipschitz correspondence $g: X'\to X$ such that the composition $g\circ f$ has bounded displacement. 
\end{definition}

By applying the Axiom of Choice, we can always replace a coarse Lipschitz correspondence $f: (X,d)\to (X',d')$ 
with a coarse Lipschitz map $f': (X,d)\to (X',d')$ within bounded distance from $f$. With this in mind, if a coarse 
Lipschitz correspondence $f: (X,d)\to (X',d')$ admits a coarse left inverse, then $f$ is within bounded distance from a quasiisometric embedding $f': (X,d)\to (X',d')$. However, the converse is in general false, even in the setting of maps between 
finitely-generated groups equipped with word metrics.

We now specialize these concepts to the context of group homomorphisms. We note that each continuous homomorphism of groups with left-invariant proper metrics is always coarse Lipschitz. 
Suppose in the remainder of this section that $\Ga$ is a finitely generated group 
and $G$ is a connected Lie group equipped with a left invariant metric.

\begin{definition}
We say that for a homomorphism $\rho: \Ga\to G$,
a correspondence $r: G\to \Ga$ is a {\em coarse retraction} if $r$ is a coarse left inverse to $\rho$. 
A subgroup $\Ga< G$ is a {\em coarse retract} if the inclusion map $\Ga\embed G$ admits a coarse retraction.

Similarly, we say that a  homomorphism $\rho: \Ga\to G$ admits a {\em coarse equivariant retraction} if there exists a coarse Lipschitz retraction $r: G\to \Ga$ such that 
$$
r(h g)= r(h) r(g), \quad \forall h\in \rho(\Ga).  
$$
Accordingly, a subgroup $\Ga< G$ is a {\em coarse equivariant retract} if the inclusion homomorphism 
$\Ga\embed G$ admits a coarse equivariant retraction. 
\end{definition}

More generally, given an isometric action of $\rho: \Ga\acts X$ on a metric space $X$, 
we say that a coarse retraction $r: X\to \Ga$ is a {\em coarse equivariant retraction} if 
$$
r(\ga x)= \ga r(x), \quad \forall \ga\in \Ga, \quad x\in X.  
$$
In the case when $X=G/K$ is the symmetric space associated with a connected semisimple Lie group $G$, a homomorphism $\Ga\to G$ admits a coarse equivariant retraction iff the isometric action of $\Ga$ on $X$ defined via $\rho$ admits a coarse equivariant retraction. 
Similarly, a subgroup $\Ga< G$ is a coarse retract iff the orbit map $\Ga \to \Ga x\subset X$ admits a coarse left-inverse.

\section{Symmetric spaces}

\subsection{Basics}

In the paper we assume that the reader is familiar with basics of  symmetric spaces of noncompact type (denoted by $X$ throughout the paper), their isometry groups, visual boundaries and Tits boundaries. We refer the reader to \cite{Eberlein, qirigid, KLM} for the required background. 

In what follows, $G$ will be a connected semisimple Lie group, $K < G$ its maximal compact subgroup, the stabilizer of a base point in $X$ which will be denoted $o$ or $p$. Then $X\cong G/K$. We let $B< G$ denote a Borel subgroup. All maximal flats in $X$ are isometric to a {\em model} flat $F_{mod}$, which is isometric to a euclidean space ${\mathbb E}^n$, where $n$ is the rank of $X$. The model flat comes equipped with a (finite) Weyl group, denoted $W$. This group fixes the origin $0\in F_{mod}$, 
viewing $\Fmod$ as a vector space. We will use the notation $\amod$ for the visual boundary of $F_{mod}$; we will identify $\amod$ with the unit sphere 
in $\Fmod$ equipped with the angular metric. The sphere $\amod$ is the {\em model spherical apartment} for the group $W$. A fundamental domain for the action of $W$ on $F_{mod}$ is a certain convex cone $\Delta=\Delta_{mod}\subset F_{mod}$, the {\em model euclidean Weyl chamber} of $W$; its visual boundary is the model spherical Weyl chamber $\simod$, which is a spherical simplex in 
$\amod$. We let $\iota: \simod\to\simod$ denote the {\em opposition involution}, also known as the {\em standard involution}, of $\simod$; it equals $-w_0$, where $w_0\in W$ is the element sending $\simod$ to the opposite chamber in the model apartment $\amod$. We let $R\subset F_{mod}^*$ 
denote the root system of $X$, $\al_1,...,\al_n$ will denote simple roots with respect to $\Delta$; 
$$
\Delta=\{x\in F_{mod}: \al_i(x)\ge 0, i=1,...,n\}. 
$$

For a face $\taumod$ of $\simod$, we will frequently use the notation  
$+\taumod=\taumod$ and $-\taumod=\iota\taumod$. 

The space $X$ has the $\Delta$-valued ``distance function'' $d_\Delta$, which is the complete $G$-congruence invariant of pairs of points in $X$, see \cite{KLM}.  

$\ol X=X\sqcup\geo X$ will denote the {\em visual compactification} of $X$ with respect to its Riemannian metric, equipped with the visual topology, and $\tits X$ the Tits boundary of $X$, which is the visual boundary together with the Tits metric $\tangle$.
The Tits boundary carries a natural structure as a piecewise spherical simplicial complex. 
For a simplex $\tau$ in $\tits X$ we will use the notation $\interior(\tau)$ for the open simplex in $\tits X$
which is the complement in $\tau$ to the union of its proper faces. 

 We will denote by $xy$ the oriented geodesic segment in $X$ connecting a point $x$ to a point $y$; similarly, $x\xi$ will denote the geodesic ray from $x\in X$ asymptotic to the point $\xi\in \geo X$.

We will always use the notation $\tau, \hat\tau$ to indicate that the simplices $\tau, \hat\tau$ 
in $\tits X$ are opposite (antipodal), i.e., are swapped by a Cartan involution of $X$. Each simplex, of course, has a continuum of antipodal simplices. Simplices $\tau, \hat\tau$ are called {\em $x$-opposite} if the Cartan involution $s_x$ fixing $x$ sends $\tau$ to $\hat\tau$. Similarly, points $\xi, \hat\xi\in \geo X$ are $x$-opposite if $s_x$ swaps $\xi, \hat\xi$.

$\theta: \tits X\to \simod$ will denote the {\em type map}, i.e. the canonical projection of the Tits building to the model chamber. For distinct points $x, y\in X$ we let $\theta(xy)\in\simod$ denote the {\em type of the direction} of the oriented segment $xy$, i.e. the unit vector 
in the direction of the vector $d_\Delta(x,y)$. 

For each face $\taumod$ of $\simod$ one defines the {\em
flag manifold} $\Flagt$, which is the set of all simplices of type $\taumod$
in $\tits X$. Equipped with the {\em visual topology}, $\Flagt$ is a
homogeneous manifold homeomorphic to $G/P$, where $P$ is
a parabolic subgroup of $G$ stabilizing a face of type $\taumod$. The full flag manifold
$G/B= \Flag(\simod)$ is naturally identified with the Furstenberg boundary
$\DF X$ of $X$.  

For a simplex $\hat\tau\in\Flagit$, 
we let $C(\hat\tau)\subset \Flagt$  denote the subset consisting of simplices antipodal to $\hat\tau$. This subset is open and dense in $\Flagt$ and is called an {\em open Schubert stratum (or cell)} in $\Flagt$. 

For a point $x\in X$, $\Si_x X$ denotes the {\em space of directions} at $x$, i.e., the unit sphere in the tangent space $T_x X$. Similarly, for a spherical building ${\mathrm B}$ or a subcomplex ${\mathrm C}\subset {\mathrm B}$, and a point $\xi\in {\mathrm C}$, we let $\Si_\xi{\mathrm C}$ denote the space of directions of ${\mathrm C}$ at $\xi$. 

For a subset $Y\subset X$ we let $\geo Y$ denote the {\em visual boundary} of $Y$, i.e., 
its accumulation set in the visual boundary of $X$. A set $Y\subset X$ is said to be {\em asymptotic} to a subset $Z\subset \geo X$ if $Z\subset \geo Y$. 

For a subset $Z\subset \geo X$ we let $V(x, Z)\subset X$ denote the 
union of geodesic rays $x\zeta$ for all $\zeta\in Z$.
In the special case when $Z=\tau$ is a simplex\footnote{This means 
a simplex with respect to the spherical Tits building structure on $\geo X$.}
in $\geo X$,
then $V(x, \tau)$ is the {\em Weyl sector in $X$ with tip $x$ and base $\tau$}. 
A Weyl sector whose base is a chamber in $\geo X$ is a {\em (euclidean) Weyl chamber} in $X$. 

Two Weyl sectors $V(x_1,\tau)$ and $V(x_2,\tau)$ are {\em strongly asymptotic} 
if for any $\eps>0$ there exist points $y_i\in V(x_i,\tau)$ such that the subsectors $V(y_1,\tau)$ and $V(y_2,\tau)$ are $\eps$-Hausdorff close.

A sequence $x_i\in V(x,\tau)$ (where $\tau$ has the type $\taumod$) is {\em $\taumod$-regular} if it diverges from the boundary of $V(x,\tau)$, i.e., from the subsectors $V(x,\tau')$ for all proper faces $\tau'$ of $\tau$. 

For distinct points $x, y\in X$ and $\xi\in \geo X$ we let $\angle_x(y,\xi)$
denote the angle between the geodesic segment $xy$ and the geodesic ray 
$x\xi$ at the point $x\in X$.

$b_\eta$ will denote the Busemann function (defined with respect to the usual Riemannian  metric on $X$) associated with a point $\eta$ in the visual boundary of $X$. The gradient $\nabla b_\eta(x), x\in X,$ is the unit vector  tangent to the geodesic ray $x\eta$  and pointing away from $\eta$. 

$d$ will denote the standard distance function on $X$,  $Hb_{\eta}$ a closed horoball in $X$, which is a sublevel set $\{b_\eta\le t\}$ for the Riemannian Busemann function $b_\eta$. 

For a chamber $\si\subset\geo X$ we let $H_{\si}$ 
denote the associated {\em horocyclic subgroup}, 
the unipotent radical of the Borel subgroup $B_{\si}$ of $G$ stabilizing $\si$. Similarly, for a simplex $\tau$ in $\tits X$ we  we let $H_{\tau}$ 
denote the associated horocyclic subgroup, 
the unipotent radical in the parabolic subgroup $P_{\tau}$ of $G$ stabilizing $\tau$.

Elements of $H_{\tau}$ preserve the strong asymptote classes of geodesic rays $x\xi$, $\xi\in \interior(\tau)$ and hence the strong asymptote classes of sectors $V(x,\tau)$. Furthermore, 
$H_{\tau}$ acts transitively on the set of sectors $V(x',\tau)$ strongly asymptotic to the given sector $V(x,\tau)$.

\subsection{Parallel sets, stars and cones} 

\subsubsection{Parallel sets}
\label{sec:parset}

Let $s\subset\tits X$ be an isometrically embedded (simplicial) unit sphere. 
We denote by $P(s)\subset X$ the {\em parallel set} associated to $s$, which 
can be defined as the 
union of maximal flats $F\subset X$ asymptotic to $s$, 
$s\subset\geo F$. 
Alternatively, 
one can define it as the 
union of flats $f\subset X$ with ideal boundary $\geo f=s$. 

The parallel set is a totally geodesic subspace 
which splits metrically as the product 
\begin{equation}
\label{eq:parsplit}
P(s)\cong f\times CS(s),
\end{equation}
of any of these flats and a symmetric space $CS(s)$ 
called its {\em cross section}.  Accordingly, 
the ideal boundary of the parallel set is a metric suspension 
\begin{equation}
\label{eq:bdparsersusp}
\tits P(s)\cong\tits f\circ\tits CS(s) .
\end{equation}
It coincides with the subbuilding $\B(s)\subset\geo X$ 
consisting of the union of all apartments $a\subset\geo X$ containing $s$,
\begin{equation*}
\B(s)=\geo P(s).
\end{equation*}

It is immediate that parallel sets are nonpositively curved symmetric spaces.
However, they do not have noncompact type as their Euclidean de Rham factors 
are nontrivial. 
The factor $f$ in the splitting (\ref{eq:bdparsersusp})
of the parallel set is then the Euclidean de Rham factor
and the cross section $CS(s)$ has trivial euclidean de Rham factor,
i.e.\ it is a symmetric space of noncompact type.

For a pair of antipodal simplices $\tau_+, \tau_-\subset\geo X$
there exists a unique minimal singular sphere $s=s(\tau_-, \tau_+)\subset \geo X$ containing them.
We denote $P(\tau_-, \tau_+):=P(s(\tau_-, \tau_+))$; this parallel set is the union of (maximal) flats $F\subset X$ whose ideal boundaries contain $\tau_-\cup \tau_+$. In order to simplify the notation, we will denote $\B(s(\tau_-, \tau_+))$ simply by $\B(\tau_-, \tau_+)$.  Given two antipodal points 
$\zeta_+, \zeta_-\in \tits X$ we let $\B(\zeta_-, \zeta_+)$ denote the subbuilding $\B(\tau_-, \tau_+)$, 
where $\tau_\pm$ are the antipodal simplices satisfying $\zeta_\pm \in \interior(\tau_\pm)$. 

We will use the notation $CS(\tau_-, \tau_+, p)$ for the cross-section $CS(s)$ passing through the point $p\in P(s)$. Similarly, $f(\tau_-,\tau_+,p)$ will denote the flat in $P(\tau_-,\tau_+)$ which is parallel to the euclidean factor of $P(\tau_-, \tau_+)$ and contains $p$. 

We will use the notation $T(s)=T(\tau_-, \tau_+)$ for the group of transvections 
along the flat $f$: This group is the same for all flats parallel to $f$ and depends only on $s$.

\subsubsection{Stars}
\label{sec:star}

\begin{dfn}[Stars]
\label{dfn:star}
Let $\tau\subset\tits X$ be a simplex. 
We define the {\em star}\footnote{$\st(\tau)$ is also known as the {\em residue} of $\tau$.} 
$\st(\tau)$ of the open simplex $\interior(\tau)$ as the subcomplex of $\tits X$ 
consisting of all simplices intersecting the open simplex $\interior(\tau)$ nontrivially (i.e., containing $\tau$). In other words, $\st(\tau)$ is the smallest subcomplex of $\tits X$ containing all chambers $\si$ such that $\tau\subset \si$.

We define the {\em open star} $\ost(\tau)\subset\geo X$ 
as the union of all open simplices  whose closure intersects $\interior(\tau) $ nontrivially. For the model simplex $\tau_{mod}$, we will use the notation 
$\ost(\tau_{mod})$ to denote its open star in the simplicial complex consisting of faces of $\si_{mod}$. 

For a point $\xi\in \tits X$ we let $\st(\xi)\subset\geo X$ denote the star of the simplex spanned by 
$\xi$, i.e. the unique simplex $\tau$ such that $\xi\in \interior(\tau)$. 
\end{dfn}

\medskip 
Note that $\ost(\tau)$ is an open subset of the simplex $\si_{mod}$; it 
does not include any open faces of $\tau$ except for the interior of $\tau$. Furthermore, $\D\st(\tau)=\st(\tau)-\ost(\tau)$ is the union of all panels $\pi$ of type $\theta(\pi)\not\supset\tau_{mod}$ which are contained in a chamber with face $\tau$. 

\begin{lem}
\label{lem:starintapt}The star $\st(\bar\tau)$ of a simplex $\bar\tau\subset a_{mod}$ 
is a convex subset of $a_{mod}$. Furthermore,  $\st(\bar\tau)$ equals the intersection of the simplicial hemispheres $\bar h\subset a_{mod}$ 
such that  $\interior(\bar\tau)\subset\interior{\bar h}$. 
\end{lem}
\proof
If a hemisphere $\bar h$ contains a simplex $\bar\tau$, but does not contain it in its boundary, 
then all chambers containing this simplex as a face belong to the (closed) hemisphere. 
Vice versa, if a chamber $\bar\si$ does not contain $\bar\tau$ as a face, 
then there exists a wall which separates $\bar\si$ from $\bar\tau$. 
\qed

\medskip
Similarly, the star $\st(\tau)$ of a simplex $\tau\subset\tits X$ 
is a convex subset of $\tits X$. 
One can represent it as the intersection of all simplicial $\pihalf$-balls 
which contain $\interior(\tau)$ in their interior. One can represent $\st(\tau)$ also as the intersection of fewer balls:

\begin{lem}[Convexity of stars]
\label{lem:starconv}
(i) 
Let $\tau\subset\tits X$ be a simplex. 
Then $\st(\tau)$ equals the intersection of the simplicial $\pihalf$-balls 
whose interior contains $\interior(\tau)$. 

(ii)
For any simplex $\hat\tau$ opposite to $\tau$, 
the star $\st(\tau)$ equals the intersection 
of the subbuilding 
$\B(\tau,\hat\tau)=\geo P(\tau,\hat\tau)$ 
with all simplicial $\pihalf$-balls 
whose interior contains $\interior(\tau)$ 
and whose center lies in this subbuilding. 
\end{lem}
\proof
(i) If a simplicial $\pihalf$-ball contains a simplex $\tau$, but does not contain it in its boundary, 
then all chambers containing this simplex as a face belong to this ball. 
Vice versa, let $\si$ be a chamber which does not contain $\tau$ as a face. 
There exists an apartment $a\subset\tits X$ which contains $\si$ and $\tau$. 
As before in the proof of Lemma~\ref{lem:starintapt}, 
there exists a simplicial hemisphere $h\subset a$ containing $\tau$ but not $\si$. 
Then the simplicial $\pihalf$-ball with the same center as $h$ contains $\tau$ but not $\si$. 

(ii) Note first that $\st(\tau)\subset \B(\tau,\hat\tau)$. 
Then we argue as in part (i), observing that 
if $\si\subset \B(\tau,\hat\tau)$ then $a$ can be chosen inside $\B(\tau,\hat\tau)$. 
\qed

\subsubsection{Stars and ideal boundaries of cross sections}
\label{sec:idbdcrsct}

Let $\nu\subset\geo X$ be a simplex.
We say that two chambers $\si_1,\si_2\supset\nu$ are {\em $\nu$-antipodal}
if there exists a segment connecting interior points of $\si_1,\si_2$ 
and passing through an interior point of $\nu$.

The link $\Si_{\nu}\st(\nu)$ carries a natural structure 
as a topological spherical building,
and is naturally isomorphic as such to 
$\geo CS(\hat\nu,\nu)$ for any $\hat\nu$ opposite to $\nu$.
Chambers $\si\subset\st(\nu)$ correspond to chambers in $\Si_{\nu}\st(\nu)$,
and pairs of $\nu$-opposite chambers to pairs of opposite chambers. 
It follows that for every chamber $\si\supset\nu$,
the set of $\nu$-opposite chambers in $\st(\nu)$
is open and dense as a subset of $\stF(\nu)$.
Here, we denote by $\stF(\nu)\subset\DF X$ the subset of chambers containing $\nu$.

\subsubsection{Weyl cones}

Given a simplex $\tau$ in $\tits X$ and a point $x\in X$, the union $V(x, \st(\tau))$ of all rays $x\zeta, \zeta\in \st(\tau)$ is called the {\em Weyl cone} with the tip $x$ and the base $\st(\tau)$. 
Below we will prove that 
Weyl cones $V(x,\st(\tau))$ are convex. We begin with 

\begin{lemma}
For every $x\in P(\tau,\hat\tau)$,  
the Weyl cone $V(x,\st(\tau))$ is contained in the parallel set 
$P(\tau,\hat\tau)$.
\end{lemma}
\proof Consider a chamber $\si$ in $\tits X$ containing $\tau$.  The Weyl sector $V(x,\si)$ is contained in a (unique maximal) flat $F\subset X$. 
Since $\tau, \hat \tau$ are antipodal with respect to $x$, $\tau\cup \hat\tau\subset \geo F$. Therefore, $F\subset P(\tau,\hat\tau)$. \qed

\begin{prop}[Convexity of Weyl cones]
\label{prop:wconeconv}
Let $\hat\tau$ be the simplex opposite to $\tau$ with respect to $x$. 
Then the Weyl cone $V(x,\st(\tau))$ is the intersection of the parallel set $P(\tau,\hat\tau)$
with the horoballs which are centered at $\geo P(\tau,\hat\tau)$ and contain $V(x,\st(\tau))$. In particular, $V(x, \st(\tau))$ is a closed convex subset of $X$.  
\end{prop}
\proof
One inclusion is clear. 
We must prove that each point $y\in P(\tau,\hat\tau)- V(x,\st(\tau))$ 
is not contained in one of these horoballs. 
There exists a maximal flat $F\subset P(\tau,\hat\tau)$ containing $x$ and $y$. 
(Any two points in a parallel set lie in a common maximal flat.)
We extend the oriented segment $xy$ to a ray $x\eta$ inside $F$. 

As in the proof of Lemma~\ref{lem:starconv}, 
there exists $\zeta\in\geo F$ such that 
$\bar B(\zeta,\pihalf)$ contains $\st(\tau)$ 
but does not contain $\eta$. 
Then the horoball $Hb_{\zeta,x}$ intersects $F$ in a half-space 
which contains $x$ in its boundary hyperplane but does not contain $\eta$ in its ideal boundary. 
Therefore does not contain $y$.
By convexity, 
$V(x,\st(\tau))\subset Hb_{\zeta,x}$. 
\qed

\medskip
The following consequence will be important for us.
\begin{cor}[Nested cones]
\label{cor:nestcone}
If $x'\in V(x,\st(\tau))$, then 
$V(x',\st(\tau))\subset V(x,\st(\tau))$. 
\end{cor}

Let $xy\subset X$ be an oriented $\tau_{mod}$-regular geodesic segment. 
Then we define the simplex $\tau=\tau(xy)\subset\geo X$ as follows: 
Forward extend the segment $xy$ to the geodesic ray $x\xi$, 
and let $\tau$ be the unique face of type $\tau_{mod}$ of $\tits X$ such that $\xi\in \st(\tau)$.

\begin{dfn}[Diamond]
We define the {\em $\taumod$-diamond} of a $\taumod$-regular segment $x_-x_+$ as 
\begin{equation*}
\diamot(x_-,x_+)=
V(x_-,\st(\tau_+))\cap V(x_+,\st(\tau_-))
\subset P(\tau_-,\tau_+)
\end{equation*}
where $\tau_{\pm}=\tau(x_{\mp}x_{\pm})$. 
\end{dfn}

Thus, every diamond is a convex subset of $X$.

\subsubsection{Shadows at infinity and strong asymptoticity of Weyl cones}
\label{sec:shadw}

This material is taken from \cite[\S 4.1]{coco15}.

For a simplex 
$\tau_-\in\Flagit$
and a point $x\in X$,
we consider the function 
\begin{equation}
\label{eq:distfrpar}
\tau\mapsto d(x,P(\tau_-,\tau))
\end{equation}
on the open Schubert stratum $C(\tau_-)\subset\Flagt$.
We denote by $\tau_+\in C(\tau_-)$ the simplex $x$-opposite to $\tau_-$.
\begin{lem}
\label{lem:contpr}
The function (\ref{eq:distfrpar}) is continuous and proper.
\end{lem}
\proof
This follows from the fact that 
$C(\tau_-)$ and $X$ are homogeneous spaces for the parabolic subgroup $P_{\tau_-}$.
Indeed, 
continuity follows from the continuity of the function 
$$g\mapsto d(x,P(\tau_-,g\tau_+))=d(g^{-1}x,P(\tau_-,\tau_+))$$
on $P_{\tau_-}$ which factors through the orbit map $P_{\tau_-}\to C(\tau_-),g\mapsto g\tau_+$.

Regarding properness,
note that a simplex $\tau\in C(\tau_-)$ is determined by any point $y$ contained in the parallel set $P(\tau_-,\tau)$,
namely as the simplex $y$-opposite to $\tau_-$.
Thus,
if $P(\tau_-,\tau)\cap B(x,R)\neq\emptyset$ for some fixed $R>0$,
then there exists $g\in P_{\tau_-}$ such that $\tau=g\tau_+$ and $d(x,gx)<R$. 
In particular, $g$ lies in a compact subset.
This implies properness.
\qed

\medskip
Moreover, 
the function (\ref{eq:distfrpar}) has a unique minimum zero in $\tau_+$. 

We define the following open subsets of $C(\tau_-)$ 
which can be regarded as {\em shadows} of balls in $X$ with respect to $\tau_-$. 
For $x\in X$ and $r>0$, we put
\begin{equation*}
U_{\tau_-,x,r}:=\{\tau\in C(\tau_-) | d(x,P(\tau_-,\tau))<r\} .
\end{equation*}
The next fact expresses the uniform strong asymptoticity
of asymptotic Weyl cones. 
\begin{lem}
\label{lem:expconvsect}
For $r,R>0$ exists $d=d(r,R)>0$ such that:

If $y\in V(x,\st(\tau_-))$ with $d(y,\D V(x,\st(\tau_-)))\geq d(r,R)$, 
then 
$U_{\tau_-,x,R}\subset U_{\tau_-,y,r}$.
\end{lem}
\proof
If $U_{\tau_-,x,R}\not\subset U_{\tau_-,y,r}$ 
then there exists $x'\in B(x,R)$ such that $d(y,V(x',\st(\tau_-)))\geq r$. 
Thus, if the assertion is wrong, 
there exist a sequence $x_n\to x_{\infty}$ in $\ol B(x,R)$ 
and a sequence $y_n\to\infty$ in $V(x,\st(\tau_-))$ 
such that $d(y_n,\D V(x,\st(\tau_-)))\to+\infty$
and $d(y_n,V(x_n,\st(\tau_-)))\geq r$.  

Let $\rho:[0,+\infty)\to V(x,\tau_-)$ be a geodesic ray with initial point $x$ and asymptotic to an interior point of $\tau_-$.
Then the sequence $(y_n)$ eventually enters every Weyl cone $V(\rho(t),\st(\tau_-))$. 
Since the distance function $d(\cdot,V(x_n,\st(\tau_-)))$ is convex and bounded, and hence non-increasing 
along rays asymptotic to $\st(\tau_-)$, 
we have that 
\begin{equation*}
R\geq d(x,V(x_n,\st(\tau_-))) 
\geq d(\rho(t),V(x_n,\st(\tau_-)))\geq d(y_n,V(x_n,\st(\tau_-)))\geq r
\end{equation*}
for $n\geq n(t)$. 
It follows that 
\begin{equation*}
R\geq d(\rho(t),V(x_{\infty},\st(\tau_-)))\geq r
\end{equation*}
for all $t$. 
However, the ray $\rho$ is strongly asymptotic to $V(x_{\infty},\st(\tau_-))$, a contradiction. 
\qed

\subsubsection{Some spherical building facts}
\label{sec:sphbl}

We discuss some facts from spherical building geometry.
In this paper,
they are applied to the visual boundary $\geo X$ 
equipped with its structure as a thick spherical Tits building.

First recall the following lemma, cf.\ the first part of \cite[Lemma 3.10.2]{qirigid}.\footnote{The statement of the second part contains a typo: $H$ should be replaced by $\ol B(\eta,\pihalf)\cap A$.}

\begin{lemma}\label{lem:antipode}
In a spherical building $\B$ every point $\xi\in \B$ has an antipode in every apartment $a\subset \B$,
and hence for every simplex $\tau\subset B$ there is an opposite simplex $\hat\tau\subset a$.
\end{lemma}

We need the more precise statement
that a point has {\em several} antipodes in an apartment
unless it lies in the apartment itself:

\begin{lem}[cf.\ {\cite[Sublemma 5.20]{morse}}]
\label{lem:oneantip}
Let $\xi$ be a point in a spherical building $\B$ and let $a\subset \B$ be an apartment. 
If $\xi$ has only one antipode in $a$, 
then $\xi\in a$.
\end{lem}
\proof
Suppose that $\xi\not\in a$ and let $\hat\xi\in a$ be an antipode of $\xi$. 
We choose a``generic'' segment $\xi\hat\xi$ of length $\pi$ tangent to $a$ at $\hat\xi$ as follows.
The suspension $\B(\xi,\hat\xi)\subset \B$ contains an apartment $a'$ with the same unit tangent sphere at $\hat\xi$, $\Si_{\hat\xi}a'=\Si_{\hat\xi}a$. 
Inside $a'$ there exists a segment $\xi\hat\xi$ whose interior does not meet simplices of codimension $\geq2$. 
Hence $\hat\xi\xi$ leaves $a$ at an interior point $\eta\neq\xi,\hat\xi$ of a panel $\pi\subset a$,
i.e.\ $a\cap\xi\hat\xi=\eta\hat\xi$ and $\pi\cap\xi\hat\xi=\eta$,
and $\eta\xi$ initially lies in a chamber adjacent to $\pi$ but not contained in $a$. 
Let $s\subset a$ be the wall (codimension one singular sphere) containing $\pi$. 
By reflecting $\hat\xi$ at $s$, one obtains a second antipode for $\xi$ in $a$, 
contradiction. 
\qed

\medskip
We will also need the following fact:
\begin{lem}
\label{lem:intapts}
Suppose that the spherical building $\B$ is thick.
Then for any simplex $\tau\subset \B$, 
the intersection of all apartments containing $\tau$ equals $\tau$.
\end{lem}
\proof
Suppose first that $\tau$ is a chamber. 
Let $\xi\in \B-\tau$.
For a generic point $\eta\in\interior(\tau)$, the segment $\eta\xi$ leaves $\tau$ 
through an interior point of a panel.
By thickness, there exists an apartment $a\subset \B$ such that 
$a\cap\eta\xi=\tau\cap\eta\xi$.
Then $\xi\not\in a$.
This shows the assertion in the case when $\tau$ is a chamber. 
If $\tau$ is an arbitrary simplex, it follows that the intersection of all apartments containing $\tau$ 
is contained in the intersection of all chambers containing $\tau$,
which equals $\tau$.
\qed

\section{Regularity and contraction}
\label{app:regcontr}

This material is taken from \cite{morse}.
We include it for completeness.

In this section, we define a class of discrete subgroups of semisimple Lie groups 
by a certain 
asymptotic {\em regularity} condition
which in rank one just amounts to discreteness, but in higher rank is strictly stronger. 
The condition will be formulated in two equivalent ways.
First dynamically in terms of the action on a flag manifold, 
then geometrically in terms of the orbits in the symmetric space.
This class of subgroups contains the class of Anosov subgroups.

\subsection{Contraction}

Consider the action $$G\acts\Flagt$$
on the flag manifold of type $\taumod$.
We recall that for a simplex $\tau_-$ of type $-\taumod:=\iota\taumod$ 
we denote by $C(\tau_-)\subset\Flagt$ the open dense $P_{\tau_-}$-orbit.

We introduce the following dynamical conditions for sequences and subgroups in $G$:
\begin{dfn}[Contracting sequence, cf.\ {\cite[Def.\ 6.1]{coco15}}]
\label{def:contracting_sequence}
A sequence $(g_n)$ in $G$ is {\em $\taumod$-con\-trac\-ting} 
if there exist simplices $\tau_{\pm}\in\Flagpmt$ such that 
\begin{equation}
\label{eq:contrtau}
g_n|_{C(\tau_-)}\to\tau_+
\end{equation} 
uniformly on compacts as $n\to+\infty$.
\end{dfn}

\begin{dfn}[Convergence type dynamics, cf.\ {\cite[Def.\ 6.10]{morse}} and {\cite[Def.\ 2.2]{anosov}}]
\label{def:conv}
A subgroup $\Ga<G$ is called a {\em $\taumod$-convergence subgroup}
if every sequence $(\ga_n)$
of pairwise distinct elements in $\Ga$ contains a $\taumod$-contracting subsequence.
\end{dfn}

Note that $\taumod$-contracting sequences diverge to infinity
and therefore $\taumod$-convergence subgroups are necessarily discrete.

The contraction property exhibits a symmetry:
\begin{lem}[Symmetry]
\label{lem:contrsym}
Property (\ref{eq:contrtau}) is equivalent to the dual property that 
\begin{equation}
\label{eq:contrtaudual}
g^{-1}_n|_{C(\tau_+)}\to\tau_-
\end{equation} 
uniformly on compacts as $n\to+\infty$.
\end{lem}
\proof
We give a simplified version of \cite[Proposition 6.5]{coco15} and its proof.

Suppose that (\ref{eq:contrtau}) holds but (\ref{eq:contrtaudual}) fails.
Equivalently, after extraction
there exists a sequence 
$\xi_n\to\xi\neq\tau_-$ in $\Flagit$ 
such that $g_n\xi_n\to\xi'\in C(\tau_+)$.
Since $\xi\neq\tau_-$,
there exists $\hat\tau_-\in C(\tau_-)$ not opposite to $\xi$.
(For instance, take an apartment in $\geo X$ containing $\tau_-$ and $\xi$, 
and let $\hat\tau_-$ be the simplex opposite to $\tau_-$ in this apartment.)
Hence there is a sequence $\tau_n\to\hat\tau_-$ in $\Flagt$ such that $\tau_n$ is not opposite to $\xi_n$ for all $n$.
It can be obtained e.g.\ by taking a sequence $h_n\to e$ in $G$ such that $\xi_n=h_n\xi$
and putting $\tau_n=h_n\hat\tau_-$.
Since $\hat\tau_-\in C(\tau_-)$,
condition (\ref{eq:contrtau}) implies that $g_n\tau_n\to\tau_+$.
It follows that $\tau_+$ is not opposite to $\xi'$, 
because $g_n\tau_n$ is not opposite to $g_n\xi_n$ 
and being opposite is an open condition.
This contradicts that $\xi'\in C(\tau_+)$.
Therefore,
condition (\ref{eq:contrtau}) implies (\ref{eq:contrtaudual}). 
The converse implication follows by replacing the sequence $(g_n)$ with
$(g_n^{-1})$. 
\qed

\begin{lem}[Uniqueness]
\label{lem:contruniq}
The simplices $\tau_{\pm}$ in (\ref{eq:contrtau}) are uniquely determined. 
\end{lem}
\proof
Suppose that besides (\ref{eq:contrtau}) we also have that $g_n|_{C(\tau'_-)}\to\tau'_+$
with other simplices $\tau'_{\pm}$.
Since the subsets $C(\tau_-)$ and $C(\tau'_-)$ are open dense in $\Flagit$,
it follows that their intersection is nonempty 
and hence $\tau'_+=\tau_+$.
Using the equivalent dual conditions (\ref{eq:contrtaudual})
we similarly obtain that $\tau'_-=\tau_-$.
\qed

\subsection{Regularity}
\label{sec:reg}

The second set of asymptotic properties concerns the geometry of orbits in $X$.

We first consider sequences in the model euclidean Weyl chamber $\De=\De_{mod}$.
\begin{dfn}[cf.\ {\cite[Def.\ 4.4]{coco15}}]
\label{def:pureg}
A sequence $(\de_n)$ in $\De$ is 

(i) {\em $\taumod$-regular} if it drifts away from $V(0,\D\st(\taumod))\subset\D\De$,
$$ d(\de_n,V(0,\D\st(\taumod))) \to+\infty .$$

(ii) {\em $\taumod$-pure} if it is contained in a tubular neighborhood of the sector $V(0,\taumod)$
and drifts away from its boundary $\D V(0,\taumod)=V(0,\D\taumod)$,
$$ d(\de_n,V(0,\D\taumod)) \to+\infty .$$
\end{dfn}
Note that $(\de_n)$ is $\taumod$-regular/pure iff 
$(\iota\de_n)$ is $\iota\taumod$-regular/pure.

We extend these notions to sequences in $X$ and $G$:
\begin{dfn}[Regular and pure, cf.\ {\cite[Def.\ 4.4]{coco15}}]
(i) A sequence $(x_n)$ in $X$ is {\em $\taumod$-regular}, respectively, {\em $\taumod$-pure}
if for some (any) base point $o\in X$ the sequence of $\De$-distances $d_{\De}(o,x_n)$ 
has this property.

(ii) A sequence $(g_n)$ in $G$ is {\em $\taumod$-regular}, respectively, {\em $\taumod$-pure}
if for some (any) point $x\in X$ the orbit sequence $(g_nx)$ in $X$ has this property.

(iii) A subgroup $\Ga<G$ is {\em $\taumod$-regular}
if all sequences of pairwise distinct elements in $\Ga$ have this property.
\end{dfn}

That these properties are independent of the base point and stable under bounded perturbation of the sequences,
is due to the triangle inequality
$|d_{\De}(x,y)-d_{\De}(x',y')|\leq d(x,x')+d(y,y')$.

A sequence $(g_n)$ is $\taumod$-regular/pure iff 
the inverse sequence $(g_n^{-1})$ is $\iota\taumod$-regular/pure.
In particular,
a subgroup $\Ga<G$ is $\taumod$-regular iff it is $\iota\taumod$-regular.
It then follows that $\Ga$ is also $\taumod'$-regular
where $\taumod'$ is the span of $\taumod$ and $\iota\taumod$. 

The face type of a pure sequence is uniquely determined.
Moreover,
a $\taumod$-regular sequence is $\taumod'$-regular for every face type $\taumod'\subseteq\taumod$,
because $\ost(\taumod')\supseteq\ost(\taumod)$.

Note that $\taumod$-regular subgroups are in particular discrete.
If $\rank(X)=1$, then discreteness is equivalent to ($\simod$-)regularity.
In higher rank, {\em regularity} can be considered as a {\em strengthening of discreteness}:
A discrete subgroup $\Ga<G$ may not be $\taumod$-regular for any face type $\taumod$;
this can happen e.g.\ for free abelian subgroups of transvections of rank $\geq2$.

We observe furthermore:
\begin{lem}
\label{lem:obspureg}
(i) $\taumod$-pure sequences are $\taumod$-regular.

(ii) Every sequence, which diverges to infinity, 
contains a $\taumod$-pure subsequence for some face type $\taumod\subseteq\simod$.
\end{lem}
\proof
Assertion (i) is a direct consequence of the definitions,
and (ii) follows by induction on face types.
\qed

\medskip
Note also that a sequence, which diverges to infinity, 
is $\taumod$-regular if and only if it contains $\numod$-pure subsequences only for face types $\numod\supseteq\taumod$. 
(We will not use this fact.)

\medskip
The lemma implies in particular,
that every sequence $\ga_n\to\infty$ in a discrete subgroup $\Ga<G$ contains a subsequence 
which is $\taumod$-regular, even $\taumod$-pure, for some $\taumod$.

\begin{rem}
\label{rem:regfins}
Pureness and regularity have natural interpretations in terms of Finsler geometry, 
see Proposition~\ref{prop:puregstr}.
\end{rem}

\subsection{Contraction implies regularity}
\label{sec:contrimplrg}

In order to relate the contraction and regularity properties,
it is useful to consider the $G$-action on flats in $X$.

We recall that ${\mathcal F}_{\taumod}$ denotes the space of flats $f\subset X$ of type $\taumod$.
Two flats $f_{\pm}\in{\mathcal F}_{\taumod}$ are {\em dynamically related} 
with respect to a sequence $(g_n)$ in $G$,
$$f_-\stackrel{(g_n)}{\sim}f_+ ,$$
if there exists a sequence of flats $f_n\to f_-$ such that $g_nf_n\to f_+$.
Moreover, the action of $(g_n)$ on ${\mathcal F}_{\taumod}$ is {\em proper} 
iff there are no dynamical relations with respect to subsequences,
cf.\ section~\ref{sec:discdynrl}.

Dynamical relations between singular flats yield dynamical relations between maximal ones:
\begin{lem}
\label{lem:dynrsingmx}
If $f_{\pm}\in{\mathcal F}_{\taumod}$ are flats such that $f_-\stackrel{(g_n)}{\sim}f_+$,
then for every maximal flat $F_+\supset f_+$ there exist a maximal flat $F_-\supset f_-$
and a subsequence $(g_{n_k})$ 
such that $F_-\stackrel{(g_{n_k})}{\sim}F_+$.
\end{lem}
\proof
Let $f_n\to f_-$ be a sequence in ${\mathcal F}_{\taumod}$
such that $g_nf_n\to f_+$.
Then there exists a sequence of maximal flats $F_n\supset f_n$ such that $g_nF_n\to F_+$.
The sequence $(F_n)$ is bounded,\footnote{Equivalently, 
all flats in the sequence intersect a fixed bounded subset of $X$.}
because the sequence $(f_n)$ is,
and hence subconverges to a maximal flat $F_-$.
\qed

\medskip
For pure sequences there are dynamical relations between singular flats of the corresponding type 
with respect to suitable subsequences:
\begin{lem}
\label{lem:improp}
If $(g_n)$ is $\taumod$-pure,
then the action of $(g_n)$ on ${\mathcal F}_{\taumod}$ is not proper.

More precisely,
there exist simplices $\tau_{\pm}\in\Flagt$ such that 
for every flat $f_+\in{\mathcal F}_{\taumod}$ asymptotic to $\tau_+$
there exist a flat $f_-\in{\mathcal F}_{\taumod}$ asymptotic to $\tau_-$
and a subsequence $(g_{n_k})$ such that 
$$f_-\stackrel{(g_{n_k})}{\sim}f_+.$$

\end{lem}
\proof
We prove the stronger statement. 
By pureness, there exists a sequence of simplices $(\tau_n)$ in $\Flagt$
such that for any point $x\in X$ we have 
$$ \sup_n d(g_nx,V(x,\tau_n))  <+\infty .$$
There exists a subsequence $(g_{n_k})$ such that 
$\tau_{n_k}\to\tau_+$ and $g_{n_k}^{-1}\tau_{n_k}\to\tau_-$.

Let $f_+\in{\mathcal F}_{\taumod}$ be asymptotic to $\tau_+$.
We choose $x\in f_+$ and consider the sequence of flats $f_{n_k}\in{\mathcal F}_{\taumod}$ 
through $x$ asymptotic to $\tau_{n_k}$.
Then $f_{n_k}\to f_+$.
The sequence of flats $(g_{n_k}^{-1}f_{n_k})$ is bounded.
Therefore, after further extraction, we obtain convergence $g_{n_k}^{-1}f_{n_k}\to f_-$.
\qed

\medskip
By a diagonal argument one can also show that the subsequences $(g_{n_k})$ 
in the two previous lemmas 
can be made independent of the flats $F_+$ respectively $f_+$.

For contracting sequences, 
the possible dynamical relations between maximal flats are restricted as follows:
\begin{lem}
\label{lem:contrdynrel}
Suppose that $(g_n)$ is {\em $\taumod$-contracting} with (\ref{eq:contrtau}),
and that 
$F_-\stackrel{(g_n)}{\sim}F_+$
for maximal flats $F_{\pm}\in{\mathcal F}$. 
Then $\tau_{\pm}\subset\geo F_{\pm}$.
\end{lem}
\proof
Compare \cite[\S 5.2.4]{morse}.

The visual boundary sphere of any maximal flat $F\subset X$
contains a simplex $\tau$ opposite $\tau_-$.\footnote{See Lemma \ref{lem:antipode}.} 
Suppose that $\tau_-\not\subset\geo F_-$.
Then $\geo F_-$ contains {\em several}, 
i.e.\ at least two different simplices $\hat\tau_-,\hat\tau'_-$ opposite $\tau_-$,
cf.\ Lemma~\ref{lem:oneantip}.

Let $F_n\to F_-$ be a sequence in ${\mathcal F}$ such that $g_nF_n\to F_+$.
There are sequences of simplices $\tau_n,\tau'_n\subset\geo F_n$ 
such that $\tau_n\to\hat\tau_-$ and $\tau'_n\to\hat\tau'_-$.
After extraction,
we also have convergence 
$g_n\tau_n\to\hat\tau_+$ and $g_n\tau'_n\to\hat\tau'_+$
to {\em different} limit simplices $\hat\tau_+,\hat\tau'_+\subset\geo F_+$. 
On the other hand,
in view of $\hat\tau_-,\hat\tau'_-\in C(\tau_-)$, 
the contraction property (\ref{eq:contrtau}) implies that 
$g_n\tau_n\to\tau_+$ and $g_n\tau'_n\to\tau_+$,
convergence to the same simplex, 
a contradiction.
Thus, $\tau_-\subset\geo F_-$.

Considering the inverse sequence $(g_n^{-1})$
yields that also $\tau_+\subset\geo F_+$,
cf.\ Lemma~\ref{lem:contrsym}.
\qed

\medskip
Combining the previous lemmas, we obtain:
\begin{lem}
\label{lem:contrpuresq}
If a sequence in $G$ is $\taumod$-contracting and $\numod$-pure, then $\taumod\subseteq\numod$.
\end{lem}
\proof
We denote the sequence by $(g_n)$.
According to Lemmas~\ref{lem:improp} and~\ref{lem:dynrsingmx},
by $\numod$-purity, there exist simplices $\nu_{\pm}\in\Flagn$
such that for every maximal flat $F_+$ asymptotic to $\nu_+$ there exist a maximal flat $F_-$ asymptotic to $\nu_-$
and a subsequence $(g_{n_k})$ such that 
$$F_-\stackrel{(g_{n_k})}{\sim}F_+$$
Lemma~\ref{lem:contrdynrel} implies that $\tau_+\subset\geo F_+$.
Varying $F_+$, it follows that $\tau_+\subseteq\nu_+$.\footnote{See  Lemma~\ref{lem:intapts}.}
\qed

\medskip
From these observations, we conclude:
\begin{prop}[Contracting implies regular]
\label{prop:contrimpreg}
If a sequence $(g_n)$ is $\taumod$-contracting, then it is $\taumod$-regular.
\end{prop}
\proof
Suppose that 
$(g_n)$ is not $\taumod$-regular.
Then after extraction $(g_n)$ is $\numod$-pure for a face type $\numod\not\supset\taumod$.
However, this contradicts the last lemma.
\qed

\subsection{Regularity implies contraction}

We now prove a converse to Proposition~\ref{prop:contrimpreg}.
Since contractivity involves a convergence condition,
we can expect regular sequences to be contracting only after extraction.

Consider a $\taumod$-regular sequence $(g_n)$ in $G$.
After fixing a point $x\in X$,
there exist simplices $\tau_n^{\pm}\in\Flagpmt$
(unique for large $n$)
such that 
\begin{equation}
\label{eq:shad}
g_n^{\pm1}x \in V(x,\st(\tau_n^{\pm})) .
\end{equation}
Here we use that the sequence $(g_n^{-1})$ is $\iota\taumod$-regular,
compare the comment after Definition~\ref{def:pureg}.
\begin{lem}
\label{lem:flconvcontr}
If $\tau_n^{\pm}\to\tau_{\pm}$ in $\Flagpmt$,
then $(g_n)$ is $\taumod$-contracting with (\ref{eq:contrtau}).
\end{lem}
\proof
Compare the proof of \cite[Proposition 5.14]{morse}.

Since $x\in g_nV(x,\st(\tau_n^-))=V(g_nx,\st(g_n\tau_n^-))$,
it follows together with $g_nx \in V(x,\st(\tau_n^+))$
that the Weyl cones $V(g_nx,\st(g_n\tau_n^-))$ and $V(x,\st(\tau_n^+))$ 
lie in the same parallel set, namely in $P(g_n\tau_n^-,\tau_n^+)$, and face in opposite directions. 
In particular, the simplices $g_n\tau_n^-$ and $\tau_n^+$ are $x$-opposite,
and thus $g_n\tau_n^-$ converges to the simplex $\hat\tau_+$ which is $x$-opposite to $\tau_+$,
$$ g_n\tau_n^- \to\hat\tau_+ .$$
Since the sequence $(g_n^{-1}x)$ is $\iota\taumod$-regular,
it holds that 
$$ d(g_n^{-1}x,\D V(x,\st(\tau_n^-))) \to+\infty  $$
According to Lemma~\ref{lem:expconvsect} in appendix~\ref{sec:shadw},
for any $r,R>0$
the inclusion of shadows
$$ U_{\tau_n^-,x,R} \subset U_{\tau_n^-,g_n^{-1}x,r} $$
holds for $n\geq n(r,R)$.
Therefore there exist positive numbers $R_n\to+\infty$ and $r_n\to0$ such that 
$$ U_{\tau_n^-,x,R_n} \subset U_{\tau_n^-,g_n^{-1}x,r_n} $$
for large $n$,
equivalently
\begin{equation}
\label{eq:mpfshd}
g_nU_{\tau_n^-,x,R_n} \subset U_{g_n\tau_n^-,x,r_n} .
\end{equation}
Since $\tau_n^-\to\tau_-$ and $R_n\to+\infty$,
the sequence of shadows $U_{\tau_n^-,x,R_n}\subset C(\tau_n^-)\subset\Flagt$ {\em exhausts} $C(\tau_-)$
in the sense that every compactum in $C(\tau_-)$ is contained in $U_{\tau_n^-,x,R_n}$ for large $n$.\footnote{Indeed, 
for fixed $R>0$ we have Hausdorff convergence 
$U_{\tau_n^-,x,R}\to U_{\tau_-,x,R}$ in $\Flagt$,
which immediately follows e.g.\ from the transitivity of the action $K_x\acts\Flagit$
of the maximal compact subgroup $K_x<G$ fixing $x$.
Furthermore,
the shadows $U_{\tau_-,x,R}$ exhaust $C(\tau_-)$ as $R\to+\infty$,
cf.\ the continuity part of Lemma~\ref{lem:contpr}.}
On the other hand,
since $g_n\tau_n^-\to\hat\tau_+$ and $r_n\to0$,
the shadows 
$U_{g_n\tau_n^-,x,r_n}$ {\em shrink}, i.e.\ Hausdorff converge to the point $\tau_+$.\footnote{Indeed, 
$U_{g_n\tau_n^-,x,r}\to U_{\hat\tau_+,x,r}$ in $\Flagt$ for fixed $r>0$,
and $U_{\hat\tau_+,x,r}\to\tau_+$ as $r\to0$,
using again the continuity part of Lemma~\ref{lem:contpr}
and the fact that the function (\ref{eq:distfrpar}) assumes the value zero only in $\tau_+$.}
Together with these observations on exhaustion and shrinking of shadows, 
(\ref{eq:mpfshd}) shows that 
$$ g_n|_{C(\tau_-)} \to \tau_+ $$
uniformly on compacts,
i.e.\ the (sub)sequence $(g_n)$ is $\taumod$-contracting. 
\qed

\medskip
With the lemma, we can add a converse to Proposition~\ref{prop:contrimpreg} and obtain:
\begin{prop}
\label{prop:regequivcontr}
The following properties are equivalent for sequences in $G$:

(i) Every subsequence contains a $\taumod$-contracting subsequence. 

(ii) The sequence is $\taumod$-regular.
\end{prop}
\proof
This is a direct consequence of the lemma. 
For the implication (ii)$\Ra$(i) one uses the compactness of flag manifolds,
and for the implication (i)$\Ra$(ii) the following fact:
If a sequence is not $\taumod$-regular, then it contains a $\numod$-pure subsequence
for a $\numod\not\supset\taumod$,
and hence in particular a subsequence 
none of whose subsequences is $\taumod$-regular.
This subsequence then cannot have any $\taumod$-contracting subsequence
by Proposition~\ref{prop:contrimpreg}.
\qed

\medskip
We conclude for subgroups:
\begin{thm}
\label{thm:regimplcontrgp}
A discrete subgroup $\Ga<G$ is $\taumod$-regular iff it is a $\taumod$-convergence subgroup.
\end{thm}
\proof
By definition, 
$\Ga$ is $\taumod$-regular iff every sequence $\ga_n\to\infty$ in $\Ga$ is $\taumod$-regular,
and $\taumod$-convergence iff every sequence has a $\taumod$-contracting subsequence.
According to the proposition, these conditions are equivalent.
\qed

\subsection{Convergence at infinity and limit sets}
\label{sec:conv}

The discussion in the preceeding two sections 
leads to 
a natural notion of convergence at infinity 
for regular sequences in $X$ and $G$,
compare \cite[\S 5.3]{morse}.
As regularity, 
it can be expressed both in terms of orbits in $X$ and dynamics on flag manifolds.

Let first $(g_n)$ be a $\taumod$-regular sequence in $G$.
Flexibilizing condition (\ref{eq:shad}),
we choose points $x,x'\in X$ and consider a sequence $(\tau_n)$ in $\Flagt$ 
such that 
\begin{equation}
\label{eq:bdddstfrcona}
\sup_nd\bigl(g_nx,V(x',\st(\tau_n))\bigr) < +\infty.
\end{equation}
Note that the condition is independent of the choice of the points $x$ and $x'$.\footnote{Recall that 
the Hausdorff distance of asymptotic Weyl cones $V(y,\st(\tau))$ and $V(y',\st(\tau))$
is bounded by the distance $d(y,y')$ of their tips.}

\begin{lem}
\label{lem:regseqindepch}
The accumulation set of $(\tau_n)$ in $\Flagt$ 
only depends on $(g_n)$.
\end{lem}
\proof
Let $(\tau'_n)$ be another sequence 
such that $d(g_nx,V(x',\st(\tau'_n)))$ is uniformly bounded. 
Assume that after extraction $\tau_n\to\tau$ and $\tau'_n\to\tau'$.
We must show that $\tau=\tau'$.

We may suppose that $x'=x$.
There exist bounded sequences $(b_n)$ and $(b'_n)$ in $G$
such that 
$$g_nb_nx\in V(x,\st(\tau_n)) 
\qquad\hbox{ and }\qquad 
g_nb'_nx\in V(x,\st(\tau'_n))$$
for all $n$.
Note that the sequences $(g_nb_n)$ and $(g_nb'_n)$ in $G$ are again $\taumod$-regular.
By Lemma~\ref{lem:flconvcontr},
after further extraction 
they are $\taumod$-contracting with 
$$g_nb_n|_{C(\tau_-)}\to\tau
\qquad\hbox{ and }\qquad 
g_nb'_n|_{C(\tau'_-)}\to\tau'$$
uniformly on compacts for some $\tau_-,\tau'_-\in\Flagt$.
Moreover, we may assume convergence $b_n\to b$ and $b'_n\to b'$.
Then 
$$g_n|_{C(b\tau_-)}\to\tau
\qquad\hbox{ and }\qquad 
g_n|_{C(b'\tau'_-)}\to\tau'$$
uniformly on compacts.
With Lemma~\ref{lem:contruniq} it follows that $\tau=\tau'$.
\qed

\medskip
In view of the lemma, 
we define,
cf.\ \cite[Def.\ 5.26]{morse}:
\begin{dfn}[Flag convergence of sequences in $G$]
A $\taumod$-regular sequence $(g_n)$ in $G$ is said to 
{\em $\taumod$-flag converge} 
to a simplex $\tau\in\Flagt$,
$$ g_n\to\tau ,$$
if $\tau_n\to\tau$ in $\Flagt$
for some sequence $(\tau_n)$ in $\Flagt$ 
satisfying (\ref{eq:bdddstfrcona}).
\end{dfn}

We can now characterize contraction in terms of flag convergence.
We rephrase Lemma~\ref{lem:flconvcontr}
and show that its exact converse holds as well:

\begin{lem}
\label{lem:flconvcontrconvreph}
For a sequence $(g_n)$ in $G$,
the following are equivalent:

(i) $(g_n)$ is $\taumod$-contracting 
with 
$g_n|_{C(\tau_-)}\to\tau_+$ uniformly on compacts.

(ii) $(g_n)$ is $\taumod$-regular and $g_n^{\pm1}\to\tau_{\pm}$.
\end{lem}
\proof
The implication (ii)$\Ra$(i) follows from Lemma~\ref{lem:flconvcontr}.
Suppose that (i) holds.
Then the sequence $(g_n)$ is $\taumod$-regular by Proposition~\ref{prop:contrimpreg}. 
Let $(\tau_n^{\pm})$ be sequences satisfying (\ref{eq:shad}).
We must show that $\tau_n^{\pm}\to\tau_{\pm}$.
Otherwise, after extraction we obtain that
$\tau_n^{\pm}\to\tau'_{\pm}$
with $\tau'_+\neq\tau_+$ or $\tau'_-\neq\tau_-$.
Then also 
$g_n|_{C(\tau'_-)}\to\tau'_+$
by Lemma~\ref{lem:flconvcontr},
and Lemma~\ref{lem:contruniq} implies that $\tau'_{\pm}=\tau_{\pm}$,
a contradiction.
\qed

\medskip
Vice versa,
we can characterize flag convergence in terms of contraction
and thus give an alternative dynamical definition of flag convergence:

\begin{lem}
\label{lem:flconvitcontr}
For a sequence $(g_n)$ in $G$, 
the following are equivalent:

(i) $(g_n)$ is $\taumod$-regular and $g_n\to\tau$.

(ii) There exists a bounded sequence $(b_n)$ in $G$ and $\tau_-\in\Flagit$
such that 
$g_nb_n|_{C(\tau_-)}\to\tau$ uniformly on compacts.

(iii) There exists a bounded sequence $(b'_n)$ in $G$ 
such that 
$b'_ng_n^{-1}|_{C(\tau)}$ converges to a constant map uniformly on compacts.
\end{lem}
\proof
(ii)$\Ra$(i):
According to the previous lemma 
the sequence 
$(g_nb_n)$ is $\taumod$-regular and $\taumod$-flag converges, 
$g_nb_n\to\tau$.
Since $d(g_nx,g_nb_nx)$ is uniformly bounded,
this is equivalent to $(g_n)$ being $\taumod$-regular and $g_n\to\tau$.

(i)$\Ra$(ii):
Due to the compactness of flag manifolds,
there exists a bounded sequence $(b'_n)$ in $G$ such that $(b'_ng^{-1}_n)$ 
$\iota\taumod$-flag converges,
$b'_ng^{-1}_n\to\tau_-\in\Flagit$.
We put $b_n={b'_n}^{-1}$.
Since also $(g_nb_n)$ is $\taumod$-regular and $g_nb_n\to\tau$,
it follows 
from the previous lemma 
that $g_nb_n|_{C(\tau_-)}\to\tau$ uniformly on compacts.

The equivalence (ii)$\Leftrightarrow$(iii) with $b'_n=b_n^{-1}$ follows from Lemma~\ref{lem:contrsym}.
\qed

\medskip
We extend the notion of flag convergence to sequences in $X$.

Let now $(x_n)$ be a $\taumod$-regular sequence in $X$. 
We choose again a base point $x\in X$ and consider a sequence $(\tau_n)$ in $\Flagt$ such that 
\begin{equation}
\label{eq:bdddstfrcon}
\sup_nd\bigl(x_n,V(x,\st(\tau_n))\bigr) < +\infty
\end{equation}
analogous to (\ref{eq:bdddstfrcona}).
As before, the condition is independent of the choice of the point $x$, and:

\begin{lem}
The accumulation set of $(\tau_n)$ in $\Flagt$ only depends on $(x_n)$.
\end{lem}
\proof
Let $(g_n)$ be a sequence in $G$ such that the sequence $(g_n^{-1}x_n)$ in $X$ is bounded. 
Then $(g_n)$ is $\taumod$-regular and its $\taumod$-flag accumulation set in $\Flagt$ equals 
the accumulation set of the sequence $(\tau_n)$.
This reduces the assertion to Lemma~\ref{lem:regseqindepch}.
\qed

\medskip
We therefore can define:
\begin{dfn}[Flag convergence of sequences in $X$]
A $\taumod$-regular sequence $(x_n)$ in $X$ is said to 
{\em $\taumod$-flag converge} 
to a simplex $\tau\in\Flagt$,
$$ x_n\to\tau ,$$
if $\tau_n\to\tau$ in $\Flagt$ for some sequence $(\tau_n)$ in $\Flagt$ satisfying (\ref{eq:bdddstfrcon}).
\end{dfn}

For any $\taumod$-regular sequence $(g_n)$ in $G$ and any point $x\in X$,
we have 
$g_n\to\tau$ iff $g_nx\to\tau$.

Flag convergence and flag limits are stable under bounded perturbations of sequences:
\begin{lem}
(i)
For any $\taumod$-regular sequence $(g_n)$ and any bounded sequence $(b_n)$ in $G$,
the sequences $(g_n)$ and $(g_nb_n)$ have the same $\taumod$-flag accumulation sets in $\Flagt$.

(ii)
If $(x_n)$ and $(x'_n)$ are $\taumod$-regular sequences in $X$ such that $d(x_n,x'_n)$ is uniformly bounded,
then both sequences have the same $\taumod$-flag accumulation set in $\Flagt$.
\end{lem}
\proof
(i)
The sequence $(g_nb_n)$ is also regular and satisfies condition (\ref{eq:bdddstfrcona}) iff $(g_n)$ does.

(ii)
The sequence $(x'_n)$ satisfies condition (\ref{eq:bdddstfrcon}) iff $(x'_n)$ does. 
\qed

\begin{rem}
\label{rem:flconvfins}
There is a natural topology on the bordification $X\sqcup\Flagt$
which induces $\taumod$-flag convergence.
Moreover,
the bordification embeds into a natural Finsler compactification of $X$,
compare Remark~\ref{rem:regfins}.
\end{rem}

From flag convergence, 
we obtain for discrete subgroups the following notion of limit sets in flag manifolds,
compare  \cite[Def.\ 6.9]{coco15} and \cite[Def.\ 5.32]{morse}:

\begin{dfn}[Flag limit set]
For a discrete subgroup $\Ga<G$ we define its {\em $\taumod$-limit set} 
$$\Lat\subset\Flagt$$
as the set of possible limit simplices of 
$\taumod$-flag converging $\taumod$-regular sequences in $\Ga$,
equivalently,
as the set of simplices $\tau_+$ as in (\ref{eq:contrtau})
for all $\taumod$-contracting sequences in $\Ga$.
\end{dfn}

The limit set is $\Ga$-invariant and closed,
as a diagonal argument shows.

A discrete subgroup is {\em $\taumod$-antipodal} for a $\iota$-invariant face type $\taumod$,
if $\Lat$ is antipodal in the sense that its elements are pairwise antipodal. 

\subsection{Uniform regularity}

We introduce stronger forms of the regularity conditions discussed in section~\ref{sec:reg}.

We first consider sequences in the model euclidean Weyl chamber $\De$.
\begin{dfn}
\label{def:unifreg}
A sequence $\de_n\to\infty$ in $\De$ is 
{\em uniformly $\taumod$-regular} if it drifts away from $V(0,\D\st(\taumod))\subset\D\De$
at a linear rate with respect to $\|\de_n\|$,
$$ \liminf_{n\to+\infty} d(\de_n,V(0,\D\st(\taumod))) / \|\de_n\| > 0.$$
\end{dfn}

We extend these notions to sequences in $X$ and $G$:
\begin{dfn}[Uniformly regular]
(i) A sequence $(x_n)$ in $X$ is {\em uniformly $\taumod$-regular}
if for some (any) base point $o\in X$ the sequence of $\De$-distances $d_{\De}(o,x_n)$ 
has this property.

(ii) A sequence $(g_n)$ in $G$ is {\em uniformly $\taumod$-regular}
if for some (any) point $x\in X$ the orbit sequence $(g_nx)$ in $X$ has this property.

(iii) A subgroup $\Ga<G$ is {\em uniformly $\taumod$-regular} if all sequences of pairwise distinct elements in $\Ga$ have this property.
\end{dfn}

For a subgroup $\Ga<G$ uniform $\taumod$-regularity is equivalent to the property that 
the visual limit set $\La\subset\geo X$ is contained in the union $\geo^{\taumod-reg} X$ 
of the open $\taumod$-stars.

A subgroup $\Ga<G$ is uniformly $\taumod$-regular iff it is uniformly $\iota\taumod$-regular.

\section{Finsler compactifications of symmetric spaces}
\label{sec:fico}

Let $X=G/K$ be a symmetric space of noncompact type.

\subsection{Finsler metrics}
\subsubsection{The Riemannian distance}

We denote by $d^{Riem}$ the $G$-invariant Riemannian distance on $X$.

Let $xy\subset X$ be an oriented geodesic segment. 
The Busemann functions $b_{\xi}$ for $\xi\in\geo X$ have slope $\geq-1$ along $xy$,
because they are 1-Lipschitz.
Therefore 
\begin{equation}
\label{ineq:ribu}
b_{\xi}(x)-b_{\xi}(y) \leq d^{Riem}(x,y)
\end{equation}
with equality iff $y\in x\xi$.
Therefore, the Riemannian distance can be represented in the form:
\begin{equation}
\label{eq:riemmbu}
d^{Riem}(x,y) = \max_{\xi\in\geo X} \bigl(b_{\xi}(x)-b_{\xi}(y)\bigr) 
\end{equation}

\subsubsection{Finsler distances}
\label{sec:finsdist}

We fix a type $\bar\theta\in\simod$
and now work only with Busemann functions $b_{\xi}$ of this type, $\theta(\xi)=\bar\theta$. 
There is the following sharper bound for the slopes of such Busemann functions along segments:

\begin{lem}
\label{lem:slpbd}
The slope of a Busemann function $b_{\xi}$ of type $\theta(\xi)=\bar\theta$ 
along a non-degenerate oriented segment $xy\subset X$ is 
$\geq-\cos\angle(\theta(xy),\bar\theta)$
with equality in some point, equivalently, along the entire segment, iff $y\in V(x,\st(\xi))$.
\end{lem}
\proof
The slope of $b_{\xi}|_{xy}$ in an interior point $z\in xy$ equals 
$-\cos\angle_z(y,\xi)$. 
The angle $\angle_z(y,\xi)$ assumes its minimal value $\angle(\theta(zy),\bar\theta)=\angle(\theta(xy),\bar\theta)$
iff the segment $zy$ and the ray $z\xi$ are contained in a euclidean Weyl chamber 
with tip at $z$,
equivalently,
if $xy$ and $x\xi$ are contained in a euclidean Weyl chamber with tip at $x$,
equivalently, 
if $y\in V(x,\st(\xi))$.
In this case, 
the slope of $b_{\xi}$ equals $\angle(\theta(xy),\bar\theta)$ along the entire segment $xy$.
\qed

\medskip
We define the $G$-invariant {\em $\bar\theta$-Finsler distance} $d^{\bar\theta}:X\times X\to [0,+\infty)$ by 
\begin{equation}
\label{eq:fimbu}
d^{\bar\theta}(x,y) := \max_{\theta(\xi)=\bar\theta} \bigl( b_{\xi}(x)-b_{\xi}(y) \bigr) 
\end{equation}
where the maximum is taken over all ideal points $\xi\in\geo X$ with type $\theta(\xi)=\bar\theta$.
By analogy with (\ref{ineq:ribu}), we have the inequality 
\begin{equation}
\label{ineq:bincr}
b_{\xi}(x)-b_{\xi}(y) \leq d^{\bar\theta}(x,y) 
\end{equation}
for all $\xi\in\geo X$ with $\theta(\xi)=\bar\theta$.
According to the lemma,
equality holds iff $y\in V(x,\st(\xi))$.

The triangle inequality is clearly satisfied for $d^{\bar\theta}$.
In view of $\diam(\simod)\leq\pihalf$ and the lemma we have semipositivity, i.e.\ $d^{\bar\theta}\geq0$.
Regarding symmetry, we have the identity
\begin{equation}
\label{eq:fisy}
d^{\iota\bar\theta}(y,x)  = d^{\bar\theta}(x,y)
\end{equation}
and hence $d^{\bar\theta}$ is symmetric iff $\iota\bar\theta=\bar\theta$.
To see (\ref{eq:fisy}) we note that, according to the lemma,
$b_{\xi}$ has maximal decay along $xy$ iff 
$b_{\hat\xi}$ has maximal decay along $yx$,
where $\hat\xi\in\geo X$ denotes the ideal point 
which is $x$-opposite to $\xi$
and has type $\iota\bar\theta$. 

The distance $d^{\bar\theta}$ can be derived from the vector valued $\De$-distance $d_{\De}$ by
composing it with the linear functional $l_{\bar\theta}=-b_{\bar\theta}$ on $\Fmod\supset\De_{mod}$
(normalized at the origin):
$$ d^{\bar\theta} =l_{\bar\theta}\circ d_{\De} $$
Let $F\subset X$ be a maximal flat.
The restriction of the distance $d^{\bar\theta}$ to $F$
can be written intrinsically as 
\begin{equation}
\label{eq:distffl}
d^{\bar\theta}(x,y) = \max_{\xi\in\geo F,\theta(\xi)=\bar\theta} \bigl( b_{\xi}(x)-b_{\xi}(y) \bigr) 
\end{equation}
for $x,y\in F$,
because equality holds in (\ref{ineq:bincr}) if $\xi$ lies in a chamber $\si\subset\geo F$ 
with $y\in V(x,\si)$.
The restriction of $d^{\bar\theta}$ to a maximal flat is thus the translation invariant pseudo-metric 
associated to the $W$-invariant polyhedral seminorm on $\Fmod$ given by
\begin{equation}
\label{eq:plynrm}
\| \cdot \|_{\bar\theta} =  \max_{w\in W} (l_{\bar\theta}\circ w^{-1}) .
\end{equation}
If the seminorm $\| \cdot \|_{\bar\theta}$ is a norm,
equivalently, if $d^{\bar\theta}$ is a (non-symmetric) metric, 
then $d^{\bar\theta}$ is equivalent to the Riemannian distance $d^{Riem}$.
We describe when this is the case:

\begin{lem}[Positivity]
\label{lem:pos}
The following are equivalent:

(i) $d^{\bar\theta}$ is a (non-symmetric) metric.

(ii) The radius of $\simod$ with respect to $\bar\theta$ is $<\pihalf$.

(iii) $\bar\theta$ is not contained in a factor of a nontrivial spherical join decomposition of $\simod$. 
\end{lem}
\proof
The equivalence (i)$\Leftrightarrow$(ii) is immediate in view of (\ref{eq:plynrm}). 

To see (ii)$\Leftrightarrow$(iii),
consider the spherical join decomposition of $\simod$ into its irreducible factors $\simod^i$.
These have diameter $<\pihalf$.
We work now in the model apartment and represent directions by unit vectors
which we orthogonally decompose into their $\simod^i$-components.
Any two vectors in $\De_{mod}$ with nontrivial $\simod^i$-components for some $i$
have angle $<\pihalf$.
This yields the implication (iii)$\Ra$(ii). 
The converse direction is clear. 
\qed

\medskip
In particular,
$d^{\bar\theta}$ is a (non-symmetric) metric 
if $\bar\theta$ is regular or if $X$ is irreducible.

If $d^{\bar\theta}$ is only a pseudo-metric,
then $X$ splits as a product $X_1\times X_2$ such that $d^{\bar\theta}$ is degenerate 
precisely in the $X_2$-direction and induces an honest (non-symmetric) metric on $X_1$.

\subsubsection{Geodesics}
\label{sec:geodesics}

We first analyze when equality holds in the triangle inequality for the Finsler distance. 
Let $\taumod$ denote the face type spanned by $\bar\theta$, i.e.\ $\bar\theta\in\interior(\taumod)$.

\begin{lem}\label{lem:triangle_equality}
A triple of points $x,y,z\in X$ satisfies
$$ d^{\bar\theta}(x,z) + d^{\bar\theta}(z,y) =d^{\bar\theta}(x,y) $$
iff it is contained in a parallel set $P(\tau_-,\tau_+)$
for a pair of opposite simplices $\tau_{\pm}\in\Flagpmt$
and $z$ lies in the $\taumod$-diamond 
determined by $x,y$:
\begin{equation}
\label{eq:zindiamo}
z\in \diamot(x,y)= V(x,\st(\tau_+))\cap V(y,\st(\tau_-)) 
\end{equation}
\end{lem}
\proof Assume that the equality holds.
Let $\tau_+\in\Flagt$ be a simplex such that $y\in V(x,\st(\tau_+))$, 
and let $\xi_+\in\tau_+$ be the ideal point with type $\theta(\xi_+)=\bar\theta$.
Then $b_{\xi_+}$ has maximal decay along $xy$,
cf.\ Lemma~\ref{lem:slpbd}.
From 
$$ d^{\bar\theta}(x,y)=
d^{\bar\theta}(x,z) + d^{\bar\theta}(z,y) \geq 
\bigl(b_{\xi_+}(x)-b_{\xi_+}(z)\bigr)  + \bigl(b_{\xi_+}(z)-b_{\xi_+}(y)\bigr)  
= b_{\xi_+}(x)-b_{\xi_+}(y)  =  d^{\bar\theta}(x,y)  $$
it follows that $b_{\xi_+}$ must have maximal decay also along the segments $xz$ and $zy$.
This implies that $z\in V(x,\st(\tau_+))$,
again by the same lemma.
Furthermore, 
$b_{\xi_-}$ has maximal decay along $yx$ 
for the ideal point $\xi_-$ which is $x$-opposite to $\xi_+$
and therefore contained in the simplex $\tau_-$ $x$-opposite to $\tau_+$.
It follows that also $z\in V(y,\st(\tau_-))$.

Conversely, if (\ref{eq:zindiamo}) holds,
then $b_{\xi_+}$ has maximal decay along $xy$, $yz$ and $xz$,
and hence the equality is satisfied. 
\qed

\medskip
It follows 
that the (pseudo-)metric space $(X, d^{\bar\theta})$ is a {\em geodesic} space.
The Riemannian geodesics in $X$ are also $d^{\bar\theta}$-geodesics,
but besides these there are other $d^{\bar\theta}$-geodesics,
due to the non-strict convexity of balls for the norm $\| \cdot \|_{\bar\theta}$.

The lemma yields a precise description of all $d^{\bar\theta}$-geodesics:
A path $c:I\to X$ is an (unparametrized) $d^{\bar\theta}$-geodesic 
iff it is contained in a parallel set $P(\tau_-,\tau_+)$ with $\tau_{\pm}\in\Flagpmt$
and
$$ \hbox{$c(t')\in V(c(t),\st(\tau_+))$, equivalently, $c(t)\in V(c(t'),\st(\tau_-))$} $$
for all $t<t'$ in $I$,
i.e. $c$ {\em drifts} towards $\tau_+$ and away from $\tau_-$. 
As a consequence, a geodesic $c:[t_-,t_+]\to X$ 
is contained the diamond $\diamot(c(t_-),c(t_+))$ determined by its endpoints.
Moreover, we obtain the {\em Finsler geometric interpretation of diamonds}, 
namely the diamond $\diamot(x,y)$ is the union of all $d^{\bar\theta}$-geodesics $xy$.

\medskip
The most relevant case for this paper is when $\bar\theta$ is {\em regular}, $\bar\theta\in\interior(\simod)$.
The above discussion then specializes as follows:
The pair of simplices $\tau_{\pm}$ in the lemma becomes a pair of opposite chambers $\si_{\pm}$,
the parallel set $P(\tau_-,\tau_+)$ becomes a maximal flat $F(\si_-,\si_+)$, 
the Weyl cones $V(\cdot,\st(\tau_{\pm}))$ become euclidean Weyl chambers $V(\cdot,\si_{\pm})$.
Thus (\ref{eq:zindiamo}) simplifies to 
\begin{equation*}
z\in V(x,\si_+)\cap V(y,\si_-) ,
\end{equation*}
and a $d^{\bar\theta}$-geodesic $c:I\to X$ is contained in a maximal flat $F(\si_-,\si_+)$
and
$$ \hbox{$c(t')\in V(c(t),\si_+)$, equivalently, $c(t)\in V(c(t'),\si_-)$} $$
for all $t<t'$ in $I$.

\subsection{Finsler compactifications}

Throughout this section we assume that the type $\bar\theta$ is regular, $\bar\theta\in\interior(\simod)$.
In particular, $d^{\bar\theta}$ is a metric. 

\subsubsection{Definition}

If one applies the horoboundary construction, cf. section~\ref{sec:horob},
to the Riemannian distance $d^{Riem}$ on $X$,
one obtains the visual compactification
\begin{equation}
\label{eq:viscomp}
\ol X = X \sqcup \geo X .
\end{equation}
The ideal boundary points are represented by Busemann functions,
i.e. the horofunctions are in this case precisely the Busemann functions. 

We define the {\em $\bar\theta$-Finsler compactification} of $X$
as the compactification 
\begin{equation}
\label{eq:finscomp}
\ol X^{\bar\theta} = X \sqcup \geo^{\bar\theta} X .
\end{equation}
obtained by applying the horoboundary construction 
to the Finsler distance $d^{\bar\theta}$.

\subsubsection{Horofunctions}
\label{sec:mix}

For a chamber $\si\subset\geo X$, 
let $\theta_{\si}\in\si$ denote the point of type $\bar\theta$.
The associated Busemann function $b_{\theta_{\si}}$ is well-defined up to additive constant,
and the Busemann function $b_{\theta_{\si}}-b_{\theta_{\si}}(x)$ normalized in a point $x\in X$ is well-defined.

According to our definition (\ref{eq:fimbu}) of the $d^{\bar\theta}$-distance, we have 
$$ d^{\bar\theta}_x:= d^{\bar\theta}(\cdot,x) = \max_{\si} \bigl( b_{\theta_{\si}}-b_{\theta_{\si}}(x) \bigr) $$
where the maximum is taken over all chambers $\si$.
For a simplex $\tau\subset\geo X$ and a point $x\in X$,
we consider the ``mixed'' Busemann function
\begin{equation}
\label{eq:mbf}
b^{\bar\theta}_{\tau,x} := 
\max_{\si\supset\tau}\bigl(b_{\theta_{\si}}-b_{\theta_{\si}}(x)\bigr)
\end{equation}
normalized in $x$, 
the maximum being taken only over the chambers which contain $\tau$ as a face. 
We will see that these are precisely the horofunctions for $\ol X^{\bar\theta}$.

On a euclidean Weyl chamber with tip at $x$,
the function $d^{\bar\theta}_x$ agrees with one of the Busemann functions occuring in the maximum:
If $\si,\hat\si\subset\geo X$ are $x$-opposite chambers,
equivalently, if $x\in F(\hat\si,\si)$,
then
$$ d^{\bar\theta}(\cdot,x) =b_{\theta_{\si}}-b_{\theta_{\si}}(x) $$
on $V(x,\hat\si)$,
cf.\ Lemma~\ref{lem:slpbd}.
Thus, on a Weyl cone with tip at $x$,
the function $d^{\bar\theta}_x$ reduces to a maximum over a subfamily of Busemann functions:
If $\tau,\hat\tau\subset\geo X$ are $x$-opposite simplices,
equivalently, if $x\in P(\hat\tau,\tau)$,
then 
$$ d^{\bar\theta}(\cdot,x) =\max_{\si\supset\tau}\bigl(b_{\theta_{\si}}-b_{\theta_{\si}}(x)\bigr) 
=b^{\bar\theta}_{\tau,x}$$
on $V(x,\st(\hat\tau))$.
For the normalized distance functions,
we observe:
If $o\in X$ is a base point and if $x$ lies in the Weyl sector $V(o,\tau)$,
equivalently, 
if $o\in V(x,\hat\tau)=\cap_{\hat\si\supset\hat\tau}V(x,\hat\si)$, 
then the difference 
$b_{\theta_{\si}}(o)-b_{\theta_{\si}}(x)$
has the same value for all chambers $\si\supset\tau$,
and hence the function $d^{\bar\theta}_x$ normalized in $o$ is given by 
\begin{equation}
\label{eq:ndonwc}
d^{\bar\theta}_x - d^{\bar\theta}_x(o) 
=\max_{\si\supset\tau}\bigl(b_{\theta_{\si}}-b_{\theta_{\si}}(o)\bigr) 
=b^{\bar\theta}_{\tau,o}
\end{equation}
on the Weyl cone $V(x,\st(\hat\tau))$.

With these observations
we are prepared for understanding the horofunctions,
i.e.\ the limits of (normalized) distance functions $d^{\bar\theta}_x$ as $x\to\infty$.
We first show that the mixed Busemann functions $b^{\bar\theta}_{\tau,x}$
are horofunctions.

\begin{lem}
\label{lem:limalsct}
Let $o\in X$ and $\tau\in\Flagt$.
If $(x_n)$ is a $\taumod$-regular sequence in the Weyl sector $V(o,\tau)$,
then 
$$ d^{\bar\theta}_{x_n} - d^{\bar\theta}_{x_n}(o) \to 
b^{\bar\theta}_{\tau,o}$$
uniformly on compacts in $X$.
\end{lem}
\proof
{\em Step 1.}
Let $\hat\tau$ be the simplex $o$-opposite to $\tau$.
We first note that the claimed convergence holds uniformly on compacts in $P(\hat\tau,\tau)$.
In fact, on every such compact, eventually equality holds.
This follows from (\ref{eq:ndonwc}) and because the $\taumod$-regularity of the sequence $(x_n)$ implies that 
the Weyl cones $V(x_n,\st(\hat\tau))$ exhaust $P(\hat\tau,\tau)$ as $n\to+\infty$.

{\em Step 2.}
To verify the convergence on all of $X$,
we use the action of the horocyclic subgroup $H_{\tau}$.
The Busemann functions $b_{\xi}$ centered at ideal points $\xi\in\st(\tau)$
are $H_{\tau}$-invariant,
$b_{\xi}\circ h^{-1}=b_{h\xi}=b_{\xi}$ for $h\in H_{\tau}$,
and hence also the mixed Busemann functions,
$$ b^{\bar\theta}_{\tau,x}\circ h^{-1}= b^{\bar\theta}_{\tau,x} .$$
By step 1, it holds for $(h,x)\in H_{\tau}\times P(\hat\tau,\tau)$ that 
\begin{equation}
\label{eq:fllfr}
\underbrace{d^{\bar\theta}_{hx_n}(hx)}_{d^{\bar\theta}_{x_n}(x)} - d^{\bar\theta}_{x_n}(o) 
\to \underbrace{b^{\bar\theta}_{\tau,o}(hx)}_{b^{\bar\theta}_{\tau,o}(x)} ,
\end{equation}
and the convergence is uniform on $H_{\tau}\times A$ for $A\subset P(\hat\tau,\tau)$ compact. 
Note that 
$$ d^{\bar\theta}(x_n,hx_n)=d^{\bar\theta}(h^{-1}x_n,x_n) \to0 $$
as $n\to+\infty$ locally uniformly in $h$, because $(x_n)$ is a $\taumod$-regular sequence in $V(o,\tau)$
and therefore drifts away from $\D V(o,\tau)$.
Hence 
$$ d^{\bar\theta}_{hx_n}-d^{\bar\theta}_{x_n} \to0 $$
uniformly on $X$
due to the triangle inequality,
compare (\ref{ineq:osc}).
It follows from (\ref{eq:fllfr}) that 
$$ d^{\bar\theta}_{x_n}(hx) - d^{\bar\theta}_{x_n}(o) 
\to b^{\bar\theta}_{\tau,o}(hx) $$
locally uniformly in $(h,x)$, i.e.
$$ d^{\bar\theta}_{x_n} - d^{\bar\theta}_{x_n}(o) 
\to b^{\bar\theta}_{\tau,o} $$
locally uniformly on $X$, as claimed.
\qed

\medskip
We show next that,
vice versa,
there are no other horofunctions besides the mixed Busemann functions $b^{\bar\theta}_{\tau,x}$:

\begin{lem}
Let $x_n\to\infty$ be a divergent sequence in $X$.
Then, after extraction, there exist a simplex $\tau\subset\geo X$ and a point $p\in X$ such that 
$$ d^{\bar\theta}_{x_n} - d^{\bar\theta}_{x_n}(p) \to b^{\bar\theta}_{\tau,p}$$
uniformly on compacts in $X$.
\end{lem}
\proof
We reduce the assertion to the previous lemma 
using the action of the maximal compact subgroup $K<G$ fixing a base point $o\in X$.
There exists a sequence $(k_n)$ in $K$ such that the sequence $(k_nx_n)$ 
is contained in a fixed euclidean Weyl chamber $V(o,\si)$.
After extraction, we may assume that $k_n\to e$.
We may assume moreover
that the sequence $(x_n)$, equivalently, $(k_nx_n)$ is $\taumod$-pure for some face type $\taumod$.
Let $\tau\subset\si$ be the face of type $\taumod$.
After further extraction, there exists a point $p\in V(o,\si)$ 
such that $(k_nx_n)$ approaches the Weyl sector $V(p,\tau)\subset V(o,\si)$,
i.e.\ there exists a 
sequence $(y_n)$ in $V(p,\tau)$
such that $d(k_nx_n,y_n)\to0$.
The sequence $(y_n)$ is then also $\taumod$-pure, and in particular $\taumod$-regular.
By the previous lemma,
$d^{\bar\theta}_{y_n} - d^{\bar\theta}_{y_n}(p) \to b^{\bar\theta}_{\tau,p}$
uniformly on compacts, and hence 
$$ d^{\bar\theta}_{k_nx_n} - d^{\bar\theta}_{k_nx_n}(p) \to b^{\bar\theta}_{\tau,p} $$
uniformly on compacts.
Since 
$b^{\bar\theta}_{\tau,p}\circ k_n\to b^{\bar\theta}_{\tau,p}$
uniformly on compacts, 
the assertion follows. 
\qed

\medskip
Thus,
the horofunctions for $\ol X^{\bar\theta}$ are precisely the mixed Busemann functions $b^{\bar\theta}_{\tau,p}$.

We now discuss some of their properties.
As already mentioned,
the functions $b^{\bar\theta}_{\tau,p}$ are invariant under the horocyclic subgroup $H_{\tau}$.
Moreover, they are invariant up to additive constants under transvections towards $\tau$:
\begin{lem}
\label{lem:transvtaubeh}
For a transvection $t$ with axes asymptotic to 
$\xi\in\tau$
it holds that 
$b^{\bar\theta}_{\tau,p}\circ t^{-1}\equiv b^{\bar\theta}_{\tau,p}$.\footnote{Recall 
that the notation $f\equiv g$ for functions $f,g$ means that  $f-g$ is a constant.}
\end{lem}
\proof
For every chamber $\si\supset\tau$
the function $b_{\theta_{\si}}\circ t^{-1}-b_{\theta_{\si}}$ is constant,
because $t$ fixes $\st(\tau)\supset\si$ 
and therefore $b_{\theta_{\si}}\circ t^{-1}\equiv b_{\theta_{\si}}$. 
Furthermore, 
the difference $b_{\theta_{\si}}\circ t^{-1}-b_{\theta_{\si}}$ is independent of $\si$
because, along a $t$-axis, $b_{\theta_{\si}}$ is linear with slope $-\cos\angle(\theta(\xi),\bar\theta)$
independent of $\si$.
It follows that 
$b^{\bar\theta}_{\tau,p}\circ t^{-1}-b^{\bar\theta}_{\tau,p}=b_{\theta_{\si}}\circ t^{-1}-b_{\theta_{\si}}$
is constant.
\qed

\medskip
Our next aim is to distinguish the functions $b^{\bar\theta}_{\tau,p}$ from each other. 

Let $\hat\tau$ be the simplex $p$-opposite to $\tau$, 
and let $CS(p)=CS(\tau,\hat\tau,p)$ 
denote the cross section of the parallel set $P(\tau,\hat\tau)$ through $p$. 

\begin{lem}
\label{lem:uniqmax}
$b^{\bar\theta}_{\tau,p}>0$ on $CS(p)-\{p\}$.
\end{lem}
\proof
Let $p\neq q\in CS(p)$.
We need to find a chamber $\si\supset\tau$ such that 
$$ b_{\theta_{\si}}(q) > b_{\theta_{\si}}(p)  .$$
The latter holds if 
\begin{equation*}
\angle_p(q,\theta_{\si})>\pihalf 
\end{equation*}
because then the convex function $b_{\theta_{\si}}$ strictly increases along $pq$.

Let $F\subset P(\tau,\hat\tau)$ be a maximal flat containing $pq$.
Then $\tau\subset\geo F$,
and we denote by $\eta\in\geo F$ the ideal point with $q\in p\eta$.
We will show that $\geo F$
contains a chamber $\si\supset\tau$ with the desired property,
equivalently, 
with the property that 
\begin{equation}
\label{eq:angineqds}
\tangle(\eta,\theta_{\si})>\pihalf .
\end{equation}
Let $\bar\theta'\in\interior(\taumod)$ be the nearest point projection of $\bar\theta\in\interior(\simod)$
to $\taumod$.
We denote by $\theta'_{\tau}\in\tau$ the point of type $\bar\theta'$.
The arcs $\theta'_{\tau}\theta_{\si}$ in $\geo F$ for $\si\supset\tau$ 
are perpendicular to $\tau$.
Note that $\tau\subset\D B(\eta,\pihalf)$, because $\eta\in\geo CS(p)$,
and hence also the arc $\theta'_{\tau}\eta$ is perpendicular to $\tau$.
Property (\ref{eq:angineqds}) is equivalent to 
\begin{equation}
\label{eq:angineqdsrf}
\angle_{\theta'_{\tau}}(\eta,\theta_{\si})>\pihalf . 
\end{equation}
Since the type $\bar\theta$ is regular,
the directions $\oa{\theta'_{\tau}\theta_{\si}}$ for the chambers $\si\supset\tau$ correspond to a regular Weyl orbit 
in the spherical Coxeter complex associated to the link of $\tau$ in $\tits X$.
If (\ref{eq:angineqdsrf}) would fail for all chambers $\si\supset\tau$,
then this regular Weyl orbit would be contained in a closed hemisphere,
which is impossible.\footnote{In a spherical Coxeter complex without sphere factor
no regular Weyl orbit is contained in a closed hemisphere, compare Lemma~\ref{lem:pos}.}
\qed

\medskip
Based on these properties, we can now distinguish the functions $b^{\bar\theta}_{\tau,p}$ from each other.

\begin{lem}[Distinguishing horofunctions]
\label{lem:distlifct}
$b^{\bar\theta}_{\tau,p}\equiv b^{\bar\theta}_{\tau',p'}$ 
iff $\tau=\tau'$ 
and the sectors $V(p,\tau)$ and $V(p',\tau)$ are strongly asymptotic.
\footnote{Note that two asymptotic Weyl sectors $V(x, \tau)$ and $V(x', \tau)$ are strongly asymptotic
iff their tips $x,x'$ lie in the same orbit of the closed subgroup of $P_{\tau}$ 
which is generated by $H_{\tau}$
and the transvections along lines asymptotic to $\tau$.}
\end{lem}
\proof ``$\Rightarrow$'':
We first show that $\tau$ can be read off the asymptotics of $b^{\bar\theta}_{\tau,p}$.

For a chamber $\si$ and an ideal point $\xi\in\geo X$ it holds that 
$\tangle(\xi,\theta_{\si})=\angle(\theta(\xi),\bar\theta)$ iff $\xi\in\si$,
because $\bar\theta\in\inte(\simod)$.
Lemma~\ref{lem:slpbd} therefore implies that 
along a ray $x\xi$ the function $b_{\theta_{\si}}$ has  slope $\equiv-\cos\angle(\theta(\xi),\bar\theta)$
if $\xi\in\si$,
and slope $>-\cos\angle(\theta(\xi),\bar\theta)$ everywhere if $\xi\not\in\si$.
It follows that 
$b^{\bar\theta}_{\tau,p}$ has slope $\equiv-\cos\angle(\theta(\xi),\bar\theta)$ along $x\xi$ 
iff this is the case for the $b_{\theta_{\si}}$ for all chambers $\si\supset\tau$,
i.e.\ iff $\xi\in\cap_{\si\supset\tau}\si=\tau$.
Hence, 
$b^{\bar\theta}_{\tau,p}\equiv b^{\bar\theta}_{\tau',p'}$ implies that $\tau=\tau'$.

Due to $H_{\tau}$-invariance,
$b^{\bar\theta}_{\tau,hp}=b^{\bar\theta}_{\tau,p}$ for $h\in H_{\tau}$,
we may replace $p,p'$ by points in their $H_{\tau}$-or\-bits.
We can therefore assume that $p,p'\in P=P(\hat\tau,\tau)$ with $\hat\tau\in C(\tau)$,
because every $H_{\tau}$-orbit intersects the parallel set $P$ (exactly once).
In view of Lemma~\ref{lem:transvtaubeh},
we may furthermore replace $p,p'$ by their images under transvections 
along lines parallel to the euclidean factor of $P$,
and can thus assume that they lie in the same cross section of $P$,
$p'\in CS(p)$.
Now Lemma~\ref{lem:uniqmax} implies that $p=p'$. 

``$\Leftarrow$'':
This follows from the invariance of the equivalence classes $[b^{\bar\theta}_{\tau,\cdot}]$
under $H_{\tau}$ and the transvections along lines asymptotic to $\tau$.
\qed

\begin{cor}
The points in $\geo^{\bar\theta} X$
one-to-one correspond to the strong asymptote classes of Weyl sectors in $X$.
\end{cor}

Note that all Finsler boundary points are limits of sequences along Weyl sectors,
and in particular limits of sequences along Finsler geodesic rays.
Hence all horofunctions are {\em Busemann functions}, as defined in section~\ref{sec:horob}.

\subsubsection{Convergence at infinity}
\label{sec:idpts}

We fix a base point $o\in X$
and denote by $K$ the maximal compact subgroup of $G$ fixing $o$.

We first study the convergence at infinity of divergent sequences in $X$.
Since a divergent sequence always contains pure subsequences of some face type, 
we can restrict to this case.

Let $(x_n)$ be a $\taumod$-pure sequence in $X$.
There exists a sequence $(\tau_n)$ in $\Flagt$ such that $x_n\in V(o,\st(\tau_n))$. 
Let $\hat\tau_n\in\Flagit$ denote the simplices $o$-opposite to $\tau_n$.
Then $V(o,\st(\tau_n))\subset P(\hat\tau_n,\tau_n)$. 
Due to pureness,
there exists a bounded sequence of points $p_n\in CS(\hat\tau_n,\tau_n,o)$ such that 
$x_n\in f(\hat\tau_n,\tau_n,p_n)$. 

\begin{prop}[Convergence at infinity]
\label{prop:convinfi}
Then $(x_n)$ converges in $\ol X^{\bar\theta}$ 
iff there is convergence $\tau_n\to\tau$ in $\Flagt$ and $p_n\to p\in CS(\hat\tau,\tau,o)$ in $X$.
In this case, 
$$x_n\to [b^{\bar\theta}_{\tau,p}] .$$
\end{prop}
\proof
Suppose that $\tau_n\to\tau$ and $p_n\to p$.
We write 
$\tau_n=k_n\tau$ and $\hat\tau_n=k_n\hat\tau$ 
with $k_n\to e$ in $K$. 
The sequence of points 
$$k_n^{-1}x_n\in f(\hat\tau,\tau,\underbrace{k_n^{-1}p_n}_{\to p})\subset P(\hat\tau,\tau)$$
is also $\taumod$-pure 
and contained in a tubular neighborhood of the sector $V(o,\tau)$.
It follows that 
$d(k_n^{-1}x_n,V(p,\tau))\to0$.
Lemma~\ref{lem:limalsct} then implies that 
$$d^{\bar\theta}_{k_n^{-1}x_n}-d^{\bar\theta}_{k_n^{-1}x_n}(p)
\to b^{\bar\theta}_{\tau,p}$$
uniformly on compacts,
and furthermore that 
$$d^{\bar\theta}_{x_n}-d^{\bar\theta}_{x_n}(p)\to b^{\bar\theta}_{\tau,p} ,$$
i.e.\ $x_n\to [b^{\bar\theta}_{\tau,p}]$. 
The converse direction follows from this direction, Lemma~\ref{lem:distlifct} distinguishing horofunctions
and the compactness of flag manifolds. 
\qed

\medskip
Finsler and flag convergence for divergent sequences in $X$ are related as follows:
\begin{cor}[Finsler and flag convergence]
\label{cor:finsflgcnv}
If $x_n\to [b^{\bar\theta}_{\tau,p}]$ with $\theta(\tau)=\taumod$,
then $(x_n)$ is $\taumod$-pure and $\taumod$-flag converges, $x_n\to\tau$.
\end{cor}
\proof
If $(x_n)$ were not $\taumod$-pure, we could extract a $\numod$-pure subsequence 
for a different face type $\numod\neq\taumod$.
By the lemma, after further extraction,
$(x_n)$ Finsler converges to a boundary point $[b^{\bar\theta}_{\nu,q}]$ with $\theta(\nu)=\numod$.
However, $[b^{\bar\theta}_{\nu,q}]\neq[b^{\bar\theta}_{\tau,p}]$ according to Lemma~\ref{lem:distlifct}, 
a contradiction. 
Hence $(x_n)$ must be $\taumod$-pure.
Since $x_n\in V(o,\st(\tau_n))$ and $\tau_n\to\tau$, again due to the lemma and Lemma~\ref{lem:distlifct},
the definition of flag convergence implies that $x_n\to\tau$.
\qed

\medskip
We will also use the following fact:

\begin{lem}
\label{lem:limsamstr}
Let $(x_n)$ and $(x'_n)$ be sequences in $X$ 
which are bounded distance apart
and converge at infinity,
$x_n\to [b^{\bar\theta}_{\tau,p}]$ and $x'_n\to [b^{\bar\theta}_{\tau',p'}]$.
Then $\tau=\tau'$.
\end{lem}
\proof
Since the functions $d^{\bar\theta}_{x_n}-d^{\bar\theta}_{x'_n}$ 
are uniformly bounded independently of $n$,
also $b^{\bar\theta}_{\tau,p}-b^{\bar\theta}_{\tau',p'}$ is bounded.
This implies that $\tau=\tau'$,
compare the first part of the proof of Lemma~\ref{lem:distlifct}.
\qed

\begin{rem}
If $\theta(\tau)\neq\simod$,
then the limit points 
$[b^{\bar\theta}_{\tau,p}]$ and $[b^{\bar\theta}_{\tau',p'}]$ are in general different. 
\end{rem}

We now discuss
the convergence of sequences {\em at} infinity. 

For face types $\taumod\subset\numod$,
every boundary point of type $\numod$ is a limit of boundary points of type $\taumod$:

\begin{lem}
\label{lem:canapp}
Let $\tau\subset\nu$ be faces in $\geo X$,
and let $(x_n)$ be a $\theta(\nu)$-regular sequence in a sector $V(p,\nu)$.
Then 
$[b^{\bar\theta}_{\tau,x_n}]\to[b^{\bar\theta}_{\nu,p}]$.
\end{lem}
\proof
Using Lemma~\ref{lem:limalsct},
we can approximate the boundary points $[b^{\bar\theta}_{\tau,x_n}]$ by points in $X$:
There exist points $y_n\in V(x_n,\tau)\subset V(p,\nu)$ 
such that 
$$d^{\bar\theta}_{y_n}-d^{\bar\theta}_{y_n}(x_n)-b^{\bar\theta}_{\tau,x_n}\to0$$
uniformly on compacts.
The sequence $(y_n)$ is also $\theta(\nu)$-regular,
and the same lemma yields that 
$$d^{\bar\theta}_{y_n}-d^{\bar\theta}_{y_n}(p)\to b^{\bar\theta}_{\nu,p} $$
uniformly on compacts. 
It follows that $[b^{\bar\theta}_{\tau,x_n}]\to[b^{\bar\theta}_{\nu,p}]$. 
\qed

\medskip
The next result yields necessary conditions for the convergence of sequences at infinity:
\begin{lem}
\label{lem:inbigrstr}
If
$$[b^{\bar\theta}_{\tau_n,x_n}]\to[b^{\bar\theta}_{\nu,p}]$$
and $\theta(\tau_n)=\taumod$ for all $n$, 
then $\taumod\subseteq\theta(\nu)$
and $\tau_n\to\tau\subseteq\nu$. 
\end{lem}
\proof
We may assume without loss of generality that 
$x_n\in CS(\hat\tau_n,\tau_n,o)$ and $p\in CS(\hat\nu,\nu,o)$
where $\hat\tau_n$ is $o$-opposite to $\tau_n$ and $\hat\nu$ is $o$-opposite to $\nu$.

As in the proof of the previous lemma,
we approximate using Lemma~\ref{lem:limalsct} the points $[b^{\bar\theta}_{\tau_n,x_n}]$ at infinity 
by points $y_n\in V(x_n,\tau_n)$ such that still
$$y_n\to[b^{\bar\theta}_{\nu,p}].$$
The latter holds if the growth 
\begin{equation*}
d(y_n,V(x_n,\D\tau_n)) \to+\infty
\end{equation*}
is sufficiently fast.
Sufficiently fast growth implies moreover that 
$y_n\in V(o,\st(\tau_n))$,
and hence that there exist chambers $\si_n\supseteq\tau_n$ such that 
$y_n\in V(o,\si_n)$.

After extraction,
we may assume that $(y_n)$ is $\taumod'$-pure for some face type $\taumod'\subseteq\simod$.
Invoking sufficiently fast growth 
again, it follows that $\taumod'\supseteq\taumod$.

Consider the faces $\tau_n\subseteq\tau'_n\subseteq\si_n$ of type $\theta(\tau'_n)=\taumod'$,
and denote by $\hat\tau'_n$ the simplices $o$-opposite to $\tau'_n$.
There exists a bounded sequence $(x'_n)$ of points $x'_n\in CS(\hat\tau'_n,\tau'_n,o)$
such that 
$y_n\in f(\hat\tau'_n,\tau'_n,x'_n)$.
After further extraction,
we may assume convergence $\tau'_n\to\tau'$ and $x'_n\to x'$. 
Then $y_n\to[b^{\bar\theta}_{\tau',x'}]$
by Proposition~\ref{prop:convinfi}, 
and hence $[b^{\bar\theta}_{\tau',x'}]=[b^{\bar\theta}_{\nu,p}]$.
In particular, $\taumod'=\theta(\nu)$ and $\tau'=\nu$.
It follows that $\tau_n\to\tau\subseteq\nu$,
i.e. the assertion holds for the subsequence. 

Returning to the original sequence of points $[b^{\bar\theta}_{\tau_n,x_n}]$, 
our argument shows
that every subsequence has a subsequence for which the assertion holds.
Consequently, 
$\taumod\subseteq\theta(\nu)$ 
and the sequence of simplices $\tau_n$ can only accumulate at the face $\tau\subseteq\nu$ of type $\taumod$.
In view of the compactness of $\Flagt$,
it follows that $\tau_n\to\tau$.
\qed

\medskip
Our discussion of sequential convergence implies that the Finsler compactification 
does not depend on the regular type $\bar\theta$.

\begin{prop}[Type independence of Finsler compactification]
\label{prop:independence}
For any two regular types $\bar\theta,\bar\theta'\in\interior(\simod)$,
the identity map $\id_X$ extends to a $G$-equivariant homeomorphism 
$$\ol X^{\bar\theta} \to \ol X^{\bar\theta'}$$
sending 
$[b^{\bar\theta}_{\tau,p}]\mapsto[b^{\bar\theta'}_{\tau,p}]$ at infinity.
\end{prop}
\proof
The extension of $\id_X$ sending $[b^{\bar\theta}_{\tau,p}]\mapsto[b^{\bar\theta'}_{\tau,p}]$ 
is a $G$-equivariant bijection $\ol X^{\bar\theta} \to \ol X^{\bar\theta'}$.
The conditions given in Proposition~\ref{prop:convinfi} for sequences $x_n\to\infty$ in $X$ to converge at infinity
do not depend on the type $\bar\theta$,
i.e.\ $x_n\to[b^{\bar\theta}_{\tau,p}]$ in $\ol X^{\bar\theta}$ 
iff $x_n\to[b^{\bar\theta'}_{\tau,p}]$ in $\ol X^{\bar\theta'}$.
A general point set topology argument now implies that the extension is a homeomorphism,
see Lemma~\ref{lem:topologylemma}.
\qed

\medskip
We therefore will from now on mostly use the notation $\ol X^{Fins}$ for $\ol X^{\bar\theta}$.

\subsubsection{Stratification and $G$-action}\label{sec:Stratification}
\label{sec:strat}

For every face type $\taumod\subseteq\simod$,
we define the {\em stratum} at infinity
\begin{equation}
\label{eq:strat} 
S_{\taumod}= \{[b^{\bar\theta}_{\tau,p}] : \theta(\tau)=\taumod,  p\in X\} .
\end{equation}
Furthermore, we put $S_{\emptyset}=X$. 
We define the {\em stratification} of $\ol X^{\bar\theta}$ as 
$$ \ol X^{\bar\theta}=\bigsqcup_{\emptyset\subseteq\taumod\subseteq\simod}S_{\taumod} .$$
In the sequel, when talking about the stratification, we will also admit $\emptyset$ as a face type.

Lemmas~\ref{lem:canapp} and \ref{lem:inbigrstr} yield for the closures of strata:
\begin{equation}
\label{eq:strat-sm}
\ol S_{\taumod} = \bigsqcup_{\numod\supseteq\taumod} S_{\numod}
\end{equation}
The stratum $S_{\emptyset}=X$ is open dense,
while the stratum $S_{\simod}=\DF X$ is closed 
and contained in the closure of every other stratum. 

\medskip
The continuous extension of the {\em $G$-action} on $X$ to $\ol X^{\bar\theta}$
is given at infinity by
$$ g\cdot[b^{\bar\theta}_{\tau,p}] = [b^{\bar\theta}_{\tau,p}\circ g^{-1}] = [b^{\bar\theta}_{g\tau,gp}].$$
The $G$-orbits are precisely the strata $S_{\taumod}$.

The stabilizer of a boundary point $[b^{\bar\theta}_{\tau,p}]$ is the semidirect product 
$$ H_{\tau}  \rtimes \bigl(T(\hat\tau,\tau)\times K_{f(\hat\tau,\tau,p)}\bigr) $$
where $H_{\tau}\subset P_{\tau}$ is the horocyclic subgroup,
$\hat\tau\in C(\tau)$ a simplex opposite to $\tau$,
$f(\hat\tau,\tau,p)$ the singular flat through $p$ with visual boundary sphere $s(\hat\tau,\tau)$,
$K_{f(\hat\tau,\tau,p)}<G$ its pointwise stabilizer,
and $T(\hat\tau,\tau)$ the group of transvections along $f(\hat\tau,\tau,p)$,
cf.\ section~\ref{sec:mix}.

\medskip
We will use the following observation 
concerning the dynamics of $G\acts\ol X^{\bar\theta}$:
\begin{lem}
\label{lem:zoom}
Every open subset $O\subset\ol X^{\bar\theta}$ which intersects the closed stratum,
$O\cap\DF X\neq\emptyset$, 
sweeps out the entire space,
$G\cdot O=\ol X^{\bar\theta}$.
\end{lem}
\proof
The $G$-orbit $\DF X$ is in the closure of every $G$-orbit $S_{\taumod}$.
\qed

\medskip
For the strata at infinity,
there are the natural $G$-equivariant {\em fibrations} 
of homogeneous $G$-spaces
\begin{equation}
\label{eq:fibstr}
S_{\taumod} \lra \Flagt
\end{equation}
by the forgetful maps $[b^{\bar\theta}_{\tau,p}]\mapsto\tau$.
The fiber
\begin{equation}
\label{eq:smallstrat} 
X_{\tau}= \{[b^{\bar\theta}_{\tau,p}] : p\in X \} 
\end{equation}
over $\tau\in\Flagt$ is naturally identified 
with the space of strong asymptote classes of Weyl sectors $V(x,\tau)$,
cf. Lemma~\ref{lem:distlifct},
which is in turn naturally identified with the cross section of the parallel set $P(\tau,\hat\tau)$
for any simplex $\hat\tau\in C(\tau)$.
We will refer to the fibers $X_{\tau}$ as {\em small strata}.
Again according to Lemmas~\ref{lem:canapp} and \ref{lem:inbigrstr},
we have that:
\begin{equation}
\label{eq:stratf}
\ol X_{\tau} = \bigsqcup_{\nu\supseteq\tau} X_{\nu}
\end{equation}
Note that for different simplices 
$\tau_1,\tau_2$ of the same type $\taumod$,
it holds that 
\begin{equation}
\label{eq:disjstrcls}
\ol X_{\tau_1}\cap\ol X_{\tau_2}=\emptyset 
\end{equation}
because every simplex in $\geo X$ has at most one face of type $\taumod$.

\begin{rem}
One can show that 
the closure $\ol X_{\tau}$ is naturally identified with the regular Finsler compactification of $X_{\tau}$.
\end{rem}

The discussion of convergence at infinity in the previous section 
yields the following {\em characterization of pureness and regularity}
for divergent sequences in $X$ in terms of their accumulation set in the Finsler boundary:

\begin{prop}[Pureness and regularity]
\label{prop:puregstr}
Let $x_n\to\infty$ be a divergent sequence. Then:

(i) $(x_n)$ is $\taumod$-pure iff it accumulates at a compact subset of the stratum $S_{\taumod}$.

(ii) $(x_n)$ is $\taumod$-regular iff it accumulates at the stratum closure $\ol S_{\taumod}$.
\end{prop}
\proof
(i) 
If $(x_n)$ is $\taumod$-pure,
then Proposition~\ref{prop:convinfi} implies that it accumulates at a compact subset of $S_{\taumod}$.
Otherwise,
if $(x_n)$ is not $\taumod$-pure,
then it contains a $\numod$-pure subsequence for another face type $\numod\neq\taumod$
and therefore has, by the same proposition, accumulation points in $S_{\numod}$,
i.e.\ outside $S_{\taumod}$.

(ii)
If $(x_n)$ is $\taumod$-regular, then all $\numod$-pure subsequences have type $\numod\supseteq\taumod$,
and the assertion therefore follows from (i) and (\ref{eq:strat-sm}).
\qed

\medskip
Similarly, we can characterize {\em flag convergence}, 
compare Corollary~\ref{cor:finsflgcnv} above:
\begin{prop}[Flag convergence]
\label{prop:relconvflfi}
A $\taumod$-regular sequence $(x_n)$ $\taumod$-flag converges, $x_n\to\tau\in\Flagt$,
iff it accumulates at the small stratum closure $\ol X_{\tau}$.
\end{prop}
\proof
By the previous proposition, 
$(x_n)$ accumulates at $\ol S_{\taumod}$.

Suppose that we have Finsler convergence $x_n\to[b^{\bar\theta}_{\nu,p}]$.
By Corollary~\ref{cor:finsflgcnv}, 
$(x_n)$ is $\numod$-pure
with $\numod=\theta(\nu)$ and $\numod$-flag converges, $x_n\to\nu$.
Necessarily $\numod\supseteq\taumod$, because $(x_n)$ is $\taumod$-regular.
It follows that we also have $\taumod$-flag convergence $x_n\to\tau_{\nu}\subset\nu$
to the face $\tau_{\nu}$ 
of type $\taumod$.
Furthermore, $[b^{\bar\theta}_{\nu,p}]\in\ol X_{\tau_{\nu}}$, cf.\ (\ref{eq:stratf}).

Thus,
if $x_n\to\tau$
then all accumulation points of $(x_n)$ in $\geo^{Fins}X$ must lie in $\ol X_{\tau}$.
On the other hand,
if $(x_n)$ does not $\taumod$-flag converge to $\tau$, then after extraction it $\taumod$-flag converges
to some other $\tau'\neq\tau$ 
and has Finsler accumulation points in the small stratum closure $\ol X_{\tau'}$ disjoint from $\ol X_{\tau}$,
cf.\ (\ref{eq:disjstrcls}).
\qed

\subsubsection{Maximal flats and Weyl sectors}
\label{sec:compmaxfl}

We first discuss maximal flats $F\subset X$.
We start by showing that 
their extrinsic closure in $\ol X^{\bar\theta}$ coincides with 
their intrinsic Finsler compactification, 
cf.\ the general discussion in the end of section~\ref{sec:horob}.

According to Proposition~\ref{prop:convinfi},
$$ \D^{\ol X^{\bar\theta}} F= \{[b^{\bar\theta}_{\nu,x}] : \nu\subset\geo F, x\in F\} \subset\geo^{\bar\theta}X $$
where $\D^{\ol X^{\bar\theta}} F\subset\geo^{\bar\theta}X$ means the ``boundary'' $\ol F-F$ of $F$ inside $\ol X^{\bar\theta}$.

\begin{lem}
For simplices $\nu,\nu'\subset\geo F$ and points $x,x'\in F$ it holds that:
If $b^{\bar\theta}_{\nu,x}|_F\equiv b^{\bar\theta}_{\nu',x'}|_F$,
then $b^{\bar\theta}_{\nu,x}\equiv b^{\bar\theta}_{\nu',x'}$.
\end{lem}
\proof
One proceeds as in the proof of Lemma~\ref{lem:distlifct}.
The asymptotics of $b^{\bar\theta}_{\nu,x}|_F$ and $b^{\bar\theta}_{\nu',x'}|_F$
allow to read off $\nu$ and $\nu'$,
and thus $\nu=\nu'$.
Furthermore,
Lemma~\ref{lem:uniqmax} implies that the singular flats spanned by the Weyl sectors 
$V(x,\nu)$ and $V(x',\nu)$ coincide,
equivalently, these Weyl sectors intersect.
Hence they are strongly asymptotic,
which implies that $b^{\bar\theta}_{\nu,x}\equiv b^{\bar\theta}_{\nu',x'}$. 
\qed

\medskip
It follows that there is a natural inclusion
$$ \ol F^{\bar\theta} \subset \ol X^{\bar\theta}$$
of Finsler compactifications.

The stratification of $\ol X^{\bar\theta}$ induces the {\em stratification}
$$ \ol F^{\bar\theta}= \bigsqcup_{\taumod} S^F_{\taumod}$$
by the stratum $S^F_{\emptyset}=F$ and the strata at infinity
$$ S^F_{\taumod}=S_{\taumod}\cap\geo^{\bar\theta}F 
= \{[b^{\bar\theta}_{\tau,x}] : \tau\subset\geo F,\theta(\tau)=\taumod ,x\in F\} .
$$
The fibration (\ref{eq:fibstr})
restricts to the finite decomposition 
$$ S^F_{\taumod}=\bigsqcup_{\tau\subset\geo F,\theta(\tau)=\taumod} X^F_{\tau}$$
into the small strata 
\begin{equation}\label{eq:strata_of_bar_F}
X^F_{\tau}= \{[b^{\bar\theta}_{\tau,x}] : x\in F \} 
\end{equation}
at infinity.
They are euclidean spaces 
which are canonically identified with the affine subspaces of $F$
perpendicular to the Weyl sectors $V(x,\tau)$.
By analogy with (\ref{eq:stratf}), the closures of small strata decompose as
\begin{equation}
\label{eq:stratfmfl}
\ol X^F_{\tau} = \bigsqcup_{\tau\subseteq\nu\subseteq\geo F} X^F_{\nu} .
\end{equation}
{\em Dynamics at infinity:}
The subgroup $T_F<G$ of transvections along $F$
restricts to the group of translations on $F$.
Unlike for the visual boundary,
the induced action
$$T_F\acts\geo^{\bar\theta}F$$
on the Finsler boundary is nontrivial in higher rank.
Its orbits are the small strata $X^F_{\tau}$, 
since 
$t[b^{\bar\theta}_{x,\tau}]=[b^{\bar\theta}_{tx,\tau}]$
for $t\in T_F$.
The stabilizer $\Stab_G(F)$ of $F$ in $G$ acts on $F$ by the affine Weyl group.
Its orbits at infinity are the big strata $S^F_{\taumod}$.

It is worth pointing out how the description of the {\em convergence at infinity}
simplifies for divergent sequences contained in maximal flats.
Proposition~\ref{prop:convinfi} reduces to:

\begin{prop}[Convergence at infinity for maximal flats]
\label{lem:convdescrfl}
Suppose that $(x_n)$ is a $\taumod$-pure sequence in $F$ 
which for some simplex $\tau\subset\geo F$ of type $\taumod$ 
is contained in tubular neighborhoods of the sectors $V(\cdot,\tau)$.

Then $(x_n)$ converges in $\ol F^{\bar\theta}$ iff the flats $f(\hat\tau,\tau,x_n)\subset F$ Hausdorff converge.
In this case,
$$x_n\to [b^{\bar\theta}_{\tau,p}] $$
with $p\in F$ such that 
$f(\hat\tau,\tau,x_n)\to f(\hat\tau,\tau,p)$.
\end{prop}

The discussion for {\em Weyl sectors} $V=V(p,\tau)$ is analogous:
We have 
\begin{equation}
\label{eq:fbdws}
\D^{\ol X^{\bar\theta}}V=
\{[b^{\bar\theta}_{\nu,x}] : \nu\subseteq\tau, x\in V\} \subset\geo^{\bar\theta}X .
\end{equation}
Again, horofunctions uniquely extend from $V$ to $X$:
\begin{lem}
For simplices $\nu,\nu'\subseteq\tau$ and points $x,x'\in V$ it holds that:
If $b^{\bar\theta}_{\nu,x}|_V\equiv b^{\bar\theta}_{\nu',x'}|_V$,
then $b^{\bar\theta}_{\nu,x}\equiv b^{\bar\theta}_{\nu',x'}$.
\end{lem}
\proof
As in the case of maximal flats, cf.\ the previous lemma,
one can recognize $\nu$ and $\nu'$ from the asymptotics of $b^{\bar\theta}_{\nu,x}|_V$ and $b^{\bar\theta}_{\nu',x'}|_V$,
and thus sees that $\nu=\nu'$.
Then Lemma~\ref{lem:uniqmax} implies that 
$V(x,\nu)$ and $V(x',\nu)$ intersect and hence are strongly asymptotic.
\qed

\medskip
We thus have the natural inclusion
$$ \ol V^{\bar\theta} \subset \ol X^{\bar\theta} .$$
Furthermore,
we have the {\em stratification}
\begin{equation}
\label{eq:strtsctcls}
\ol V^{\bar\theta}= \bigsqcup_{\nu\subseteq\tau} X^V_{\nu}
\end{equation}
by $X^V_{\emptyset}=V$ and the (small) strata at infinity
$$ X^V_{\nu}=X_{\nu}\cap\geo^{\bar\theta}V
= \{[b^{\bar\theta}_{\nu,x}] : x\in V\} .$$
The stratum closures decompose as
\begin{equation*}
\ol X^V_{\nu} = \bigsqcup_{\nu\subseteq\nu'\subseteq\tau} X^V_{\nu'} .
\end{equation*}

\subsubsection{Action of maximal compact subgroups}

Let $o\in X$ be a base point and $K<G$ the maximal compact subgroup fixing it.
Let $V=V(o,\si)$ be a euclidean Weyl chamber in $X$ with tip at $o$.
We recall some facts about the action $K\acts X$:

(i) 
$V$ is a {\em cross section} for the action,
i.e. every $K$-orbit intersects $V$ exactly once.

(ii) {\em Point stabilizers.}
The fixed point set in $V$ of any element $k\in K$ 
is a Weyl sector $V(o,\tau)$,
where $\emptyset\subseteq\tau\subseteq\si$ is the face fixed by $k$.
In other words, if $k$ fixes a point $p\in V$,
then it fixes the smallest Weyl sector $V(o,\tau)$ containing it. 
(Here, we put $V(o,\emptyset):=\{o\}$.)

We now establish analogous properties for the action $K\acts\ol{X}^{\bar\theta}$.

\begin{lem}[Cross section]
\label{lem:cross-section}
$\ol V^{\bar\theta}\subset\ol X^{\bar\theta}$ is a cross section for the action of $K\acts\ol{X}^{\bar\theta}$. 
\end{lem} 
\proof 
Since $K\cdot\ol V^{\bar\theta}$ is compact and contains $K\cdot V=X$,
and since $X$ is dense in $\ol{X}^{\bar\theta}$,
it holds that $K\cdot\ol V^{\bar\theta}=\ol X^{\bar\theta}$,
i.e.\ every $K$-orbit in $\ol{X}^{\bar\theta}$ intersects $\ol V^{\bar\theta}$.
We must show that every $K$-orbit in $\geo^{\bar\theta} X$ intersects $\geo^{\bar\theta} V$ only once. 
Suppose that 
\begin{equation}
\label{eq:ktransp}
k\cdot [b^{\bar\theta}_{\tau,p}] = [b^{\bar\theta}_{k\tau,kp}] = [b^{\bar\theta}_{\tau',p'}]
\end{equation}
for $k\in K$ and Weyl sectors $V(p,\tau),V(p',\tau')\subset V$, 
cf.\ (\ref{eq:fbdws}).
Then, 
in particular, $k\tau=\tau'$.
Since $\tau,\tau'\subseteq\si$,
this implies that $\tau=\tau'$ and $k\tau=\tau$.

It follows that $k$ preserves the parallel set $P(\hat\tau,\tau)$ 
where $\hat\tau$ denotes the simplex $o$-opposite to $\tau$.
The sectors $V(kp,\tau),V(p',\tau)\subset P(\hat\tau,\tau)$ are strongly asymptotic,
because $[b^{\bar\theta}_{\tau,kp}] = [b^{\bar\theta}_{\tau,p'}]$
This means that they intersect.
Let $q\in V(kp,\tau)\cap V(p',\tau)$.
Then $q,k^{-1}q\in V$ and hence $k^{-1}q=q$,
because $V$ is a cross section for the action $K\acts X$.
It follows that the sectors $V(p,\tau)$ and $V(kp,\tau)$ intersect
and hence are strongly asymptotic. 
Thus, $k\cdot [b^{\bar\theta}_{\tau,p}] =[b^{\bar\theta}_{\tau,kp}] = [b^{\bar\theta}_{\tau,p}]$.
\qed

\begin{lem}
\label{lem:stabs}
Let $k\in K$ and $V(p,\tau)\subset V$. 
The following are equivalent:

(i) $k$ fixes $[b^{\bar\theta}_{\tau,p}]\in\ol V^{\bar\theta}$.

(ii) $k$ fixes $V(p,\tau)$ pointwise.

(iii) $k$ fixes pointwise the smallest Weyl sector $V(o,\nu)$ containing $V(p,\tau)$. 
\end{lem}
\proof (i)$\Ra$(ii):
In the proof of the previous lemma,
we saw that the sectors $V(p,\tau)$ and $V(kp,\tau)$ intersect.
Since $k$ also preserves the cross sections of $P(\hat\tau,\tau)$,
it follows that $k$ fixes $p$.
The converse direction (ii)$\Ra$(i) is trivial. 

(ii)$\Leftrightarrow$(iii):
This is clear, because 
the fixed point set of $k$ on $V$ is a sector $V(o,\nu)$.
\qed

\medskip
Let $K_{\tau}=P_{\tau}\cap K$ denote the stabilizer of the simplex $\tau$ in $K$,
and put $K_{\emptyset}=K$.

\begin{cor}[Point stabilizers in compactified euclidean Weyl chambers]
\label{cor:stabs}
(i) 
The points in 
$\ol V^{\bar\theta}$ fixed by $K_{\tau}$ are 
precisely the points in $\ol{V(o,\tau)}^{\bar\theta}$.

(ii) The points with stabilizer equal to $K_{\tau}$
are precisely the points in 
$$\ol{V(o,\tau)}^{\bar\theta}-\bigcup_{\emptyset\subseteq\nu\subsetneq\tau}\ol{V(o,\nu)}^{\bar\theta}.$$
\end{cor}

\begin{notation}
In view of Proposition~\ref{prop:independence} we will from now on 
denote the Finsler compactification $\ol X^{\bar\theta}$
for $\bar\theta\in\interior(\simod)$ by $\ol X^{Fins}$. 
\end{notation}

\section{Coxeter groups and their regular polytopes}
\label{sec:cox}

\subsection{Basics of polytopes}

We refer the readers to \cite{Grunbaum} and \cite{Ziegler} for a detailed treatment of polytopes. 
In what follows, $V$ will denote a euclidean vector space, i.e. a finite-dimensional real vector space equipped with an inner product $(x,y)$. We will use the notation $V^*$ for the dual vector space,
and for $\la\in V^*$ and $x\in V$ we let $\<\la, x\>= \la(x)$. The inner product on $V$ defines 
the inner product, again denoted $(\la, \mu)$, on the dual space.

A {\em polytope} $B$ in $V$ is a compact convex subset equal to the intersection of finitely many closed half-spaces. Note that we do not require $B$ to have nonempty interior. 
The {\em affine span} $\<B\>$ of $B$ is the minimal affine subspace of $V$ containing $B$. 
The topological frontier of $B$ in its affine span is the boundary $\partial B$ of $B$.  
A {\em facet} of $B$ is a codimension one face of $\partial B$. 

Each polytope $B$ has a {\em face poset} ${\mathcal F}_B$.
It is the poset  whose elements are the faces of $B$ 
with the order given by the inclusion relation. Two polytopes are {\em combinatorially isomorphic} if there is an isomorphism of their posets. Such an isomorphism necessarily preserves the dimension of faces. Two polytopes $B$ and $B'$
 are {\em combinatorially homeomorphic} if there exists a (piecewise linear) homeomorphism $h: B\to B'$ which sends faces to faces. 
 
Given a polytope $B$ whose dimension equals $n=\dim(V)$, the {\em polar} (or {\em dual}) polytope of $B$ is defined as the following subset of the dual vector space: 
$$
B^*=\{\la\in V^*: \la(x)\le 1, \forall x\in B\}. 
$$
Thus, $\la\in B^*\subset V^*$ implies that the affine hyperplane $H_\la=\{\la=1\}$ is disjoint from the interior of $B$.  Moreover, $\la\in \partial B^*$ iff $H_\la$ has nonempty intersection with $B$. 
Each face $\varphi$ of $B$ determines the {\em dual face} $\varphi^*$ of $B^*$, consisting of the elements 
$\la\in B^*$ which are equal to $1$ on the entire face $\varphi$. This defines a natural bijection between the faces of $B$ and $B^*$:
$$
\star: \varphi \mapsto \varphi^*. 
$$
Under this bijection, faces have complementary dimensions:
$$
\dim(\varphi)+ \dim(\varphi^*)= n-1.  
$$
The bijection $\star$ also reverses the face inclusion:
$$
\varphi\subset \psi \iff \varphi^*\supset \psi^*. 
$$
In particular, the face poset of $\partial B^*$ is dual to the face poset of $\partial B$. If $W$ is a group of linear transformations preserving $B$, its dual action
$$
w^*(\la) = \la\circ w^{-1}
$$
on $V^*$  preserves $B^*$. 
The naturality of $\star$ implies that it is $W$-equivariant. 

A polytope $B$ is called  {\em simplicial} if its faces are simplices. 
It is called {\em simple} if it has a natural structure  of a manifold with corners: Each vertex $v$ of $B$ is contained in exactly $d$ facets, where $d$ is the dimension of $B$. Equivalently, the affine functionals defining these facets in $\<B\>$ 
have linearly independent linear parts. For each simplicial polytope, its dual is a simple polytope,
and vice versa.

\begin{lem}\label{lem:combinatorialequivalence}
Two polytopes are combinatorially isomorphic if and only if they are combinatorially homeomorphic. 
\end{lem} 
\proof One direction is trivial.
Conversely, given an isomorphism of posets,
one constructs a combinatorial homeomorphism by induction over skeleta and coning off.
\qed

\subsection{Root systems}

In this and the following sections, the euclidean vector space $V$ is identified with the model maximal flat $\Fmod$ 
for the symmetric space $X$; the root system $R\subset V^*$ is the root system of $X$. Accordingly, the Coxeter group $W$ defined via $R$ is the Weyl group of $X$. Since the symmetric space $X$ has noncompact type, $R$ spans $V^*$, 
i.e. $W$ fixes only the origin $0$ in $V$. 

Given a face $\tau$ of the spherical Coxeter complex $\geo V$, 
we define the root subsystem 
$$R_\tau\subset R$$
 consisting of all roots which vanish identically on $V(0,\tau)$.  

Each root $\al\in R$  corresponds to a {\em coroot} $\al^\vee\in V$, which is a vector such that the reflection 
$s_\al: V\to V$ corresponding to $\al$ acts on $V$ by the formula:
\begin{equation}\label{eq:reflection}
s_\al(x)= x - \<\al, x\> \al^\vee. 
\end{equation}
The group $W$ also acts isometrically on the dual space $V^*$; each reflection $s_\al\in W$ acts on $V^*$ as a reflection.
The corresponding wall is given by the equation
$$
\{\la\in V^*: \<\la, \al^\vee\>=0\}, 
$$
equivalently, this wall is $\al^\perp$, the orthogonal complement of $\al$ in $V^*$.

\medskip 
From now on, we fix a Weyl chamber $\Delta=\De_{mod}\subset V$ for the action of $W$ on $V$. 
The visual boundary of $\Delta$ is the model spherical chamber $\simod$.

\begin{notation}
We let $[n]$ denote the set $\{1,...,n\}$. 
\end{notation}

The choice of $\De$ determines the set of positive roots $R^+\subset R$ and the set of {\em simple roots} $\al_1,...,\al_n\in R^+$, where $n=\dim(V)$; 
$$
\Delta=\{x\in V: \al_i(x)=\<\al_i, x\> \ge 0, i\in [n]\}. 
$$
We will use the notation $s_i=s_{\al_i}$ for the {\em simple reflections}.
They generate $W$. 

The {\em dual chamber} to $\Delta$ is
$$
\Delta^*\subset V^*, \quad \Delta^*= \{\la\in V^*: (\al_i, \la)\ge 0, i\in [n]\}. 
$$

\begin{rem}
Note that there is another notion of a dual cone to $\Delta$ in $V^*$, 
namely the {\em root cone} $\Delta^\vee$, consisting of all $\la\in V^*$ such that the restriction of $\la$ to $\Delta$ is nonnegative. The root cone 
consists of the nonnegative linear combinations of simple roots. The root cone contains the dual chamber but, is, with few exceptions, strictly larger. 
\end{rem}

Let $B$ be a $W$-invariant polytope in $V$ with nonempty interior. We will use the notation $\Delta_B=\Delta\cap B$, $\Delta^*_{B^*}=\Delta^*\cap B^*$.

\begin{lem}\label{lem:chamber}
Suppose that $\la\in \Delta^*$ is such that $\la(x)\le 1$ for all $x\in\Delta_B$. Then $\la\in B^*$. 
\end{lem}
\proof Let $\la\in V^*$ and let $v\in\interior(\De)\subset V$.
Then $\la|_{Wv}$ is maximal in $v$ iff $\la\in\De^*$.
The assertion follows.
\qed

\subsection{Geometry of the dual ball}

We assume now that $B\subset V$ is a $W$-invariant polytope in $V$ with nonempty interior, 
such that  
$$\Delta_B=\{x\in\De: l(x)\leq1\}$$ 
where $l=l_{\bar\theta} \in \interior(\De^*)$ is a {\em regular linear functional}. 
The gradient vector of $l$ gives a direction $\bar\theta$, which is a regular point of $\simod$.

Set $l_w=w^*l=l\circ w^{-1}$, where $w\in W$. Then, 
$$
B= \bigcap_{w\in W} \{x\in V: l_w(x)  \le 1\}, 
$$
i.e. the facets of $B$ are carried by the affine hyperplanes $l_w= 1$ for $w\in W$. 

The polytope $B$ defines a (possibly non-symmetric) norm on $V$, 
namely the norm for which $B$ is the unit ball:
\begin{equation}\label{eq:norm}
||x||= ||x||_{\bar\theta}= \max_{w\in W} \left( l_w(x)\right). 
\end{equation}

\medskip
We let $\om_1,...,\om_n$ denote the nonzero vertices of the $n$-simplex $\Delta_B$. We will label these vertices consistently with the labeling of the simple roots: $\om_i$ is the unique vertex of $\Delta_B$ on which 
$\al_i$ does not vanish. Geometrically speaking, $\om_i$ is opposite to the facet $A_i$ of $\Delta_B$ carried by the wall $\al_i=0$. 

The regularity of $l$ implies: 
\begin{lem}\label{lem:simplicial}
The polytope $B$ is simplicial. Its facets are the simplices
$$
\{ x\in w\Delta : l_w(x)= 1 \}.
$$ 
For each reflection $s_i=s_{\al_i}$, the line segment $\om_i s_i(\om_i)$ is not contained in the boundary of $B$. 
\end{lem}
\proof We will prove the last statement. The segment $\om_i s_i(\om_i)$ is parallel to the vector $\al_i^\vee$. If 
$\al_i^\vee$ were to be parallel to the face $l=1$ of $B$, then $\<l, \al_i^\vee\>=0$, which implies that $l$ is singular. \qed

\begin{cor} 
Since the polytope $B$ is  simplicial, the dual polytope $B^*$ is simple. 
\end{cor}

The chamber $\Delta^*$ contains a distinguished vertex of $\Delta^*_{B^*}$, namely the linear functional 
$l$; this is the only vertex of $\Delta^*_{B^*}$ contained in the interior of $\Delta^*$. 
(The other vertices of $\Delta^*_{B^*}$ are not vertices of $B^*$.)

We now analyze the geometry of $\Delta^*_{B^*}$ in more detail.

\begin{lem}
$\Delta^*_{B^*}$ is given by the set of $2n$ inequalities $(\cdot, \al_i)\ge 0$ and $\<\cdot, \om_i\>\le 1$ for $i\in [n]$. 
\end{lem}
\proof It is clear that these inequalities are necessary for $\la\in V^*$ to belong to $\Delta^*_{B^*}$. In order to prove that they are sufficient, we have to show that each $\la$ satisfying these inequalities belongs to $B^*$. The inequalities $\<\la, \om_i\>\le 1$ show that the restriction of $\la$ to $\Delta_B$ is $\le 1$. Now, Lemma \ref{lem:chamber} shows that $\la(x)\le 1$ for all $x\in B$. \qed  

\medskip
Close to the origin, 
$\Delta^*_{B^*}$ is given by the $n$ inequalities $(\cdot, \al_i)\ge 0$,
while the other $n$ inequalities are strict.
Close to $l$, it is given by the $n$ inequalities $\<\cdot, \om_i\>\le 1$,
while the other $n$ inequalities are strict.

\medskip
We define the {\em exterior facet} $E_i\subset\De^*_{B^*}$ by the equation
$$ \langle\cdot,\om_i\rangle =1 ,$$
and the {\em interior facet} $F_j\subset\De^*_{B^*}$ as the fixed point set of the reflection $s_j$,
equivalently, by the equation 
$$ (\cdot,\al_j)=0 .$$
For subsets $I, J\subset [n]=\{1,...,n\}$ we define 
the {\em exterior faces}
$$ E_I:= \bigcap_{i\in I} E_i $$
containing $l$,
and the {\em interior faces}
$$ F_J:= \bigcap_{j\in J} F_j $$
containing the origin.
These are nonempty faces of $\Delta^*_{B^*}$ of the expected dimensions,
due to the linear independence of the $\om_i$'s, respectively, the $\al_j$'s. 

As a consequence of the last lemma, 
every face of the polytope $\Delta^*_{B^*}$ has the form
$$ E_I \cap F_J $$ 
for some subsets $\emptyset\subseteq I,J\subseteq[n]$.

\medskip
We now describe the combinatorics of the polytope $\Delta^*_{B^*}$.

\begin{lem}
For each $i=1,...,n$, $E_i\cap F_i=\emptyset$. 
\end{lem}
\proof Suppose that $\la\in \Delta^*_{B^*}$ 
is a point of intersection of these faces. Then $\la$ is a linear function fixed by the reflection $s_i$ and satisfying the equation
$\< \la, \om_i\> =1$. 
Then $\la(s_i(\om_i))=1$ as well. Thus, $\la=1$ on the entire segment connecting the vertices $\om_i$ and $s_i(\om_i)$ of $B$. Since $\la$ belongs to $B^*$, this segment has to be contained in the boundary of $B$. But this contradicts Lemma \ref{lem:simplicial}. Therefore, such a $\la$ cannot exist. \qed 

\medskip
We denote by $W_J<W$ the subgroup generated by the reflections $s_j$ for $j\in J$.
The fixed point set of $W_J$ on $\De^*_{B^*}$ equals $F_J$.

Furthermore, 
we define $\om_I$ as the face of $B$, as well as of $\De_B$,  
which is the convex hull of the vertices $\om_i$ for $i\in I$.
The dual face $\om_I^*$ of $B^*$ is given, as a subset of $B^*$, 
by the equations $\langle\cdot,\om_i\rangle=1$.
It is contained in $W_J\cdot\De^*_{B^*}$, where we put $J=[n]-I$. 
Indeed, 
the vertices of $\om_I^*$ are the functionals $l_w$
for which the dual facet $l_w=1$ of $B$ contains $\om_I$,
equivalently, for $w\in W_J$.

Note that $W_J$ preserves $\om_I$ and therefore also $\om_I^*$ 
(and acts on it as a reflection group).
The fixed point set of $W_J$ on $W_J\cdot\De^*_{B^*}$ is contained in the intersection 
$$\bigcap_{w\in W_J}w\De^*_{B^*}$$
and in particular in $\De^*_{B^*}$.
This implies that 
$$ \emptyset\neq\Fix_{W_J}(\om_I^*) \subset\De^*_{B^*} .$$
Note that $E_I=\om_I^*\cap\De^*_{B^*}$.
It follows that 
$$ E_I\cap F_J\supseteq \Fix_{W_J}(\om_I^*) \neq\emptyset .$$
In combination with the previous lemma, we conclude:
\begin{lem}
\label{lem:nonemptints}
For arbitrary subsets $\emptyset\subseteq I, J\subseteq [n]$, it holds that 
$E_I\cap F_J\ne\emptyset$ iff $I\cap J=\emptyset$.  
\end{lem}
Next, we prove the uniqueness of the labeling of the faces.
\begin{lem}
\label{lem:contnm}
If $E_I\cap F_J=E_{I'}\cap F_{J'}\neq\emptyset$, then $I=I'$ and $J=J'$.
\end{lem}
\proof
Since $E_I\cap E_{I'}=E_{I\cup I'}$ and $F_J\cap F_{J'}=F_{J\cup J'}$,
the proof reduces to the case of containment $I\subseteq I'$ and $J\subseteq J'$. 

Suppose that $j'\in J'-J$.
Then, intersecting both sides of the equality $E_I\cap F_J=E_{I'}\cap F_{J'}$
with $E_{j'}$, the previous lemma yields that 
$$ \emptyset\neq E_{I\cup\{j'\}}\cap F_J = E_{I'\cup\{j'\}}\cap F_{J'} =\emptyset, $$
a contradiction. 
Thus $J=J'$, and similarly $I=I'$.
\qed

\medskip
For the $n$-cube $[0,1]^n$,
we define similarly facets 
$E'_i=\{t_i=1\}$ and $F'_j=\{t_j=0\}$.
They satisfy the same intersection properties as in Lemmas~\ref{lem:nonemptints} and \ref{lem:contnm}.
Hence the correspondence 
$$E_I\cap F_J \stackrel{c}{\mapsto} E'_I\cap F'_J$$
is a combinatorial isomorphism between the polytopes $\De^*_{B^*}$ and $[0,1]^n$.
Lemma~\ref{lem:combinatorialequivalence} now yields:
\begin{thm}\label{thm:cube-iso}
The polytope $\Delta^*_{B^*}$ is combinatorially homeomorphic to the $n$-cube $[0,1]^n$, 
i.e. there exists a combinatorial homeomorphism 
$$\Delta^*_{B^*} \stackrel{h}{\lra}  [0,1]^n$$
inducing the bijection $c$ of face posets. 
\end{thm}

\subsection{Cube structure of the compactified Weyl chamber}

In this section we construct a canonical homeomorphism from the Finsler compactification $\ol{\Delta}^{Fins}$ of the model Weyl chamber $\Delta\subset V$, to the cube $[0,\infty]^n$. Recall that $\al_1,...,\al_n$ are the simple roots with respect to $\Delta$. Each intersection
$$
\Delta_i= \ker(\al_i)\cap \Delta
$$
is a facet of $\Delta$. 

For $x\in \Delta$ define
$$
\oa\al(x):= (\al_1(x),\ldots, \al_n(x))\in [0,\infty)^n. 
$$
This map is clearly a homeomorphism from $\Delta$ to $[0,\infty)^n$. 
We wish to extend the map $\oa{\al}$ to a homeomorphism of the compactifications. 

We recall the description of sequential convergence at infinity, see Proposition~\ref{lem:convdescrfl}.
A sequence $x_k\to\infty$ in $\De$ converges at infinity iff the following properties hold:

(i) By parts (i) and (ii) of the lemma,
there exists a face $\tau=\taumod$ of $\simod=\geo \Delta$ such that for every $\al_i\in R_{\tau}$
the sequence $(\al_i(x_k))$ converges.

(ii) By the $\taumod$-regularity assertion in part (i) of the lemma,
for the other simple roots $\al_i\not\in R_{\tau}$, we have divergence $\al_i(x_k)\to+\infty$.

In other words,
the sequence $(x_k)$ converges at infinity, iff the limit 
$$ \lim_{k\to+\infty} \oa{\al}(x_k) \in [0,\infty]^n$$
in the closed cube exists.
Moreover, Lemma~\ref{lem:convdescrfl} combined with Lemma~\ref{lem:distlifct} 
implies that the extension 
$$ \ol\De^{Fins} \stackrel{\oa{\al}}{\lra} [0,\infty]^n $$
sending 
$$ \lim_{k\to+\infty}x_k \mapsto  \lim_{k\to+\infty} \oa{\al}(x_k)$$
for sequences $(x_k)$ converging at infinity
is well-defined and bijective.
Now, Lemma~\ref{lem:topologylemma}
implies that the extension is a homeomorphism. 
Composing with the homeomorphism
$$
\kappa: [0,\infty]^n \to [0,1]^n, \quad \kappa: (t_1,\ldots, t_n)\mapsto \left(1- \frac{1}{t_1+1}, \ldots, 1- \frac{1}{t_n+1}\right)
$$
we obtain:
\begin{lem}\label{lem:matchingfaces}
The map $\kappa\circ \oa{\al}$ is a homeomorphism from $\ol{\Delta}^{Fins}$ to the cube $[0,1]^n$. It sends the compactification of each face $\ol{\Delta}^{Fins}_i, i\in [n]$,  to the face $F'_i$ of the cube  $[0,1]^n$. 
\end{lem}

\medskip
For a partition $[n]=I\sqcup J$,
we define $\emptyset\subseteq\tau_I\subseteq\simod$ as the face fixed by the reflections $s_j$ for $j\in J$.
Equivalently,
the vertices of $\tau$ are the directions of the vectors $\om_i$ for $i\in I$. 

Vice versa,
for a face $\emptyset\subseteq\tau=\taumod\subseteq\simod$,
we define the partition 
$[n]=I_{\tau}\sqcup J_{\tau}$
such that $\tau_{I_{\tau}}=\tau$,
i.e. $I_{\tau}$ indexes the vertices of $\tau$. 

Moreover, we have the sector
$\Delta_I=\cap_{i\in I}\De_i=V(0,\tau_I)\subset \Delta$ 
and its compactification 
$$
\ol{\Delta}^{Fins}_I=   \bigcap_{i\in I} \ol{\Delta}^{Fins}_i, 
$$
cf. \eqref{eq:strtsctcls}.

Recall that our vector space $V$ is the underlying vector space of the model maximal flat $F=F_{mod}$. 
We can now combine the above lemma with the homeomorphism constructed in Theorem \ref{thm:cube-iso}:

\begin{thm}\label{thm:constructphi}
There exists a homeomorphism 
$$
\ol{\Delta}^{Fins}\stackrel{\phi}{\lra} \Delta^*_{B^*}\subset B^* 
$$
satisfying the following:

1. For each partition $[n]=I\sqcup J$,
$$\phi(\ol X_{\tau_I}^{\Delta})= E_I $$
and 
$$\phi(\ol{\Delta}^{Fins}_J)= F_J.$$
In particular, $\phi(0)=0$. 

2. The map $\phi$ preserves $W$-stabilizers: $x\in \ol{\Delta}^{Fins}$ is fixed by $w\in W$ iff $\phi(x)$ is fixed by $w$. 

3. As a consequence, 
$\phi$ extends to a $W$-equivariant homeomorphism of the compactified model flat to the dual ball: 
$$
\Phi_{\Fmod}: \ol F_{mod}^{Fins}\to B^*. 
$$
\end{thm}
\proof 
Combining Theorem \ref{thm:cube-iso} and Lemma \ref{lem:matchingfaces}, 
we define
$$
\phi= h^{-1}\circ \kappa\circ \oa{\al}. 
$$
$\De^*_{B^*}$ is a cross section for the action of $W$ on $B^*$,
because $\De^*$ is a cross section for its action on $V^*$.
By Lemma \ref{lem:cross-section}, the compactified chamber $ \ol{\Delta}^{Fins} $ is a cross section for the action of $W$ on $\ol{F}^{Fins}$. We also note that for $J=[n]- I$, the fixed point sets of the subgroup $W_J< W$ 
in $\ol{\Delta}^{Fins}$ and $\De^*_{B^*}$ are precisely 
$\ol{\Delta}^{Fins}_I$ and  $F_I$,
cf. Corollary~\ref{cor:stabs}.
The last assertion of the theorem follows using Lemma \ref{lem:transformation lemma}. \qed 

\begin{rem}
One can also derive this theorem from  \cite[Proposition I.18.11]{Borel-Ji}. 
Our proof is a direct argument which avoids symplectic geometry. 
\end{rem}

\begin{rem}
We note that the paper \cite{Karlsson} computes horofunctions on finite dimensional vector spaces $V$ equipped with polyhedral norms, but does not address the question about the global topology of the associated compactification of $V$. 
See also \cite{Brill,Walsh07}.
\end{rem}  

\begin{ques}
Suppose that $||\cdot||$ is a polyhedral norm on a finite-dimensional real vector space $V$. 
Is it true that the horoclosure $\ol{V}$ of $V$ with respect to this norm, with its natural stratification,
is homeomorphic to the closed unit ball for the dual norm?
Is it homeomorphic to a closed ball for arbitrary norms?
\end{ques}

\section{Manifold-with-corners structure on the Finsler compactification}
\label{sec:mfcr}

In this section, we assume that $\bar\theta\in\interior(\simod)$.
We recall that the Finsler compactification is independent of the choice of $\bar\theta$.

In  Theorem \ref{thm:constructphi} we proved the existence of a $W$-equivariant homeomorphism 
$\Phi_F: \ol{F}^{Fins}\to B^*$.
Since $B^*$ is a simple polytope, it has a natural structure of a manifold with corners, whose strata are the faces of $B^*$. Via the homeomorphism $\Phi_F$, we then endow $\ol{F}^{Fins}$ with the structure of a manifold with corners as well. The homeomorphism $\Phi_F^{-1}$ sends each face $\tau^*$ of $B^*$ (dual to the face $\tau$ of $B$, which we will identify with the corresponding face of the Coxeter complex at infinity $\amod$) to the ideal boundary
$$
\geo^{Fins} V(0, \tau).
$$
The latter can be described as the set of strong asymptote classes of sectors $V(x,\tau)$:
$$
[V(x,\tau)]=[V(x',\tau)] \iff x\equiv x' \in F/Span(V(0,\tau)), 
$$
see Lemma \ref{lem:distlifct}. In other words, this is the stratum $X_{\tau}^F$ of 
$ \ol{F}^{Fins}$, see \eqref{eq:strata_of_bar_F}.    
The goal of this section is to extend this manifold with corners structure  from 
$\ol{F}^{Fins}$ to $\ol{X}^{Fins}$. 
We will also see that this structure matches the one of the maximal Satake compactification of $X$.

\subsection{Manifold-with-corners}
\label{sec:cor}

Let $\si\subset\geo X$ be a chamber which we view as a point in the closed stratum $\DF X$ of $\ol X^{Fins}$.
Let $o\in X$ be the fixed point of $K$.
\begin{lem}
For every neighborhood $U$ of $\si$ in $\ol{V(o,\si)}^{Fins}$
and every neighborhood $U'$ of the identity $e$ in $K$,
the subset $U'\cdot U$ is a neighborhood  of $\si$ in $\ol X^{Fins}$. 
\end{lem}
\proof
Suppose that $U'\cdot U$ is not a neighborhood.
Then there exists a sequence $\xi_n\to\si$ in $\ol X^{Fins}$ outside $U'\cdot U$.
There exist chambers $\si_n$ such that $\xi_n\in\ol{V(o,\si_n)}^{Fins}$,
and points $y_n\in V(o,\si_n)$ approximating $\xi_n$ such that $y_n\to\si$. 
Our description of sequential convergence, cf. Proposition~\ref{prop:convinfi},
implies that the sequence $(y_n)$ is $\simod$-regular and $\si_n\to\si$.
Hence there exist elements $k_n\to e$ in $K$ 
such that $k_n\si=\si_n$.
Then, due to the continuity of the $K$-action,
the points $k_n^{-1}\xi_n\in \ol{V(o,\si)}^{Fins}$ converge to $\si$.
Hence they enter the neighborhood $U$, and $(k_n)$ enters $U'$ for large $n$.
This is a contradiction.
\qed

\medskip
Suppose now that the neighborhood $U\subset\ol{V(o,\si)}^{Fins}$ is sufficiently small, say,
disjoint from the union of the compactified sectors $\ol{V(o,\tau)}^{Fins}$ 
over all proper faces $\tau\subsetneq\si$. 
Then the stabilizer of every point in $U$ equals the pointwise stabilizer $K_{\si}=K_F$ of the maximal flat 
$F\supset V(o,\si)$, see Corollary~\ref{cor:stabs}.
We consider the bijective continuous map
$$ K/K_F \times U \lra KU\subset \ol X^{Fins} $$
given by the $K$-action. 
By the previous lemma, its image $KU$ is a neighborhood of the closed stratum $S_{\simod}=\DF X$. 
After shrinking $U$ to a compact neighborhood of $\si$, 
the map becomes a homeomorphism.
After further shrinking $U$ to an open neighborhood, 
the map becomes a homeomorphism onto an open neighborhood of $\DF X$. 

Since $U$ is a manifold with corners, see Theorem \ref{thm:constructphi}, 
and $K/K_F$ is a manifold, 
we conclude with Lemma~\ref{lem:zoom}: 
\begin{thm}[Manifold-with-corners]
\label{thm:corny}
$\ol X^{Fins}$ is a manifold with corners with respect to the stratification 
by the strata $S_{\taumod}$.
In particular, the manifold-with-corners structure is $G$-in\-va\-riant.
\end{thm}

This means that the $k$-dimensional stratum of the manifold with corner structure 
equals the union of the $k$-dimensional strata $S_{\taumod}$.

\subsection{Homeomorphism to ball}\label{sec:Compactified symmetric space is a ball}

At last, we can now prove that the Finsler compactification of the symmetric space $X$ is $K$-equivariantly homeomorphic to a closed ball. Let $B^*$ be the dual ball to the unit ball $B\subset F_{mod}$ of the norm \eqref{eq:norm} on the vector space $F_{mod}$, defined via the regular vector $\bar\theta$. We will identify the dual vector space of $F_{mod}$ with $F_{mod}$ itself using the euclidean metric on $F_{mod}$.
Hence, $B^*$ becomes a unit ball in $F_{mod}$ for the {\em dual norm} 
$$
||\cdot||^*= ||\cdot||_{\bar\theta}^*
$$ 
of our original norm. 

Since $B^*\subset F_{mod}$ is $W$-invariant, the dual norm extends from $F_{mod}$ to a $G$-invariant 
function $d^*_{\bar\theta}$ on $X\times X$ by 
$$
d^*_{\bar\theta}(x,y)= ||d_\Delta(x,y)||^*_{\bar\theta} .
$$
We call the set 
$$
{B}^*(o, 1)= \{q\in X: d^*_{\bar\theta}(o,q)\le 1\}
$$
the {\em dual ball}.
It is preserved by the group $K$ since $K$ fixes the point $o$. 
As a compact star-like subset of $(X,d^{Riem})$, it is homeomorphic to the closed ball. 
We can now prove: 

\begin{thm}\label{thm:homeo-to-ball}
There exists a $K$-equivariant homeomorphism  
$$\ol{X}^{Fins}\stackrel{\Phi}{\lra} {B}^*(o, 1)$$ 
which restricts to the homeomorphism
$\phi: \ol{\Delta}^{Fins}\to \Delta^*_{B^*}$
from Theorem \ref{thm:constructphi}. 
In particular, $\ol{X}^{Fins}$ is homeomorphic to the closed ball. 
\end{thm} 
\proof We will use Lemma \ref{lem:transformation lemma} to construct $\Phi$. 
In order to do so, we have to know that $\ol{\Delta}^{Fins}$ and $\Delta^*_{B^*}$ 
are cross sections for the actions of $K$ on $\ol{X}^{Fins}$ and ${B}^*(o, 1)$,
and that $\phi$ respects the $K$-stabilizers. 

1. According to Lemma \ref{lem:cross-section}, $\ol{\Delta}^{Fins}$ is a cross section for the action of $K$ on $\ol{X}^{Fins}$. Since $K$ preserves the dual ball ${B}^*(o, 1)$ and
$$
\Delta^*_{B^*}= \Delta\cap {B}^*(o, 1),
$$
while $\Delta$ is a cross section for the action $K\acts X$, it follows that $\Delta^*_{B^*}$ is a cross section 
for the action $K\acts {B}^*(o, 1)$. 

2. The faces $\tau$, $\emptyset\subseteq\tau\subseteq\si$,
correspond to index sets $J$, $\emptyset\subseteq J_{\tau}\subseteq [n]$,
where $j\in J_{\tau}$ iff the reflection $s_j$ fixes $\tau$.
According to Corollary~\ref{cor:stabs},
the fixed point set of $K_{\tau}$ on $\ol\De^{Fins}$ equals $\ol{V(o,\tau)}^{Fins}$. 
On the other hand,
the fixed point set of $K_{\tau}$ on $\De^*_{B^*}$ equals the interior face $F_{J_{\tau}}$. 
By Theorem~\ref{thm:constructphi},
the homeomorphism $\phi$ carries $\ol{V(o,\tau)}^{Fins}$ to $F_{J_{\tau}}$. 
Therefore, $\phi$ respects the point stabilizers. 
\qed

\subsection{Relation to the maximal Satake compactification}

It turns out that the compactification $\ol{X}^{Fins}$ constructed in this paper is naturally isomorphic to the {\em maximal Satake compactification} $\ol{X}^S_{max}$. To this end, we will use the {\em dual-cell} interpretation of the maximal Satake compactification, see \cite[Ch. I.19]{Borel-Ji}.

\begin{thm}\label{thm:satake}
There is a $G$-equivariant homeomorphism of manifolds-with-corners $\ol{X}^{Fins}\to \ol{X}^S_{max}$ which extends the identity map $X\to X$. 
\end{thm} 
\proof We first observe that the group $K$ acts on both compactifications so that the cross sections for the actions are the respective compactifications of the model euclidean Weyl chamber $\De=\De_{mod}\subset F=\Fmod$. 
We therefore compare the $W$-invariant compactifications of $F_{mod}$. 
On the side of  $\ol{X}^{Fins}$, the ideal boundary of $F$ is the union of small strata $X_{\tau}^{F}$ 
as in \S \ref{sec:compmaxfl}.  
Elements of $X_{\tau}^F$ are equivalence classes $[V(x,\tau)]$ of sectors $V(x,\tau)$ in $F$. Two sectors $V(x,\tau), V(x',\tau)$ with the same base $\tau$ are equivalent iff $x, x'$ project to the same vector in $F/Span(V(0,\tau))$. 
These are exactly the strata, denoted $e(C)$, in the maximal Satake compactification of $F$, denoted by 
$\ol{F}^S_{max}$,
see \cite[Ch. I.19]{Borel-Ji}: For each sector $C=V(0,\tau)$,  the stratum $e(C)$ is 
$F/Span(C)$. We then have a $W$-equivariant bijection 
$$h: \ol{F}^{Fins}\to\ol{F}^S_{max}$$
defined via the collection of maps 
$$ [V(x,\tau)] \mapsto [x]\in e(C) .$$
For $\tau=\emptyset$, this is just the identity map $F\to F$.

In order to show that this map is a homeomorphism we note that the topology on $\ol{F}^S_{max}$ is defined via roots (see \cite[Ch. I.19]{Borel-Ji}) and on the Weyl chamber $\Delta$ in  $F$ this topology is exactly the  topology on $\ol{\Delta}^{Fins}$ described in terms of simple roots, cf. the proof of Lemma \ref{lem:matchingfaces}.

Lastly, we note that the map $h$ we described respects the stabilizers in the group $K$. Therefore, by Lemma \ref{lem:transformation lemma}, we obtain a $K$-equivariant homeomorphic extension 
$$
\ol{X}^{Fins}\to \ol{X}^S_{max}$$
of $h$,
which is also an extension of the identity map $X\to X$.
Since the identity is $G$-equivariant, the same holds for the extension.
\qed

\begin{rem}
The maximal Satake compactification is a real-analytic manifold with corners on which the group $G$ acts real-analytically, see \cite[Ch. I.19]{Borel-Ji}. Therefore, the same conclusion holds for the compactification 
$\ol{X}^{Fins}$. 
\end{rem}

\subsection{Proof of Theorem \ref{thm:comp}}

The theorem is the combination of the following results:

Part (i) is proven in section~\ref{sec:Stratification}, where we established that 
$\ol{X}^{Fins}$ is a union of strata $S_{\taumod}$ each of which is a single $G$-orbit. Thus, $G$ acts on 
$\ol{X}^{Fins}$ with finitely many orbits.

Part (ii) is proven in Theorem \ref{thm:corny}. 

Part (iii) is proven in Theorem \ref{thm:homeo-to-ball}. 
 
Part (iv) is the content of Proposition \ref{prop:independence}. 

Lastly, Part (v) is established in Theorem \ref{thm:satake}. \qed

\section{Relative position and thickenings}

\subsection{Relative position at infinity and folding order}
\label{sec:thick}

In this section, 
we review some combinatorial concepts from \cite{coco15} related to the geometry of Tits buildings. 
We will discuss here only the relative position of chambers with respect to simplices,
which is the case needed in this paper,
and refer the reader to \cite{coco15} for more general treatment.

Let $\si_0,\si\subset\geo X$ be chambers.
We view them also as points $\si_0,\si\in\DF X$.
There exists an (in general non-unique) apartment $a\subset\geo X$ containing these chambers,
$\si_0,\si\subset a$,
and a unique apartment chart $\al:\amod\to a$ such that 
$\si_0=\al(\simod)$.
We define the {\em position of $\si$ relative to $\si_0$} as the chamber 
$$ \pos(\si,\si_0) := \al^{-1}(\si)\subset\amod .$$
Abusing notation, it can be regarded algebraically as the unique element 
$$ \pos(\si,\si_0)\in W $$
such that 
$$ \si=\al\bigl(\pos(\si,\si_0)\simod\bigr) ,$$
cf. \cite[\S 3.3]{coco15}.
It does not depend on the choice of the apartment $a$. 
To see this, choose regular points $\xi_0\in\interior(\si_0)$ and $\xi\in\interior(\si)$
which are not antipodal, $\tangle(\xi,\xi_0)<\pi$.
Then the segment $\xi_0\xi$ is contained in $a$ by convexity,
and its image $\al^{-1}(\xi_0\xi)$ in $\amod$ is independent of the chart $\al$
because its initial portion $\al^{-1}(\xi_0\xi\cap\si_0)$ in $\simod$ is.

The level sets of $\pos(\cdot,\si_0)$ in $\DF X$ are the {\em Schubert cells} relative $\si_0$,
i.e. the orbits of the Borel subgroup $B_{\si_0}\subset G$ fixing $\si_0$.

More generally, 
we define the position of a chamber $\si\subset\geo X$ relative to a {\em simplex} $\tau_0\subset\geo X$
as follows.
We denote $\taumod=\theta(\tau_0)$.
Let again $a\subset\geo X$ be an apartment containing $\tau_0$ and $\si$,
and let $\al:\amod\to a$ be a chart such that 
$\tau_0=\al(\taumod)$.
We define the {\em position $\pos(\si,\tau_0)$ of $\si$ relative to $\tau_0$} as the $\Wt$-orbit 
of the chamber 
$$\al^{-1}(\si)\subset\amod.$$
It can be interpreted algebraically as a coset 
$$ \pos(\si,\tau_0)\in W_{\taumod}\backslash W .$$ 
The (strong) {\em Bruhat order} ``$\preceq$" on the Weyl group $W$ 
has the following geometric interpretation as {\em folding order}, 
cf. \cite[Def.\ 3.2]{coco15}.
For distinct elements $w_1,w_2\in W$, it holds that $$w_1\prec w_2$$
if and only if there exists a folding map $\amod\to\amod$ fixing $\simod$ and mapping
$w_2\simod\mapsto w_1\simod$,
cf. \cite[\S 3.2]{coco15}.
Here, by a {\em folding map} $\amod\to\amod$ 
we mean a type preserving continuous map which sends chambers isometrically onto chambers.

The folding order on relative positions coincides with the inclusion order on Schubert cycles,
i.e. $w_1\preceq w_2$ if and only if the Schubert cell $\{\pos(\cdot,\si_0)=w_1\}$
is contained in the closure of the Schubert cell $\{\pos(\cdot,\si_0)=w_2\}$, 
and the {\em Schubert cycles} relative $\si_0$ are the sublevel sets of $\pos(\cdot,\si_0)$. 
In the case of complex semisimple Lie groups $G$ 
this inclusion relation is a classical result of Chevalley \cite{Chevalley}; 
for the case of general semisimple Lie groups we refer the reader to 
\cite[Prop 3.14]{coco15} or, alternatively, to \cite{Mitchell, Mitchell2008}.

We also need to define the folding order more generally on positions of chambers relative to simplices $\tau_0$ 
of an arbitrary face type $\taumod$.
We say that 
$$
W_{\taumod} \bar\si_1 \preceq_{\taumod} W_{\taumod} \bar\si_2
$$ 
for chambers $\bar\si_1,\bar\si_2\subset\amod$
iff there exist $\bar\si_i'\in W_{\taumod} \bar\si_i$ such that
$$
\bar\si_1'\preceq \bar\si_2' ,
$$
equivalently, geometrically, 
if for some (any) chambers $\bar\si_i'\in W_{\taumod} \bar\si_i$ 
there exists a folding map $\amod\to\amod$ fixing $\taumod$ and mapping $\bar\si'_2$ to $\bar\si'_1$.
(Note that the elements in $W_{\taumod}$ are such folding maps.)
\begin{lem}
$\prec_{\taumod} $ is a partial order. 
\end{lem}
\proof 
Transitivity holds, because the composition of folding maps is again a folding map.

To verify reflexivity,
pick points $\xi_{mod}\in\inte(\taumod)$ and $\eta_{mod}\in\inte(\simod)$.

Let $\bar\si=w\simod\subset\amod$ be a chamber and $f:\amod\to\amod$ a folding map fixing $\taumod$.
We denote $\bar\eta=w\eta_{mod}$.
If the $f$-image of the segment $\xi_{mod}\bar\eta$ is again an unbroken geodesic segment,
then the two geodesic segments are congruent by an element of $W_{\taumod}$,
because their initial directions at $\xi_{mod}$ are.
On the other hand, if the $f$-image of $\xi_{mod}\bar\eta$ is a broken geodesic segment,
then the distance of its endpoints is strictly smaller than its length,
and consequently
$f\bar\si\not\succeq\bar\si$.
This shows that
$$W_{\taumod} \bar\si_1 \preceq_{\taumod} W_{\taumod} \bar\si_2\preceq_{\taumod} W_{\taumod} \bar\si_1
\quad\Ra\quad
W_{\taumod} \bar\si_1 =W_{\taumod} \bar\si_2$$
and hence reflexivity.
\qed

\medskip
The relative position function
$$ \pos: \Flags\times\Flagt\lra W_{\taumod}\backslash W $$ 
is {\em lower semicontinuous},
compare the discussion of closures of Schubert cycles above.

\medskip
It follows from the geometric description of the folding orders in terms of folding maps
that for face types $\numod\subset\taumod$ the order $\preceq_{\numod}$ refines the order $\preceq_{\taumod}$,
because a folding map fixing $\taumod$ fixes in particular its face $\numod$.
Thus,
for chambers $\si_1,\si_2\subset\geo X$ 
and simplices $\nu\subset\tau\subset\geo X$ of types $\numod=\theta(\nu)\subset\taumod=\theta(\tau)$ 
it holds that:
\begin{equation}\label{eq:monotonic}
\pos(\si_1,\tau)\preceq_{\taumod}\pos(\si_2,\tau)
\quad\Ra\quad
\pos(\si_1,\nu)\preceq_{\numod}\pos(\si_2,\nu)\end{equation}
We now describe the action of the longest element $w_0\in W$ on relative positions. 
Note that 
$w_0\Wt w_0^{-1}=W_{w_0\taumod}=W_{\iota\taumod}$
and $w_0$ maps $\Wt$-orbits to $W_{\iota\taumod}$-orbits. 
The action of $w_0$ therefore induces a natural map 
\begin{eqnarray*}
\begin{aligned}
\Wt\backslash W
&\buildrel w_0 \over\lra &
W_{\iota\taumod}\backslash W &&\\
\Wt w &\mapsto& w_0\Wt w &=&W_{\iota\taumod} w_0w .
\end{aligned}
\end{eqnarray*}
This map is order reversing,
\begin{equation}
\label{eq:actwnd}
\Wt w \prect \Wt w'
\quad\Leftrightarrow\quad
W_{\iota\taumod}w_0 w \succ_{\iota\taumod} W_{\iota\taumod}w_0 w' ,
\end{equation}
see \cite[\S 3.2]{coco15}.

\begin{dfn}[Complementary position]
\label{dfn:complemtau}
We define the {\em complementary position} by
$$\cpos:=w_0\pos .$$
\end{dfn}
This terminology is justified by:
\begin{lem}[{\cite[Lemma 3.16]{coco15}}]
Let $\tau,\hat\tau,\si\subset\geo X$ be two simplices and a chamber contained in an apartment $a$,
and suppose that $\tau$ and $\hat\tau$ are antipodal.
Then $\pos(\si,\hat\tau)=\cpos(\si,\tau)$.
\end{lem}
\medskip
The relation of ``complementarity'' is clearly symmetric,
$\ccpos=\pos$.
Passing to complementary relative position reverses the partial order,
cf.\ (\ref{eq:actwnd}):
\begin{equation}
\label{ineq:revordpos}
\pos(\si_1,\tau)\prec_{\theta(\tau)}\pos(\si_2,\tau) 
\quad\Leftrightarrow\quad
\cpos(\si_1,\tau)\succ_{\iota\theta(\tau)}\cpos(\si_2,\tau) 
\end{equation}

\subsection{Further properties of the folding order}

This is a technical section whose main result (Proposition~\ref{prop:dsopnstr}) 
will be used in the proof of Proposition~\ref{prop:drfinq}
which is the key to proving proper discontinuity of actions of $\taumod$-regular subgroups.

We begin with a result useful for comparing relative positions. 
\begin{lem}
\label{lem:monalsegm}
(i) Let $\si_0,\si_1,\si_2\subset\geo X$ be chambers, 
and suppose that there exists a segment $\xi_0\xi_2$ 
with $\xi_0\in\inte(\si_0)$ and $\xi_2\in\inte(\si_2)$ 
containing a point $\xi_1\in\inte(\si_1)$.
Then 
$$\pos(\si_1,\si_0)\preceq\pos(\si_2,\si_0)$$
with equality iff $\si_1=\si_2$. 

(ii) More generally, let $\si_1,\si_2\subset\geo X$ be chambers 
and let $\tau_0\subset\geo X$ be a simplex of type $\taumod$.
Suppose that there exists a segment $\xi_0\xi_2$ 
with $\xi_0\in\inte(\tau_0)$ and $\xi_2\in\inte(\si_2)$ 
containing a point $\xi_1\in\inte(\si_1)$.
Then 
$$\pos(\si_1,\tau_0)\preceq_{\taumod}\pos(\si_2,\tau_0)$$
with equality iff $\si_1=\si_2$.

\end{lem}
\proof
We prove the more general assertion (ii). 
After perturbing $\xi_2$, we can arrange that the subsegment $\xi_1\xi_2$ avoids codimension two faces.
Along this subsegment we find a gallery of chambers connecting $\si_1$ to $\si_2$.
We may therefore proceed by induction and assume that the chambers $\si_1$ and $\si_2$ are adjacent,
i.e. share a panel $\pi$ which is intersected transversally by $\xi_1\xi_2$.
Working in an apartment containing $\tau_0,\si_1,\si_2$,
the wall through $\pi$ does not contain $\tau_0$
and separates $\st(\tau_0)\cup\si_1$ from $\si_2$. 
Folding at this wall yields the desired inequality.
\qed

\medskip
We next study the values of the relative position function on stars of simplices.

We fix a reference chamber $\si_0\subset\geo X$.
Let $\nu\subset\geo X$ be a simplex.
For any interior points $\eta\in\inte(\nu)$ and $\xi_0\in\inte(\si_0)$,
the segment $\eta\xi_0$ enters the {\em interior} of a chamber $\si_-\supset\nu$,
i.e. 
$$\eta\xi_0\cap\inte(\si_-)\neq\emptyset.$$
Note that the chamber $\si_-$ does not depend on the interior points $\eta,\xi_0$.
Moreover, it is contained in any apartment containing $\si_0$ and $\nu$.
We call $\si_-$ the chamber in $\st(\nu)$ {\em pointing towards $\si_0$}.

Similarly, if $\xi_0\xi_+\supsetneq\xi_0\eta$ is an extension of the segment $\xi_0\eta$ beyond $\eta$,
then there exists a chamber $\si_+\supset\nu$ such that 
$\eta\xi_+\cap\inte(\si_+)\neq\emptyset$,
and we call $\si_+$ {\em a} chamber in $\st(\nu)$ {\em pointing away from $\si_0$}.

\medskip
Let $a\subset\geo X$ be an apartment containing $\si_0$ and $\nu$.
Then $\si_-\subset a$.
Moreover, since geodesic segments inside $a$ extend uniquely, 
there exists a {\em unique} chamber $\si_+\subset\st(\nu)\cap a$ pointing away from $\si_0$.
The chambers $\si_{\pm}\subset a$ can be characterized as follows in terms of separation from $\si_0$ by walls:
\begin{lem}
\label{lem:wallsep}
Let $\si\subset\st(\nu)\cap a$ be a chamber. 
Then 

(i) $\si=\si_+$ iff $\si$ is separated from $\si_0$ by every wall $s\subset a$ containing $\nu$. 

(ii) $\si=\si_-$ iff $\si$ is not separated from $\si_0$ by any wall $s\subset a$ containing $\nu$. 
\end{lem}
\proof
(i) Clearly, $\si_+$ is separated from $\si_0$ by every wall $s\supset\nu$
because, using the above notation, $\xi_0\xi_+\cap s=\eta$.
Vice versa, if $\si$ is separated from $\si_0$ by all such walls $s$, 
then $\si$ and $\si_+$ lie in the same hemispheres bounded by the walls $s\supset\nu$ in $a$,
and therefore must coincide. 

(ii) Similarly, $\si_-$ is not separated from $\si_0$ by any wall $s\supset\nu$
because $\xi_0\eta\cap s=\eta$,
and vice versa, if $\si$ is not separated from $\si_0$ by any wall $s\supset\nu$, 
then $\si$ and $\si_-$ lie in the same hemispheres bounded by the walls $s\supset\nu$ in $a$,
and therefore must coincide. 
\qed

\begin{rem}
\label{rem:contpan}
The assertion of the lemma remains valid if one only admits the walls $s\subset a$ such that $s\cap\si$ 
is a panel containing $\nu$.
\end{rem}

The chambers pointing towards and away from $\si_0$ in $\geo X$ can also be characterized in terms of the folding order:
\begin{lem}
\label{lem:minmaxonst}
The restriction of the function $\pos(\cdot,\si_0)$ to the set of chambers contained in $\st(\nu)$
attains a unique minimal value in $\si_-$
and a unique maximal value\footnote{Meaning that it is larger than all other values.}
precisely in the chambers pointing away from $\si_0$.
\end{lem}
\proof
Let $\si\supset\nu$ be a chamber
and let $a\subset\geo X$ be an apartment containing $\si_0$ and $\si$.
Then $\si_-\subset a$. 
Let $\si_+\subset\st(\nu)\cap a$ be the unique chamber pointing away from $\si_0$.

Still using the above notation, 
let $\xi_0\xi_+\supset\xi_0\eta$ be an extension of the segment $\xi_0\eta$ 
with endpoint $\xi_+\in\inte(\si_+)$.
Let $\xi_-\in\xi_0\eta\cap\inte(\si_-)$.
The points $\xi_-$ and $\eta$ appear in this order on the (oriented) segment $\xi_0\xi_+$. 

We now perturb the segment $\xi_0\xi_+$ to a segment $\xi_0\xi'_+$ which intersects 
$\inte(\si)$ in a point $\eta'$ close to $\eta$ and $\inte(\si_-)$ in a point $\xi'_-$ close to $\xi_-$.
The perturbation is possible because $\si\supset\nu$.
Again, 
the points $\xi'_-$ and $\eta'$ appear in this order on the perturbed segment $\xi_0\xi'_+$. 
Lemma~\ref{lem:monalsegm} implies that 
$$\pos(\si_-,\si_0)\preceq\pos(\si,\si_0)\preceq\pos(\si_+,\si_0) $$
with equality in the first (second) inequality if and only if $\si=\si_-$ ($\si=\si_+$).
The assertion of the lemma follows because $\pos(\si_+,\si_0)$ does not depend on the choice of $a$.
\qed

\medskip
We now extend the lemma to the case of relative position with respect to a simplex $\tau_0$.

We pick a chamber $\si_0\supset\tau_0$. 
For a simplex $\nu$,
the chamber $\si_-\subset\st(\nu)$ pointing towards $\si_0$ is defined as before. 

\begin{cor}
\label{cor:minmaxonst}
The restriction of the function $\pos(\cdot,\tau_0)$ to the set of chambers in $\st(\nu)$
attains a unique minimal value in $\si_-$ and a unique maximal value
precisely in the chambers pointing away from $\si_0$.
\end{cor}
\proof This is an immediate consequence of Lemma \ref{lem:minmaxonst} and \eqref{eq:monotonic}. \qed 

\medskip
Combing the corollary with the discussion in section~\ref{sec:idbdcrsct}, 
we obtain:

\begin{prop}
\label{prop:dsopnstr}
Let $\tau_0,\nu$ be simplices.
Then there exists a dense open subset of chambers in $\stF(\nu)$
where the function $\pos(\cdot,\tau_0)$ attains its unique maximal value.
\end{prop}

\subsection{Thickenings}
\label{sec:thickenings}

A {\em thickening} (of the neutral element) in $W$ is a subset $$\Th\subset W$$
which is a union of sublevels for the folding order,
i.e.\ which contains with every element $w$ also every element $w'$ satisfying $w'\prec w$,
cf. \cite[Def 4.16]{coco15}.
In the theory of posets, such subsets are called {\em ideals}.

Note that 
$$ \Th^c := w_0(W-\Th) $$
is again a thickening.
Here, $w_0\in W$ denotes the longest element of the Weyl group,
that is, the element of order two mapping $\simod$ to the opposite chamber in $\amod$.
It holds that 
$$ W =\Th\sqcup w_0\Th^c $$
and we call $\Th^c$ the thickening {\em complementary} to $\Th$.

The thickening 
$\Th\subset W$ 
is called {\em fat} if $\Th\cup w_0\Th=W$, equivalently, $\Th\supseteq\Th^c$.
It is called {\em slim} if $\Th\cap w_0\Th=\emptyset$, equivalently, $\Th\subseteq\Th^c$.
It is called {\em balanced} if it is both fat and slim, equivalently, $\Th=\Th^c$,
cf. \cite[Def 3.22]{coco15}.

For types $\tilde\vartheta_0,\tilde\vartheta\in\simod$ and a radius $r\in[0,\pi]$
we define the {\em metric thickening}
$$ \Th_{\tilde\vartheta_0,\tilde\vartheta,r} := \{w\in W:d(w\tilde\vartheta,\tilde\vartheta_0)\leq r\},$$
using the natural $W$-invariant spherical metric $d$ on $\amod$, 
cf. \cite[\S 3.4.1]{coco15}.

For a face type $\taumod\subseteq\simod$, we denote by $W_{\taumod}$ its stabilizer in $W$.
Furthermore, $\iota=-w_0:\simod\to\simod$ denotes the canonical involution of the model spherical Weyl chamber. 
\begin{lem}
(i) If $\tilde\vartheta_0\in\taumod$, then 
$W_{\taumod}\Th_{\tilde\vartheta_0,\tilde\vartheta,r}=\Th_{\tilde\vartheta_0,\tilde\vartheta,r}$.

(ii) If $\iota\tilde\vartheta_0=\tilde\vartheta_0$, then 
$\Th_{\tilde\vartheta_0,\tilde\vartheta,r}$ is fat for $r\geq\pihalf$ and slim for $r<\pihalf$.
\end{lem}
\proof
(i) For $w'\in W_{\taumod}$, we have that $w'\tilde\vartheta_0=\tilde\vartheta_0$ and hence 
$$ d(w'w\tilde\vartheta,\tilde\vartheta_0) = 
d(w\tilde\vartheta,\underbrace{{w'}^{-1}\tilde\vartheta_0}_{\tilde\vartheta_0}).$$
(ii) Since $w_0\tilde\vartheta_0=-\iota\tilde\vartheta_0=-\tilde\vartheta_0$, we have 
$$ d(w_0w\tilde\vartheta,-\tilde\vartheta_0) = 
d(w\tilde\vartheta,\underbrace{-w_0\tilde\vartheta_0}_{\tilde\vartheta_0}),$$
whence the assertion.
\qed
\begin{cor}[Existence of balanced thickenings]
If the face type $\taumod$ is $\iota$-invariant, $\iota\taumod=\taumod$,
then there exists a $W_{\taumod}$-invariant balanced thickening $\Th\subset W$.
\end{cor}
\proof
Since $\iota\taumod=\taumod$,
there exists $\tilde\vartheta_0\in\taumod$ 
such that $\iota\tilde\vartheta_0=\tilde\vartheta_0$.
Moreover, 
there exists $\tilde\vartheta\in\simod$ 
such that $d(\cdot\tilde\vartheta,\tilde\vartheta_0)\neq\pihalf$ on $W$.
(This holds for an open dense subset of types $\tilde\vartheta\in\simod$.)
According to the lemma,
the metric thickening $\Th_{\tilde\vartheta_0,\tilde\vartheta,\pihalf}$ is balanced and $W_{\taumod}$-invariant.
\qed

\medskip
Given a thickening $\Th\subset W$,
we obtain {\em thickenings at infinity} as follows.

First, we define the thickening in $\DF X$ of a chamber $\si\in\DF X$ as
$$ \ThF(\si):=\{\pos(\cdot,\si)\in\Th\}\subset\DF X .$$
It is a finite union of Schubert cycles relative $\si$.
We then define the thickening of $\si$ inside the Finsler ideal boundary 
as the ``suspension" of its thickening inside the Furstenberg boundary,
$$ \Th^{Fins}(\si):=\{[b^{\bar\theta}_{\nu,p}] : \stF(\nu)\subset\ThF(\si)\} 
=\bigcup \bigl\{X_{\nu}:\stF(\nu)\subset\ThF(\si) \bigr\} 
\subset\geo^{Fins}X$$
where $\stF(\nu)$ denotes the set of chambers containing $\nu$ as a face,
see section~\ref{sec:idbdcrsct}.
Note that 
$\ThF(\si)=\Th^{Fins}(\si)\cap\DF X$.

\begin{lem}
$\Th^{Fins}(\si)$ is compact.
\end{lem}
\proof
Consider a sequence of points $[b^{\bar\theta}_{\nu_n,p_n}]\in\Th^{Fins}(\si)$, 
and suppose that it converges in 
$\geo^{Fins}X$, 
$$[b^{\bar\theta}_{\nu_n,p_n}]\to[b^{\bar\theta}_{\mu,q}].$$
We must show that also $[b^{\bar\theta}_{\mu,q}]\in\Th^{Fins}(\si)$.

After extraction,
we may assume that all simplices $\nu_n$ have the same type $\theta(\nu_n)=\numod$.
According to Lemma~\ref{lem:inbigrstr},
$\nu_n\to\nu\subseteq\mu$.
By assumption, $\st(\nu_n)\subset\ThF(\si)$,
and we must show that $\st(\mu)\subset\ThF(\si)$.
Since $\st(\nu)\supseteq\st(\mu)$, this will follow from $\st(\nu)\subset\ThF(\si)$.

The latter follows from the closedness of $\ThF(\si)$ in $\DF X$,
because every chamber $\si'\subset\st(\nu)$ is a limit of a sequence of chambers $\si'_n\subset\st(\nu_n)$.
\qed

\begin{rem}
\label{rem:thcontr}
One can show that $\Th^{Fins}(\si)\subset\geo^{Fins}X$ is a contractible CW-complex.
In the second version of this paper (see Theorem~8.21 there) we proved that it is \v{C}ech acyclic.
\end{rem}

\begin{example}
Suppose that the Weyl group $W$ of $X$ is of type $A_2$, 
i.e.\ is isomorphic to the permutation group on 3 letters. 
Let $s_1, s_2\in W$ denote the generators which are the reflections in the walls of the positive chamber $\simod$. 
There is the unique balanced thickening
$\Th=\{e, s_1, s_2\}\subset W$. 
The  thickening $\Th^{Fins}(\si)\subset \geo^{Fins}X$ is the wedge of two closed disks connected at the point  $\si$: These disks are the visual compactifications $\ol{X}_{\tau_i}, i=1, 2$, of two rank 1 symmetric spaces $X_{\tau_i}$. Here $\tau_1, \tau_2$ are the two vertices of the edge $\si$.  
\end{example}

More generally, 
we define the thickening in $\geo^{Fins}X$ of a set of chambers $A\subset\DF X$ as
$$ \Th^{Fins}(A):=\bigcup_{\si\in A}\Th^{Fins}(\si) \subset\geo^{Fins}X .$$

\begin{lem}
If $A$ is compact, then $\Th^{Fins}(A)$ is compact.
\end{lem}
\proof
Since $\DF X$ is a homogeneous space for the maximal compact subgroup $K$, 
there exists a chamber $\si_0\in A$ and a compact subset $C\subset K$ such that $A=C\si_0$.
Then $$\Th^{Fins}(A)=C\cdot\Th^{Fins}(\si_0)$$ and is hence compact
as a consequence of the previous lemma. 
\qed

\medskip
If the thickening $\Th\subset W$ is $W_{\taumod}$-invariant,
then we can define the thickening in $\geo^{Fins}X$ of a simplex $\tau\subset\geo X$ of type $\taumod$ as
$$ \Th^{Fins}(\tau):=\Th^{Fins}(\si)\subset\geo^{Fins}X $$
for a chamber $\si\supseteq\tau$.
It does not depend on $\si$.
For a set $A\subset\Flagt$ of simplices of type $\taumod$,
we define its thickening in $\geo^{Fins}X$ as
$$ \Th^{Fins}(A):=\bigcup_{\tau\in A}\Th^{Fins}(\tau)\subset\geo^{Fins}X .$$
Again, $\Th^{Fins}(A)$ is compact if $A$ is.

\begin{lem}[Fibration of thickenings]
\label{lem:thickfib}
Let $A\subset\Flagt$ be compact, 
and suppose that 
the thickenings $\Th^{Fins}(\tau)$ of the simplices $\tau\in A$
are pairwise disjoint. 
Then the natural map 
$$ \Th^{Fins}(A)\buildrel\pi\over\lra A$$
is a fiber bundle.
\end{lem}
\proof
Regarding continuity of $\pi$, 
suppose that $\xi_n\to\xi$ in $\Th^{Fins}(A)$ and $\tau_n\to\tau$ in $A$
with $\xi_n\in\Th^{Fins}(\tau_n)$.
Then $\xi\in\Th^{Fins}(\tau)$ by 
semicontinuity of relative position, and hence $\pi(\xi)=\tau$. 

In order to show that $\pi$ is a fibration, we need to construct local trivializations. 
Fix $\tau\in A$ and an opposite simplex $\hat\tau$. 
The horocyclic subgroup $H_{\hat\tau}$ acts simply transitively on an open neighborhood of $\tau$ in $\Flagt$. 
Now, let $S\subset H_{\hat\tau}$ denote 
the closed subset consisting of all $h\in H_{\hat\tau}$ which send $\tau$ to elements of $A$. 
Then $S\tau$ is a neighborhood of $\tau$ in $A$. 
Restricting the action of $H_{\hat\tau}$ to the subset $S$,  we obtain a topological embedding 
$$ S\times\Th^{Fins}(\tau)\to\Th^{Fins}(A)$$
and a local trivialization of $\pi$ over a neighborhood of $\tau$ in $A$.
\qed

\section{Proper discontinuity}
\label{sec:pd}

Our aim is to construct domains of proper discontinuity for the action 
$$ \Ga\acts\ol X^{Fins} $$
of discrete subgroups $\Ga<G$ on the Finsler compactification $\ol X^{Fins}$ of $X$.
The proper discontinuity of an action can be rephrased as the absence of dynamical relations,
and our construction of domains results from 
studying the dynamical relations between points in $\ol X^{Fins}$
with respect to the action $G\acts\ol X^{Fins}$
and determining necessary conditions.

The Furstenberg boundary is naturally embedded in the Finsler boundary 
as the closed stratum at infinity, 
$$\DF X\subset\geo^{Fins}X .$$
The $G$-action on $\ol X^{Fins}$ 
is determined by (``fills in'') the $G$-action on $\DF X\cong\Flags$,
and our approach is based on the study 
of the dynamics of the $G$-action on its associated flag manifolds 
in \cite{coco15}.
We first recall from there a combinatorial inequality for dynamical relations in $\DF X$ 
and provide an auxiliary result regarding the dynamics of pure sequences on $\DF X$,
cf.\ section~\ref{sec:dynfu}.
We then show in section~\ref{sec:derdynrelfi} 
how dynamical relations in $\geo^{Fins}X$ 
imply dynamical relations in $\DF X$
and use this to extend the combinatorial inequality for dynamical relations from $\DF X$ to $\ol X^{Fins}$. 
With this inequality at hand,
one readily obtains domains of proper discontinuity 
by removing suitable thickenings of limit sets, see section~\ref{sec:thickenings},
i.e.\ the points which have ``sufficiently special'' position relative to some limit point.

\subsection{Dynamics on the Furstenberg boundary}
\label{sec:dynfu}

We consider the action 
$$ G\acts\DF X $$
on the Furstenberg boundary $\DF X\cong\Flags$.
Specifically,
we are interested in the dynamics of diverging sequences in $G$.
Let $(g_n)$ be a $\taumod$-contracting sequence in $G$ with 
\begin{equation}
\label{eq:contrt}
g_n|_{C(\tau_-)}\to\tau_+
\end{equation}
uniformly on compacts,
where $\tau_{\pm}\in\Flagpmt$.

\subsubsection{Dynamical relations}
\label{sec:dynrelfue}

We recall from \cite{coco15}
the following necessary condition for dynamical relations 
between points in $\DF X$ with respect to the action of $(g_n)$:\footnote{\cite[Prop.\ 6.5]{coco15} applies to 
$\taumod$-contracting sequences on arbitrary flag manifolds $\Flagn$ 
for arbitrary face types $\taumod,\numod\subseteq\simod$.}

\begin{prop}[Dynamical relation inequality in $\DF X$, cf.\ {\cite[Prop.\ 6.5]{coco15}}]
\label{prop:dynreltaun}
Let $(g_n)$ be a $\taumod$-contracting sequence in $G$ is  with (\ref{eq:contrt}).
Suppose that there is a dynamical relation 
$$ \si\stackrel{(g_n)}{\sim}\si' $$
between points $\si,\si'\in\DF X$.
Then: 
\begin{equation}
\label{eq:dynrelposcondtaun}
\pos(\si',\tau_+)\prec\cpos(\si,\tau_-) 
\end{equation}
\end{prop}
Intuitively,
this means that it cannot happen that $\si$ is far from $\tau_-$ {\em and} $\si'$ is far from $\tau_+$,
where ``far'' is to be understood as having ``generic'' relative position. 

\proof
The Furstenberg boundary $\DF X$ is naturally identified 
with the regular $G$-orbits $G\xi$ in the visual boundary $\geo X$.
The assertion is therefore equivalent to the implication (i)$\Ra$(ii) of \cite[Prop.\ 6.5]{coco15}
in the special case of regular $G$-orbits.
\qed

\subsubsection{Pure sequences}
\label{sec:pudyfue}

The action of the $\taumod$-contracting sequence $(g_n)$ in $G$ satisfying (\ref{eq:contrt})
preserves the natural fibration of flag manifolds 
$$\pi_{\taumod}:\DF X\cong\Flags\to\Flagt .$$
For a simplex $\tau\in\Flagt$,
$$\stF(\tau)=\pi_{\taumod}^{-1}(\tau)
\subset\DF X$$
is the set of chambers $\si\supset\tau$, cf.\ section~\ref{sec:idbdcrsct}.
For $\tau_-\in\Flagit$, we denote by 
$$\CF(\tau_-):=\pi_{\taumod}^{-1}(C(\tau_-))\subset\DF X$$
the set of chambers over $C(\tau_-)$,
and by 
$\D\CF(\tau_-):=\pi_{\taumod}^{-1}(\D C(\tau_-))$ its complement in $\DF X$. 

The contraction property (\ref{eq:contrt}) for the action of $(g_n)$ on the base $\Flagt$
translates into the property for the dynamics on $\DF X$ that 
the $(g_n)$-orbits in $\CF(\tau_-)$ accumulate at $\stF(\tau)$ locally uniformly,
i.e.\ that for every compact subset $A\subset\CF(\tau_-)$ the sequence of subsets $g_nA$
accumulates at (a subset of) $\stF(\tau)$.
In the terminology of \cite[Def.\ 5.8]{coco15}
this means that the sequence $(g_n)$
is {\em $(\D\CF(\tau_-),\stF(\tau))$-accumulating} on $\DF X$.

As a $\taumod$-contracting sequence, $(g_n)$ is in particular $\taumod$-regular,
cf.\ section~\ref{sec:contrimplrg}.
If we make the stronger additional assumption that $(g_n)$ is {\em $\taumod$-pure}, 
cf.\ Definition~\ref{def:pureg}, 
then its accumulation dynamics on $\DF X$ can be described more precisely.
The next result expresses that there is only bounded distortion in the direction of the $\pi_{\taumod}$-fibers:

\begin{prop}
\label{prop:puredyn}
Suppose that the sequence $(g_n)$ in $G$ is 
$\taumod$-contracting with (\ref{eq:contrt})
and $\taumod$-pure.
Then, after extraction,
$$ g_n|_{\CF(\tau_-)}\to \phi $$
uniformly on compacts
to an open continuous limit map 
$$ \phi: \CF(\tau_-)\to\stF(\tau_+).$$
Moreover,
for every $\hat\tau_-\in C(\tau_-)$
the restriction 
$$\phi|_{\stF(\hat\tau_-)}: \stF(\hat\tau_-)\to\stF(\tau_+)$$
is given by the restriction of an element in $G$,
and hence is an (algebraic) homeomorphism.
\end{prop}
\proof
We fix a base point $o\in X$.

We first note that 
we can replace the sequence $(g_n)$ by a sequence of transvections. 
Indeed, 
let $b_n\to b$ and $b'_n\to b'$ be converging sequences in $G$,
and put $\tilde g_n:=b_ng_nb'_n$.
Then the $\taumod$-pureness of $(g_n)$ is equivalent to the $\taumod$-pureness of $(\tilde g_n)$,
and 
(\ref{eq:contrt}) to locally uniform convergence 
$\tilde g_n|_{C(b'^{-1}\tau_-)}\to b\tau_+$.
Furthermore, the locally uniform convergence 
$g_n|_{\CF(\tau_-)}\to \phi$ translates into $\tilde g_n|_{\CF(b'^{-1}\tau_-)}\to b\phi b'$ 
with limit map $b\phi b':\CF(b'^{-1}\tau_-)\to\stF(b\tau_+)$.
Thus, the assertion holds for $(g_n)$ iff it holds for $(\tilde g_n)$.
Since $(g_n)$ is $\taumod$-pure and we may pass to a subsequence,
we can therefore replace $(g_n)$,
using a $KAK$-decomposition of $G$,
by a sequence of transvections $t_n$ with axes through $o$.
Moreover, we may assume that 
the sequence $(t_no)$ is contained in the Weyl sector $V(o,\tau_+)$ and $\taumod$-regular,
i.e.\ drifts away from the boundary.

In this special situation, things become explicit:
The simplex $\tau_-$ is $o$-opposite to $\tau_+$.
The horocyclic subgroup $H_{\tau_-}$ acts simply transitively on $C(\tau_-)\subset\Flagt$
and, accordingly, 
the natural map 
$$ H_{\tau_-}\times\stF(\tau_+)\lra \CF(\tau_-)$$
is a homeomorphism.
The transvections $t_n$ normalize $H_{\tau_-}$ 
and it holds that\footnote{See e.g.\ \cite[\S 2.17]{Eberlein}.}
$$ c_{t_n}\to e $$
uniformly on compacts in $H_{\tau_-}$,
where $c_g$ denotes conjugation by $g$.
Since $t_n$ acts trivially on $\stF(\tau_+)$,
we obtain for $h\in H_{\tau_-}$ and $\si\in\stF(\tau_+)$ that 
$$ t_n (h\si) = \underbrace{(t_nht_n^{-1})}_{\to e} \underbrace{(t_n\si)}_{=\si} \to\si .$$
Hence 
$$ t_n\to\phi $$
uniformly on compacts in $\CF(\tau_-)$,
with the open continuous limit map $\phi:\CF(\tau_-)\to\stF(\tau_+)$ given by 
$\phi(h\si)=\si$.
Moreover,
for $\hat\tau_-\in C(\tau_-)$, 
the restriction 
$\phi|_{\stF(\hat\tau_-)}:\stF(\hat\tau_-)\to\stF(\tau_+)$ coincides with the restriction of the unique element in $H_{\tau_-}$ 
which maps $\hat\tau_-$ to $\tau_+$.
\qed

\subsection{Dynamics on the Finsler compactification}
\label{sec:derdynrelfi}

We now consider the action 
$$ G\acts\ol X^{Fins} $$
on the Finsler compactification.

\subsubsection{From Finsler to Furstenberg dynamical relations}

We show first that 
dynamical relations in the Finsler boundary 
(with respect to the action on the entire compactification)
imply intrinsic dynamical relations in the Furstenberg boundary:
\begin{lem}
\label{dynrelfuefi}
Let $g_n\to\infty$ be a sequence in $G$.
Suppose that there is a dynamical relation 
$$ \xi\stackrel{(g_n)}{\sim} \xi' $$
in $\ol X^{Fins}$
between boundary points 
$\xi\in X_{\nu}$ with $\nu\in\Flagn$
and $\xi'\in X_{\nu'}$ with $\nu'\in\Flag_{\numod'}$. 
Then, after extraction, there exist $\nu_-\in\Flagin$, $\nu'_-\in\Flag_{\iota\numod'}$
and open continuous maps 
$$ \CF(\nu_-)\stackrel{\phi}{\lra}\stF(\nu)
\qquad\hbox{ and }\qquad
\CF(\nu'_-)\stackrel{\phi'}{\lra}\stF(\nu')$$
such that for every $\tilde\si\in \CF(\nu_-)\cap\CF(\nu'_-)$
there is the dynamical relation
\begin{equation}
\label{eq:intrdynrl}
\phi\tilde\si \stackrel{(g_n)}{\sim} \phi'\tilde\si 
\end{equation}
in $\DF X$.
\end{lem}
Note that $\CF(\nu_-)$ and $\CF(\nu'_-)$ are open dense in $\DF X$,
and hence also their intersection.

The dynamical relation (\ref{eq:intrdynrl}) is meant to hold intrinsically {\em inside} $\DF X$,
i.e. there exists a sequence $(\si_n)$ in $\DF X$, and not just in $\ol X^{Fins}$, such that 
$\si_n\to\phi\tilde\si$ and $g_n\si_n\to\phi'\tilde\si$.

\proof
By assumption, and since $X$ is dense in $\ol X^{Fins}$,
there exists a sequence $(x_n)$ in $X$ such that 
\begin{equation*}
x_n\to\xi 
\qquad\hbox{ and }\qquad 
g_nx_n\to\xi' .
\end{equation*}
We fix a base point $o\in X$ and write 
$$ x_n=a_no
\qquad\hbox{ and }\qquad 
g_nx_n=b_no $$
such that 
$g_n=b_na_n^{-1}$.
Then the sequences $(a_n)$ and $(b_n)$ in $G$ are $\numod$-pure, respectively, $\numod'$-pure 
and we have flag convergence 
$$a_n\to\nu
\qquad\hbox{ and }\qquad
b_n\to\nu' ,$$
see Propositions~\ref{prop:puregstr} and~\ref{prop:relconvflfi}. 
After extraction, we obtain that also the inverse sequences flag converge,
$$a_n^{-1}\to\nu_- 
\qquad\hbox{ and }\qquad
b_n^{-1}\to\nu'_- $$
with $\nu_-\in\Flagin$ and $\nu'_-\in\Flag_{\iota\numod'}$,
cf.\ Lemma~\ref{lem:contrsym}.
The sequences $(a_n)$ and $(b_n)$ are then contracting on the appropriate flag manifolds, 
$$a_n|_{C(\nu_-)}\to\nu
\qquad\hbox{ and }\qquad
b_n|_{C(\nu'_-)}\to\nu' $$
uniformly on compacts.

Due to pureness,
we get more precise information about the accumulation dynamics of these sequences on $\DF X$.
Proposition~\ref{prop:puredyn} yields that 
$$ a_n|_{\CF(\nu_-)} \to \phi
\qquad\hbox{ and }\qquad
b_n|_{\CF(\nu'_-)}\to\phi' $$
uniformly on compacts
with open continuous limit maps 
$$ \phi: \CF(\nu_-)\to\stF(\nu)
\qquad\hbox{ and }\qquad
\phi': \CF(\nu'_-)\to\stF(\nu').$$
Then for $\tilde\si\in\CF(\nu_-)\cap\CF(\nu'_-)$, it holds that 
$$a_n\tilde\si\to\phi\tilde\si
\qquad\hbox{ and }\qquad
g_na_n\tilde\si=b_n\tilde\si\to\phi'\tilde\si ,$$
i.e.\ we obtain the dynamical relation 
$$ \phi\tilde\si \stackrel{(g_n)}{\sim} \phi'\tilde\si $$
inside $\DF X$.
\qed

\begin{rem}[Extrinsic versus intrinsic dynamical relations in $\DF X$]
In the special case $\numod=\numod'=\simod$ this result says 
that all {\em extrinsic} dynamical relations in $\DF X$, as a subset of $\ol X^{Fins}$,
are already {\em intrinsic}.
More precisely,
if $\si,\si'\in\DF X$ and $(x_n)$ is a sequence in $X$ such that 
$x_n\to\si$ and $g_nx_n\to\si'$ 
in $\ol X^{Fins}$,
then there also exists a sequence $(\si_n)$ in $\DF X$ 
such that 
$\si_n\to\si$ and $g_n\si_n\to\si'$.
\end{rem}

We deduce the following consequence from the technical statement in the last lemma:
\begin{cor}
\label{cor:dynrldnspn}
After extraction, the sequence $(g_n)$ satisfies:

If $O\subseteq\stF(\nu)$ and $O'\subseteq\stF(\nu')$ are dense open subsets,
then there exist $\si\in O$ and $\si'\in O'$ which are intrinsically dynamically related in $\DF X$ 
with respect to $(g_n)$.
\end{cor}
\proof
Since $\phi$ is open and continuous, 
the subset $\phi^{-1}(O)$ is dense open in $\CF(\nu_-)$, and hence also in $\DF X$.
Similarly, ${\phi'}^{-1}(O')$ is dense open in $\DF X$. 
Consequently, 
their intersection is nonempty
and contains some $\tilde\si$.
We put $\si=\phi\tilde\si$ and $\si'=\phi'\tilde\si$,
and use the lemma.
\qed

\subsubsection{Dynamical relations}

We can now extend the combinatorial inequality for intrinsic dynamical relations in the Furstenberg boundary 
(Proposition~\ref{prop:dynreltaun})
to the Finsler boundary:

\begin{prop}[Dynamical relation inequality in $\geo^{Fins}X$]
\label{prop:drfinq}
Let $(g_n)$ be a $\taumod$-contracting sequence in $G$ with (\ref{eq:contrt}). 
Suppose that there is a dynamical relation 
\begin{equation*}
\xi\stackrel{(g_n)}{\sim} \xi' 
\end{equation*}
in $\ol X^{Fins}$ 
between boundary points 
$\xi\in X_{\nu}$ with $\nu\in\Flagn$
and $\xi'\in X_{\nu'}$ with $\nu'\in\Flag_{\numod'}$.
Then 
\begin{equation}
\label{eq:dynrelposcondtauncp}
\pos(\si',\tau_+)\prec\cpos(\si,\tau_-) 
\end{equation}
for all 
$\si\in\stF(\nu)$ and $\si'\in\stF(\nu')$.
\end{prop}
\proof
For any two simplices $\tau,\nu\subset\geo X$,
the relative position $\pos(\cdot,\tau)$
has a unique maximal value on $\stF(\nu)$,
i.e.\ all other values are smaller,
and it attains this maximal value on a dense open subset,
cf.\ Proposition~\ref{prop:dsopnstr}.
Let $O\subset\stF(\nu)$ denote the open dense subset where $\pos(\cdot,\tau_-)$ is maximal
and $O'\subset\stF(\nu')$ the subset where $\pos(\cdot,\tau_+)$ is maximal.
By Corollary~\ref{cor:dynrldnspn},
after extraction, 
there exist $\si\in O$ and $\si'\in O'$ 
which are intrinsically dynamically related in $\DF X$,
$$ \si\stackrel{(g_n)}{\sim} \si' .$$
Applying Proposition~\ref{prop:dynreltaun},
we obtain that these $\si$ and $\si'$ satisfy inequality (\ref{eq:dynrelposcondtauncp}).
The inequality for arbitrary 
$\si\in\stF(\tau)$ and $\si'\in\stF(\tau')$ 
follows.\footnote{Here we use the fact that taking complementary position reverses the folding order,
cf.\ (\ref{ineq:revordpos}).}
\qed

\medskip
There are no dynamical relations in $\ol X^{Fins}$ between points in $X$.
This leaves the case of dynamical relations between points in $X$ and points in $\geo^{Fins}X$
which is easy to deal with:
\begin{lem}
\label{lem:drxfbd}
Suppose that the sequence $(g_n)$ in $G$ is 
$\taumod$-contracting with (\ref{eq:contrt}).
If there is a dynamical relation 
$$ x\stackrel{(g_n)}{\sim} \xi' $$
in $\ol X^{Fins}$ between a point $x\in X$ and a boundary point $\xi'\in\geo^{Fins}X$,
then $\xi'\in\ol X_{\tau_+}$.
\end{lem}
\proof
Let $(x_n)$ be any bounded sequence in $X$.
From $g_n\to\tau_+$ it follows that also $g_nx_n\to\tau_+$.
Hence $(g_nx_n)$ accumulates in $\ol X^{Fins}$ at $\ol X_{\tau_+}$,
cf.\ Proposition~\ref{prop:relconvflfi}.
\qed

\medskip
Similarly,
a dynamical relation $\xi\stackrel{(g_n)}{\sim}x'$ 
between a boundary point $\xi\in\geo^{Fins}X$ and a point $x'\in X$ 
implies that $\xi\in\ol X_{\tau_-}$,
as follows by applying the lemma to the inverse sequence $(g_n^{-1})$.

\subsubsection{Accumulation dynamics}

The dynamical relation inequality for the action of contracting sequences 
obtained in the previous section
can be rephrased 
in terms of accumulation at pairs of ``complementary''  thickenings 
at the attractive and repulsive fixed points in $\Flagt$,
compare the discussion in \cite[\S 5.2,6.1]{coco15} for dynamics on flag manifolds.

We refer the reader to section~\ref{sec:thickenings} for the definitions of Furstenberg and Finsler thickenings. 
Lemma~\ref{lem:drxfbd} implies:
\begin{cor}
\label{cor:drcontrsq}
Let $(g_n)$ be a $\taumod$-contracting sequence in $G$ with (\ref{eq:contrt}),
and let $\emptyset\neq\Th\subsetneq W$ be a $\Wt$-left invariant thickening.
Suppose that there is a dynamical relation 
\begin{equation*}
\ol x\stackrel{(g_n)}{\sim} \ol x' 
\end{equation*}
in $\ol X^{Fins}$ 
between points $\ol x,\ol x'\in\ol X^{Fins}$.
Then $\ol x\in (\Th^c)^{Fins}(\tau_-)$ 
or $\ol x' \in\Th^{Fins}(\tau_+)$.
\end{cor}
\proof
By our assumption on $\Th$,
we have that 
$\stF(\tau_-)\subset\ThF^c(\tau_-)$ and $\stF(\tau_+)\subset\ThF(\tau_+)$,
equivalently, 
that 
$\ol X_{\tau_-}\subset(\Th^c)^{Fins}(\tau_-)$ and $\ol X_{\tau_+}\subset\Th^{Fins}(\tau_+)$.
Lemma~\ref{lem:drxfbd} therefore implies the assertion if one of the points $\ol x,\ol x'$ lies in $X$.

We may therefore assume that $\ol x,\ol x'\in\geo^{Fins}X$
and apply Proposition~\ref{prop:drfinq}.
Let $\ol x\in X_{\nu}$ with $\nu\in\Flagn$
and $\ol x'\in X_{\nu'}$ with $\nu'\in\Flag_{\numod'}$.
Either $\pos(\cdot,\tau_+)$ takes values in $\Th$ on $\stF(\nu')$,
or the maximal value of $\pos(\cdot,\tau_+)$ on $\stF(\nu')$ is contained in 
$W-\Th=w_0\cTh$. 
In the latter case,
inequality (\ref{eq:dynrelposcondtauncp}) implies that 
the values of $\pos(\cdot,\tau_-)$ on $\stF(\nu)$ are contained in $\cTh$.

If $\pos(\cdot,\tau_+)|_{\stF(\nu')}$ takes values in $\Th$,
then $\stF(\nu')\subset\ThF(\tau_+)$,
equivalently,
$\ol X_{\nu'}\subset\Th^{Fins}(\tau_+)$,
and hence $\bar x'\in \Th^{Fins}(\tau_+)$.
Similarly,
if $\pos(\cdot,\tau_-)|_{\stF(\nu)}$ takes values in $\cTh$,
then $\ol x\in\ol X_{\nu}\subset(\Th^c)^{Fins}(\tau_-)$.
\qed

\medskip
In the language of accumulation dynamics introduced in \cite[\S 5.2]{coco15}, 
this means that the sequence $(g_n)$ is 
$((\Th^c)^{Fins}(\tau_-),\Th^{Fins}(\tau_+))$-accumulating on $\ol X^{Fins}$.

\subsubsection{Domains of proper discontinuity}

We now deduce our main results on proper discontinuity, 
compare 
\cite[\S 6.2-4]{coco15}.

For a discrete subgroup $\Ga<G$,
we define the {\em forward/backward $\taumod$-limit set} 
$$\Lat^{\pm}\subset\Flagpmt$$ 
as the set of all simplices $\tau_{\pm}$ as in (\ref{eq:contrtau}) 
for all $\taumod$-contracting sequences $(g_n=)\ga_n\to\infty$ in $\Ga$,
cf.\ \cite[Def.\ 6.9]{coco15}.
The limit sets $\Lat^{\pm}$ are $\Ga$-invariant and compact.
If $\taumod$ is $\iota$-invariant, then $\Lat^+=\Lat^-=:\Lat$.

Consider first the case of $\taumod$-regular, equivalently, $\taumod$-convergence subgroups. 
Here we obtain domains by removing suitable thickenings of the $\taumod$-limit set:
\begin{thm}[Domains of proper discontinuity for $\taumod$-convergence subgroups]
\label{thm:pdwconv}
Let $\Ga<G$ be a $\taumod$-con\-ver\-gen\-ce subgroup,
and let $\emptyset\neq\Th\subsetneq W$ be a $\Wt$-left invariant thickening.
Then the action 
\begin{equation}
\label{eq:ddctr}
\Ga\acts\ol X^{Fins}-\bigl((\Th^c)^{Fins}(\Latm)\cup\Th^{Fins}(\Latp)\bigr)
\end{equation}
is properly discontinuous. 
In particular,
if $\taumod$ is $\iota$-invariant and $\Th$ is fat,
then the action 
$$ \Ga\acts\ol X^{Fins}-\Th^{Fins}(\Lat) $$
is properly discontinuous. 
\end{thm}
\proof
Suppose that there is a dynamical relation
\begin{equation}
\label{eq:drfin}
\ol x\stackrel{(\ga_n)}{\sim} \ol x'
\end{equation}
in $\ol X^{Fins}$
with respect to a sequence $\ga_n\to\infty$ in $\Ga$.
After extraction we may assume that $(\ga_n)$ is $\taumod$-contracting with 
\begin{equation*}
\ga_n|_{C(\tau_-)}\to\tau_+
\end{equation*}
uniformly on compacts,
where $\tau_{\pm}\in\Lapmt$.
Corollary~\ref{cor:drcontrsq} implies that 
$\ol x\in (\Th^c)^{Fins}(\tau_-)\subset(\Th^c)^{Fins}(\Latm)$ 
or $\ol x' \in\Th^{Fins}(\tau_+)\subset\Th^{Fins}(\Latp)$.
This yields the first assertion.

The second assertion follows because $\Th^c\subseteq\Th$
due to fatness.\footnote{Here we use that 
if $\Th_1\subseteq\Th_2$, then $\Th_1^{Fins}\subseteq\Th_2^{Fins}$.}
\qed

\medskip
Note that the thickenings of limit sets $\Th^{Fins}(\Lat^{\pm}(\Ga))$ are $\Ga$-invariant and compact. 

This scheme of constructing domains of proper discontinuity 
applies equally well to arbitrary discrete subgroups $\Ga<G$,
compare the discussion in \cite[\S 6.4]{coco15}.
One then has to take into account the $\taumod$-limit sets for all face types $\taumod$.
There are several ways to proceed.
The most immediate family of possibilities is the following. 

\begin{thm}[Domains of proper discontinuity for discrete subgroups I]
\label{thm:pdwdsc}
Let $\Ga<G$ be a discrete subgroup,
and let $\emptyset\neq\Th_{\taumod}\subsetneq W$ be $\Wt$-left invariant thickenings
for all face types $\taumod\subseteq\simod$.
Then the action 
\begin{equation}
\label{eq:dompdlar}
\Ga\acts\ol X^{Fins}-\bigcup_{\taumod}\bigl((\Th_{\taumod}^c)^{Fins}(\Latm)\cup\Th_{\taumod}^{Fins}(\Latp)\bigr)
\end{equation}
is properly discontinuous. 
\end{thm}
\proof
The proof is the same as for the previous theorem:
Suppose that there is a dynamical relation (\ref{eq:drfin}) in $\ol X^{Fins}$. 
Then $(\ga_n)$ contains for some face type $\taumod$ 
a $\taumod$-contracting subsequence 
and it follows as before that 
$\ol x\in (\Th_{\taumod}^c)^{Fins}(\Latm)$ or $\ol x' \in\Th_{\taumod}^{Fins}(\Latp)$.
\qed

\medskip
In general, these domains of proper discontinuity can be further enlarged
by only removing the thickenings of the limit simplices arising from {\em pure} sequences in the group: 
Define the {\em pure forward/backward $\taumod$-limit set} 
$$\Lat^{pure,\pm}\subseteq\Lat^{\pm}$$
as the closure of the set 
of all simplices $\tau_{\pm}$ as in (\ref{eq:contrtau}) 
for all $\taumod$-pure $\taumod$-contracting sequences $(\ga_n)$ in $\Ga$.
As above, we conclude:

\begin{thm}[Domains of proper discontinuity for discrete subgroups II]
\label{thm:pdwdsc2}
Let $\Ga<G$ be a discrete subgroup,
and let $\emptyset\neq\Th_{\taumod}\subsetneq W$ be $\Wt$-left invariant thickenings
for all face types $\taumod\subseteq\simod$.
Then the action 
\begin{equation}
\label{eq:dompd}
\Ga\acts\ol X^{Fins}-\bigcup_{\taumod}\bigl((\Th_{\taumod}^c)^{Fins}(\Lat^{pure,-})\cup\Th_{\taumod}^{Fins}(\Lat^{pure,+})\bigr)
\end{equation}
is properly discontinuous. 
\end{thm}
\proof
Same argument as before,
taking into account that every sequence $\ga_n\to\infty$ in $\Ga$ contains a subsequence
which for some face type $\taumod$ is $\taumod$-contracting and $\taumod$-pure.
\qed

\medskip
Since the domain in (\ref{eq:dompd}) is larger than the domain in (\ref{eq:dompdlar}),
one can in general not expect the $\Ga$-action (\ref{eq:dompdlar}) to be cocompact.

If $\Ga$ is $\taumod$-regular, 
then it contains $\numod$-pure sequences only for the face types $\numod\supseteq\taumod$,
and hence only these limit sets $\Lan^{\pm}$ can be nonempty.
Since $\Wn\leq \Wt$, 
we may choose $\Th_{\numod}=\Th_{\taumod}$ for these face types,
and then the domain in (\ref{eq:dompdlar}) coincides with the domain in Theorem~\ref{thm:pdwconv}.

\subsubsection{Nonemptyness of domains of proper discontinuity at infinity}

If we assume in addition to the hypotheses of Theorem~\ref{thm:pdwconv}
that $\Ga$ is $\taumod$-antipodal,
then it is easy to see that, in higher rank, the domains (\ref{eq:ddctr})
strictly enlarge $X$:
\begin{prop}
\label{prop:connected}
Let $\Ga<G$ be a $\taumod$-antipodal $\taumod$-convergence subgroup and 
let $\Th\subset W$ be a $\Wt$-left invariant slim thickening.
If $\rank(X)\geq2$, then $\Th^{Fins}(\Lat)\subsetneq \geo^{Fins} X$.
\end{prop}
\proof 
If $\Th=\emptyset$, there is nothing to show.
Suppose therefore that $\Th\neq\emptyset$.

We consider the subcomplex $C$ of $\amod$
corresponding to the thickening $\ThF(\taumod)$ of the model chamber.
Since $\rank(X)\geq2$, $\amod$ is connected.
By slimness, $C$ does not contain all chambers of $\amod$.
Therefore there exists a panel $\bar\pi\subset C$ 
such that exactly one of two chambers in $\amod$ adjacent to it belongs to $C$.
In terms of chambers, this means that 
$\stF(\bar\pi)\cap\ThF(\taumod)\neq\emptyset$ and 
$\stF(\bar\pi)\not\subset\ThF(\taumod)$.

Let $\tau_0\in\Lat$,
and let $\pi\subset\geo X$ be a panel with $\pos(\pi,\tau_0)=\bar\pi$.
Then $\stF(\pi)$ intersects $\ThF(\tau_0)$, but is not contained in it.
It follows that $\stF(\pi)\not\subset\ThF(\tau)$ for all $\tau\in\Lat$,
because due to our assumptions of antipodality and slimness,
the thickenings $\ThF(\tau)$ for $\tau\in\Lat$ are pairwise disjoint. 
Consequently, $X_{\pi}\not\subset\Th^{Fins}(\Lat)$.
\qed

\begin{rem}
Note that nonemptyness of domains of proper discontinuity 
for the $\Ga$-actions on flag manifolds
is much harder to prove and requires additional assumptions.
For the case of actions on the Furstenberg boundary, see \cite[\S 8]{coco15}. 
\end{rem}

\section{General cocompactness results}
\label{sec:gencoco}

\subsection{General discrete topological group actions} 

The main result of this section is a cocompactness theorem for a certain class of properly discontinuous group actions $\Ga\acts \Om$, where $\Om$ is an open subset of a compact metrizable space $Z$ and the action 
extends to a topological action $$\Ga\acts Z.$$
In order to prove cocompactness, we need to impose certain assumptions on both $\Ga$ and the action. 
We assume that there exists a {\em model action} 
$$\Ga\acts Y$$ 
on a contractible simplicial complex
which is simplicial, properly discontinuous and cocompact. 
We further assume that $Y$ admits a $\Ga$-equivariant contractible metrizable compactification $\ol Y\supset Y$. 
The extended action $\Ga \acts \ol{Y}$ serves as a model for $\Ga\acts Z$. 
Examples of such model actions abound in geometric group theory. For instance, if $\Ga$ is Gromov-hyperbolic, we can take for $Y$ a suitable Rips complex of $\Ga$ and for $\ol{Y}$ the Gromov compactification of $Y$. Other examples are given by isometric properly discontinuous cocompact actions $\Ga\acts Y$ on  piecewise-Riemannian $CAT(0)$ complexes $Y$, where $\ol{Y}$ is the visual compactification of $Y$. 

Our next set of hypotheses relates the model action $\Ga\acts Y$ to the action $\Ga\acts Z$. 
We assume that there exists a $\Ga$-equivariant continuous map of triads
\begin{equation}
\label{eq:triads}
(\ol{Y}, Y, \underbrace{\ol Y - Y}_{=:\La_{mod}}) \stackrel{\widetilde f}{\lra} (Z, \Om,  \underbrace{Z -\Om}_{=:\La}) 
\end{equation}
such that the restriction $\widetilde f|_{\La_{mod}}:\La_{mod}\to\La$
is a $\Z_2$-\v{C}ech cohomology equivalence. 
 
\begin{thm}\label{thm:extended_cocompactness}
Under the above assumptions, 
and if $\Ga$ is torsion-free and $\Om$ is path connected, then $\Om/\Ga$ is compact. 
\end{thm}
\proof \medskip 
{\em Step 1: 
Passing to a model action on a manifold with boundary.} 
We replace the action $\Ga\acts Y$ with an action of $\Ga$ on a suitable manifold. 
Since $\pi_1(Y)=1$ and the action $\Ga\acts Y$ is free, the quotient space  $R:= Y/\Ga$ satisfies $\pi_1(R)\cong \Ga$. 

We ``thicken" $R$ to a closed manifold without changing the fundamental group.
To do so, we first embed $R$ as a subcomplex into the (suitably triangulated) euclidean space $E^{2n+1}$, 
where $n=\dim(R)$.  We denote by $N$ 
the regular neighborhood of $R$ in $E^{2n+1}$, and let $D=\D N$.

\begin{lem}\label{lem:RN}
$D$ is connected and $\pi_1(D)\to \pi_1(N)\cong \pi_1(R)$ is surjective.  
\end{lem}
\proof Let $N':= N - R$. We claim that the map $D\embed N'$ is a homotopy equivalence. The proof is the same as the one for the homotopy equivalence $R\to N$: Each simplex $c\subset N$ is the join $c_1 \star c_2$
of a simplex $c_1$ disjoint from $R$ (and, hence, contained in $D$) and a simplex $c_2\subset R$ (in the extreme cases, $c_1$ or $c_2$ could be empty). Now, use the straight line segments 
given by these join decompositions to homotop each $c- R$ to $c_1\subset D$. 

Since $R$ has codimension $\geq2$ in $N$, it does not separate $N'$ and each loop in $N$ is homotopic to a loop in $N'$. 
Hence, $N'$ is connected and 
$$\pi_1(D)\stackrel{\cong}{\lra} \pi_1(N')\lra \pi_1(N)$$
is surjective. 
\qed

\begin{lem}\label{lem:key}
There exists a connected closed manifold $M$ which admits a map $h: R\to M$ inducing an isomorphism 
of fundamental groups $\pi_1(R)\to \pi_1(M)$.   
\end{lem}
\proof We start with $N$ (the regular neighborhood of $R\subset E^{2n+1}$) as above. 
As noted in the proof of the previous lemma,
the inclusion $R\to N$ is a homotopy equivalence, and $N$ is a compact manifold with boundary. 
Consider the closed manifold $M$ obtained by doubling $N$ along its boundary $D$,
$$
M= N_1\cup_D N_2,
$$
where $N_1, N_2$ are two copies of $N$. 
We let $i: D\to M, i_k: N_k\to M$ denote the inclusion maps. Since $M$ is the double of $N$, 
we have the canonical retraction $r: M\to N_1$ (whose restriction to $N_2$ 
is a homeomorphism given by reflecting at $D$). Define the map $h=i_1\circ g$, 
$$
h: R \stackrel{g}{\lra} N_1\stackrel{i_1}{\lra} M,$$
where $g$ corresponds to the inclusion $R\to N$ and hence is a homotopy equivalence. 
We claim that $h$ induces an isomorphism $h_*$ of fundamental groups. 

The existence of the retraction $r$ implies the injectivity of $i_{1*}$ and hence of $h_*$. 

By Lemma \ref{lem:RN}, $D$ is connected.
Hence, the Seifert--van Kampen theorem implies that $\pi_1(M)$ is generated by the two 
subgroups $i_{k*}(\pi_1(N_k)), k=1,2$.  Since the homomorphisms
$$\pi_1(D)\to \pi_1(N_k)$$
are surjective (Lemma \ref{lem:RN}), we obtain 
$$
i_{1*}(\pi_1(N_1))= i_*(\pi_1(D))=  i_{2*}(\pi_1(N_2)). 
$$ 
Hence, both homomorphisms $i_{k*}: \pi_1(N_k) \to \pi_1(M)$ are surjective. 
The surjectivity of $h_*$   follows.   \qed

\medskip 
We let $m=2n+1$ denote the dimension of the manifold $M$ and its  universal cover $\widetilde{M}$.

\medskip
{\em Step 2.} 
We let $M$ be a manifold as in Lemma \ref{lem:key}. We consider the triads (\ref{eq:triads}) and the diagonal $\Ga$-action on their products with the universal cover $\widetilde M$. 
Dividing out the action, we obtain bundles over $M$ and $\tilde f$ induces the map
of triads of bundles
\begin{equation}
\label{eq:triadsbun} 
(\underbrace{(\ol{Y}\times\widetilde M)/\Ga}_{\bar E_{mod}},
\underbrace{(Y \times\widetilde M)/\Ga}_{E_{mod}},
\underbrace{(\La_{mod}\times\widetilde M)/\Ga}_{L_{mod}})
\stackrel{F}{\lra} 
(\underbrace{(\Si\times\widetilde M)/\Ga}_{\bar E},
\underbrace{(\Om\times\widetilde M)/\Ga}_{E},
\underbrace{(\La\times\widetilde M)/\Ga}_{L})
\end{equation}
Note that $E$ also fibers over $\Om/\Ga$ with fiber $\widetilde M$.

The map $F$ of triads of bundles satisfies: 

(i) $F|_{E_{mod}}: E_{mod}\to E$ is proper. 

(ii) $F|_{L_{mod}}: L_{mod}\to L$ is a cohomology equivalence of bundles. 

\begin{lem}
Both spaces $\bar E, \bar E_{mod}$ are metrizable. 
\end{lem}
\proof These spaces are fiber bundles with compact metrizable bases and fibers. 
Therefore, $\bar E, \bar E_{mod}$ are  both compact and Hausdorff. 
Hence, they are metrizable, for instance, by Smirnov's metrization theorem, 
because they are paracompact, Hausdorff and locally metrizable. \qed 

\medskip 
Our approach to proving Theorem~\ref{thm:extended_cocompactness} is based on the following observation.

In a noncompact connected manifold, the point represents the zero class in $H^{lf}_0$. 
Similarly, let $\iota: F\to E\stackrel{\pi}{\to} B$ be a fiber bundle over a noncompact connected manifold, where $\iota: F\to E_b$ is the homeomorphism of $F$ to the fiber $E_b=\pi^{-1}(b)$. If the base $B$ is path-connected, then the induced 
map 
$$\iota_*: H^{lf}_*(F)\to H^{lf}_*(E) $$
on locally finite homology is independent of the choice of $b$. In order to show triviality of this map provided that $B$ is noncompact, note that for each class $[\eta]\in Z^{i}_c(E)$ and each locally finite class $[\xi] \in H^{lf}_m(F)$, if $b$ is chosen so that $E_b$ is disjoint from the support of $\eta$, then $\<[\eta], [\xi]\>=0$.   
Here and in the sequel we use (co)homology with $\Z_2$-coefficients. Hence, $\iota_*=0$.

The compactness of $\Om/\Ga$ therefore follows from showing that the fiber of the bundle 
$$\widetilde M^m\to E\to\Om/\Ga$$
represents a nontrivial class in $H^{lf}_m(E)$, 
i.e. that the locally finite fundamental class $[\widetilde M]\in H^{lf}_m(\widetilde M)$
has a non-zero image under the inclusion induced map 
$H^{lf}_m(\widetilde M)\to H^{lf}_m(E)$. 

The proper map $F: E_{mod}\to E$ induces the map 
$F_*: H^{lf}_m(E_{mod})\to H^{lf}_m(E)$
which carries the class represented by the $\widetilde M$-fiber in the model $E_{mod}$ to the 
corresponding class in $E$.
It therefore suffices to show that 
\begin{equation}\label{eq:compose}
\underbrace{H^{lf}_m(\widetilde M)}_{\cong\Z_2} \stackrel{\iota_*}{\lra} H^{lf}_m(E_{mod})\stackrel{F_*}{\lra} H^{lf}_m(E)
\end{equation}
is a composition of injective maps.

{\em Step 3: Injectivity of $F_*$.}
We pass to compactly supported cohomology. 
We recall that locally finite homology (with field coefficients) 
is dual to compactly supported cohomology in the same degree via Kronecker duality.
We therefore must show that the dual map 
$$ H^m_c(E) \stackrel{F^*}{\lra} H^m_c(E_{mod}) $$
is surjective. 

We now switch the fiber direction
and regard $E$ and $E_{mod}$ as bundles over $M$.
We use their compactifications $\bar E$ and $\bar E_{mod}$ mentioned earlier
which allow us to replace compactly supported cohomology by relative cohomology. 
Since $E$ is compact and metrizable, while $L$ is compact, we have a natural isomorphism of Alexander-Spanier cohomology groups (cf. \cite[Lemma 11, p. 321]{Spanier}):
$$
H^m_c(E)\cong  {H}^m(\bar E, {L})
$$
Similarly, we have a natural isomorphism
$$H^m_c(E_{mod})\cong H^m(\bar E_{mod},L_{mod}).$$
Thus, the surjectivity of the previous map $F^*$ is equivalent to the surjectivity of the map
$$ 
{H}^m(\bar E,L) \stackrel{F^*_{rel}}{\lra} {H}^m(\bar E_{mod},L_{mod})$$
induced by the map of pairs
\begin{equation}\label{eq:pair}
(\bar E_{mod},L_{mod})\stackrel{F}{\lra} (\bar E,L) . 
\end{equation}
To verify the surjectivity of $F^*_{rel}$, we use 
the long exact cohomology sequence of $F$:
\begin{diagram}
\ldots &  {H}^{m-1}(\bar E) & \rTo &  {H}^{m-1}(L) &\rTo & {H}^m(\bar E, L) & \rTo &  {H}^m(\bar E)  & \rTo &  {H}^m(L)       & \ldots\\
          & \dTo                   &   ~  & \dTo~\cong &        &\dTo^{F^*_{rel}}   &        &\dTo^{F^*_{abs}}&        &\dTo~\cong  &         \\
\ldots &  {H}^{m-1}(\bar E_{mod}) & \rTo &  {H}^{m-1}(L_{mod}) &\rTo & {H}^m(\bar E_{mod}, L_{mod}) & \rTo^j &  {H}^m(\bar E_{mod})  & \rTo &  {H}^m(L_{mod})       & \ldots\\
\end{diagram}
A diagram chase (as in the proof of the 5-lemma) shows that the surjectivity of $F^*_{rel}$ follows from the surjectivity of $F^*_{abs}$.
Indeed, one first checks that $\ker j\subset\im F^*_{rel}$, and uses this to verify the inclusion 
$$ j^{-1}(\im F^*_{abs}) \subset \im(F^*_{rel}) .$$
To see that $F^*_{abs}$ is surjective,
we consider the map of bundles:
\begin{diagram}
\bar E_{mod} &                                    & \rTo^F &                                     & \bar E  \\
               & \rdTo^{\pi_{\bar E_{mod}}} &            & \ldTo^{\pi_{\bar E}}      &             \\
               &                                    &    M     &                                     &             \\
\end{diagram}
The fibration $\pi_{\bar E_{mod}}$ is a homotopy equivalence because its fibers $\ol{Y}$ are contractible. Let 
$$
s: M\to \bar E_{mod}
$$
denote a section. It follows that $s\circ\pi_{\bar E}$ is a left homotopy inverse for $F$,
i.e. $s\circ\pi_{\bar E}\circ F\simeq\id_{\bar E_{mod}}$. 
Thus, the induced map on cohomology $F^*_{abs}$ is surjective. 

{\em Step 4: Injectivity of $\iota_*$.}
We consider the fiber bundle 
$$
\widetilde M\to E_{mod}\to R.$$ 
The map $h: R\to M$ in Lemma \ref{lem:key} yields a section of this bundle. Since the base $R$ of the bundle is a finite CW complex and its fiber $\widetilde M$ is a connected $m$-manifold,
Lemma \ref{lem:thom} implies that the induced map 
$$
H^{lf}_m(\widetilde M) \stackrel{\iota_*}{\lra} H^{lf}_m(E_{mod}),
$$
is injective. 

This concludes the proof of Theorem \ref{thm:extended_cocompactness}. \qed

\subsection{Ha\"issinsky-Tukia conjecture for convergence actions} \label{sec:haiss}

We now apply our general cocompactness result from the previous section
(Theorem~\ref{thm:extended_cocompactness}) to the theory of abstract convergence groups.\footnote{See 
e.g.\ \cite{Bowditch} or \cite{Tukia1} for background on convergence groups.}
The following natural question is due to P.\ Ha\"issinsky \cite{Haissinsky}.
An equivalent question was asked by P.~Tukia in \cite[p. 77]{Tukia2}.
We owe the observation of the equivalence of the questions to V.~Gerasimov.  

\begin{ques}\label{ques:Haissinsky}
Let $\Ga\acts \Si$ be a convergence group action of a hyperbolic group on a metrizable compact space, and suppose that $\La\subset \Si$ is an invariant compact subset which is equivariantly homeomorphic to $\geo\Ga$. Then the action $ \Ga\acts \Om= \Si  -\La$ is properly discontinuous.
Is it always cocompact? 
\end{ques}

\begin{rem}
This is true for actions which are {\em expanding} at the limit set $\La$, cf. \cite{coco15}.
\end{rem}

The main result of this section is the following theorem which provides strong evidence for a positive answer to Question \ref{ques:Haissinsky} in the case of convergence group actions with path-connected discontinuity domains. 

\begin{thm}\label{thm:haiss}
Let $\Ga\acts \Si$ be a convergence group action of a virtually torsion-free  hyperbolic group on a metrizable compact space $\Si$, and suppose that $\La\subset \Si$ is an invariant compact subset which is equivariantly homeomorphic to $\geo\Ga$. 
Then the action $$ \Ga\acts \Si - \La$$ is cocompact provided that $\Si - \La$ has finitely many path connected components.
\end{thm}
\proof 
{\em Step 1.} After passing to a finite index subgroup of $\Ga$ preserving each connected component of 
$\Om':=\Si - \La$, it suffices to consider the case when $\Om'$ is path connected (and nonempty). It also suffices to consider the case when $\Ga$ is torsion-free. We let $Y$ be a contractible locally compact simplicial complex on which $\Ga$ acts properly discontinuously and cocompactly, e.g.\ a suitable Rips complex of $\Ga$. The Gromov compactification  $\ol{Y}$ of $Y$ 
is  contractible and metrizable, cf. \cite{Bestvina-Mess}, 
and $\ol{Y} - Y\cong\geo \Ga$ equivariantly.

\medskip 
{\em Step 2: Construction of a map of triads.}
Pick a point $x\in \Om'$ and define the orbit map
 $$
f: \Ga\to \Om', \quad \ga\mapsto \ga x.
 $$ 
 This map is injective since $\Ga$ is torsion-free and, hence, acts freely on $\Om'$.  Let 
 $f_\infty: \geo \Ga \to \La$ be an (the) equivariant homeomorphism. 
 We further let $\ol{\Ga}=\Ga \cup \geo \Ga$ denote the Gromov compactification of $\Ga$. 
We define the map
$$
\ol f: \ol\Ga \lra \Si, 
$$
whose restriction to $\Ga$ is $f$ 
and to $\geo\Ga$ is $f_{\infty}$.

\begin{lemma}\label{lem:conext}
$\ol f$ is an equivariant homeomorphism onto $\Ga x \cup \La$.  
\end{lemma} 
\proof 
We first note that the natural action $\Ga\acts \ol{\Ga}$ is a convergence action.

Suppose that $(\ga_n)$ is a sequence in $\Ga$ converging to $\xi\in \geo \Ga$; $\la=f_\infty(\xi)$. We claim that 
$$
\lim_{n\to\infty} f(\ga_n)= \la. 
$$

{\em Case 1: $\Ga$ is nonelementary.}  
Without loss of generality (in view of compactness of $\Si$ and the convergence property of the action $\Ga\acts \Si$), there exists $\la_-\in \La$ such that the sequence $\ga_n|_{\Si-\{\la_-\}}$ converges to some $\la_+\in \La$ uniformly on compacts. Since $f_\infty$ is a homeomorphism, $\ga_n$ converges to $f^{-1}_\infty(\la_+)$ uniformly on compacts in $\geo \Ga - f_{\infty}^{-1}(\la_-)$. The assumption that  $\Ga$ is nonelementary implies that 
$\geo \Ga - f_{\infty}^{-1}(\la_-)$ consists of more than one point. Therefore, in view of the convergence property for the action $\Ga\acts \ol\Ga$, it follows that  $\ga_n$ converges to $f^{-1}_\infty(\la_+)$ on 
${\Ga}$ (here we again pass to a subsequence if necessary).  
Hence, $\xi=f^{-1}_\infty(\la_+)$, $\la_+=\la$ and the continuity of $\ol f$ follows 
(cf. Lemma \ref{lem:topologylemma}). 

{\em Case 2: $\Ga$ is elementary, i.e, $\Ga\cong \Z$.} Then $\Ga$ is generated by a single {\em loxodromic} homeomorphism $\ga: \Si \to \Si$, i.e., $\La=\{\la_+, \la_-\}$. 
Tukia proved \cite[Lemma 2D]{Tukia1} that the sequence $(\ga^n)$ converges uniformly on compacts in 
$\Si - \{\la_-\}$ to $\la_+$, while the sequence $(\ga^{-n})$ converges uniformly on compacts in 
$\Si - \{\la_+\}$ to $\la_-$. This implies continuity of the map $\ol f$. \qed 

\medskip
We now amalgamate the spaces $\ol{Y}$ and $\Si$ using the homeomorphism 
$$
\ol f: \ol\Ga= \Ga\cup \geo \Ga\to \ol{\Ga x}= \Ga x \cup \La,$$
where we identify $\Ga$ with a subset of the vertex set of the complex $Y$. 
We denote by  $Z$ the result of the amalgamation. 
This space is metrizable by Urysohn's metrization theorem,
since it is Hausdorff, compact and 1st countable,
compare also Proposition~\ref{prop:metrization} which provides a different proof. 

Since $\ol f$ is $\Ga$-equivariant, the  topological action of $\Ga$ on $\ol{Y} \sqcup \Si$ 
descends to a topological action $\Ga \acts Z$. This action is properly discontinuous on 
$\Om:= Y \cup \Om'\subset \Si$ as for each compact $C\subset \Om$, its intersections with $Y$ and $\Om'$ are both compact and the actions $\Ga\acts Y, \Ga\acts \Om'$ are properly discontinuous. 
Lastly, we note that, in view of connectivity of $Y$,  since $\Om'$ is path connected, so is $\Om$. 
Since the embedding $\Om' \to \Om$ is proper, $\Om/\Ga$ is compact if and only if $\Om'/\Ga$ is compact. 
We let $\tilde f: \ol{Y}\to Z$ be the inclusion map. 

\medskip 
{\em Step 3.} 
 According to Theorem \ref{thm:extended_cocompactness}, $\Om/\Ga$ is compact. 
 Therefore, $\Om'/\Ga$ is compact as well. \qed

\begin{rem}
It is not hard to check that  $\Ga\acts Z$ is a convergence action, however,  
this is not needed for our argument.
\end{rem}

\section{Cocompactness}

We return to the discussion of discrete subgroups of Lie groups 
and their actions on Finsler compactifications. 
In section~\ref{sec:pd},
we constructed domains of proper discontinuity.
We will now prove the cocompactness of these actions
for certain classes of discrete subgroups. 

Let $\taumod$ be $\iota$-invariant. 
In sections~\ref{sec:reg} and~\ref{sec:conv} 
we defined $\taumod$-regular and $\taumod$-antipodal discrete subgroups $\Ga< G$ and the $\taumod$-limit set 
$\Lat \subset \Flagt$. 
In \cite{morse} we defined the following class of $\taumod$-antipodal $\taumod$-regular subgroups:

\begin{dfn}[Asymptotically embedded]
\label{def:asyemb}
A discrete subgroup $\Ga<G$ is {\em $\taumod$-asymp\-to\-ti\-cally embedded} 
if it is $\taumod$-regular, $\taumod$-antipodal,
intrinsically word hyperbolic 
and there is a $\Ga$-equivariant homeomorphism 
\begin{equation*}
\label{eq:mapalphatauintro}
\alpha: \geo \Ga \stackrel{\cong}{\lra}
\Lat\subset \Flagt
\end{equation*}
from its Gromov boundary onto its $\tau_{mod}$-limit set. 
\end{dfn}
We proved in \cite{morse} that a subgroup $\Ga< G$ is $\taumod$-asymptotically embedded iff 
it is $\taumod$-Anosov\footnote{I.e. $P_{\taumod}$-Anosov in the terminology of \cite{GW}.}. 

Suppose now that $\Ga< G$ is $\taumod$-regular and $\taumod$-antipodal,
and that $\Th\subset W$ is a $W_{\taumod}$-invariant balanced thickening. 
In this section, we will use the following notation:
$$\hat\Si:=\ol X^{Fins},\quad \hat\La:=\Th^{Fins}(\Lat),\quad \hat\Om:=\hat\Si-\hat\La $$
According to Theorem~\ref{thm:pdwconv},
the action 
$\Ga\acts\hat\Om$
is properly discontinuous. 
We will show that it is also cocompact
provided that $\Ga$ is $\taumod$-asymptotically embedded, 
by replacing the action $\Ga\acts \hat\Si$ with a convergence action $\Ga\acts \Si$ on a certain quotient space of $\hat\Si$ and then applying our cocompactness result for convergence actions (Theorem~\ref{thm:haiss}). 
The collapse takes place only in the thickening $\hat\La$ at infinity,
so that the action $\Ga\acts\hat\Om$ is not affected.

\subsection{Decompositions and collapses}

A {\em decomposition} ${\mathcal R}$ of a set $Z$ is an equivalence relation on $Z$. We let 
${\mathcal D}={\mathcal D}_{\mathcal R}$ denote the subset of the power set $2^Z$ consisting of the equivalence classes of ${\mathcal R}$. 

A decomposition of a Hausdorff topological space $Z$ is {\em closed} if the elements of ${\mathcal D}$ are closed subsets of $Z$; a decomposition is {\em compact} if its elements are compact subsets. Given a decomposition ${\mathcal R}$ of $Z$, one defines the quotient space 
$Z/{\mathcal R}$. Quotient spaces of closed decompositions are $T_1$ but in general not Hausdorff. 
\begin{definition}
A decomposition of $Z$ is {\em upper semicontinuous} (usc) if it is closed 
and for each $D\in {\mathcal D}$ and each open subset $U\subset Z$ containing $D$, there exists another open subset $V\subset Z$ containing $D$ such that  every $D'\in {\mathcal D}$ intersecting $V$ nontrivially is already contained in $U$. 
\end{definition}

\begin{lem}
[{\cite[Proposition 1, page 8]{Daverman}}]
\label{lem:uschar}
The following are equivalent for a closed decomposition ${\mathcal R}$ of $Z$:

(i) ${\mathcal R}$ is usc.

(ii) For every open subset $U\subset Z$,
the saturated subset 
$$U^*=\bigcup\{D\in {\mathcal D}:D\subset U\} $$
is open.

(iii) The quotient projection
$$Z\stackrel{\kappa}{\lra} Z/{\mathcal R}$$ is closed. 
\end{lem}
\proof
(i)$\Ra$(ii):
Let $x\in U$ and let $D\in{\mathcal D}$ be the decomposition subset through $x$.
The usc property implies that $U^*$ contains a neighborhood of $x$.

(ii)$\Ra$(i): Take $V=U^*$.

(ii)$\Ra$(iii):
Let $C\subset Z$ be closed, and let $U$ be the complement. 
Then $U^*=\kappa^{-1}\kappa(Z-C)$, and it follows that $\kappa(C)$ is closed.

(iii)$\Ra$(ii):
Let $U\subset Z$ be open. 
Then $U^*=\kappa^{-1}(Z/{\mathcal R}-\kappa(Z-U))$ is open. 
\qed

\medskip
Let $Z'\subset Z$ be the union of all elements of ${\mathcal D}$ which are not singletons,
and denote by ${\mathcal R}'$ the equivalence relation on $Z'$ induced by ${\mathcal R}$.
\begin{lem}\label{lem:elem-deco}
Suppose that $Z'$ is closed.
Then ${\mathcal R}$ is usc iff  ${\mathcal R}'$ is usc. 
\end{lem}
\proof
Suppose that ${\mathcal R}'$ is usc. 
Let $D\in {\mathcal D}$. If $D$ is a singleton, then $Z-Z'$ is a saturated open neighborhood of $D$.
On the other hand, if $D\subset Z'$ then $D$ has a saturated open neighborhood $V'$ in $Z'$.
It is an intersection $V'=V\cap Z'$ with an open subset $V\subset Z$ which is necessarily again saturated. 
This verifies that ${\mathcal R}$ is usc. 

Vice versa, suppose that ${\mathcal R}$ is usc.
Then the intersection of a saturated open subset in $Z$ with $Z'$ is open and saturated in $Z'$.
Hence ${\mathcal R}'$ is usc. 
\qed

\medskip
We will use the following result:
\begin{prop}
[{\cite[Proposition 2, page 13]{Daverman}}] 
\label{prop:metrization}
If $Z$ is metrizable and ${\mathcal R}$ is a compact usc decomposition of $Z$, then $Z/{\mathcal R}$ is again metrizable. 
\end{prop}
We now apply the notion of usc decompositions in the context of the Finsler thickening $\hat\La$ 
of $\Lat\subset\Flagt$.
Since $\Ga$ is $\taumod$-antipodal and the thickening $\Th$ is slim,
we obtain a compact decomposition ${\mathcal R}$ of $\hat\Si$, whose elements are singletons,
namely the points in $\hat\Om$, 
and the thickenings $\Th^{Fins}(\tau)$ of the simplices $\tau \in \Lat$. 
(One can show that the latter are contractible, 
cf.\ Remark~\ref{rem:thcontr}.)
We let 
$$
\kappa: \hat\Si\to \Si
$$
denote the quotient projection, and
$$\La:=\kappa(\hat\La)\cong\Lat, \quad \Om:=\kappa(\hat\Om)\cong\hat\Om .$$

\begin{lem}
The  decomposition ${\mathcal R}$ of $\hat\Si$ is compact usc. 
\end{lem}
\proof The restriction $\hat\La\to \La$ of $\kappa$ is a map of compact Hausdorff spaces 
and hence closed.
Thus the restriction of the decomposition ${\mathcal R}$ to $\hat\La$ is usc,
cf.\ Lemma~\ref{lem:uschar}.
Hence, by Lemma \ref{lem:elem-deco}, the decomposition ${\mathcal R}$ is usc as well. 
It is also compact.
\qed 

\begin{cor}
$\Si=\hat\Si/{\mathcal R}$ is metrizable.
\end{cor}
This corollary is relevant to us in order to do computations with \v{C}ech cohomology.

\begin{rem}
\label{rem:isoch}
We showed in the second version of this paper
(see Lemma~10.7 there) 
that $\Si$ is \v{C}ech acyclic,
compare Remark~\ref{rem:thcontr}.
\end{rem}

\subsection{Convergence action}\label{sec:con} 

Suppose that $\Ga< G$ is $\taumod$-regular and $\taumod$-antipodal.
We continue using the notation from the previous section. 
The action of $\Ga$ on $\hat\Si$ descends to a continuous action 
$\Ga\acts\Si$.

\begin{lem}\label{lem:ca}
$\Ga\acts\Si $ is a convergence action. 
\end{lem}
\proof Let $\ga_n\to \infty$ be a sequence in $\Ga$. Since the group $\Ga< G$ is a $\taumod$-convergence subgroup, 
we may assume after extraction 
that $(\ga_n)$ is $\taumod$-contracting: There exist simplices $\tau_{\pm}\in \Flagpmt$ such that $\ga_n\to \tau_+$ uniformly on compacts in $C(\tau_-)$. We claim that $\ga_n$ converges uniformly on compacts in $\Si - \kappa(\tau_-)$ to $\la_+=\kappa(\tau_+)$. It suffices to show that for each sequence $z_n\in  \Si$ converging to $z\ne \la_-=\kappa(\tau_-)$, 
we have after extraction that 
\begin{equation}\label{eq:limgn}
\lim_{n\to\infty} \ga_nz_n= \la_+. 
\end{equation}
Take $\zeta_n\in \kappa^{-1}(z_n)$. Then, after extraction, $\zeta_n\to \zeta\notin \Th^{Fins}(\tau_-)$. 
According to Corollary \ref{cor:drcontrsq}, the accumulation set of the sequence $\ga_n\zeta_n$ is contained in 
$ \Th^{Fins}(\tau_+)$. This implies \eqref{eq:limgn}. \qed 

\medskip 
We note that in view of Theorem \ref{thm:pdwconv} the group $\Ga$ acts properly discontinuously on 
$\Om$. It is also clear that $\La$ is the limit set of the action $\Ga\acts \Si$. Since $\hat\Om$ is path-connected, so is $\Om$. 

\subsection{Cocompactness}

We now make the stronger assumption that $\Ga<G$ is $\taumod$-asymptotically embedded. 
Continuing the discussion of the previous section,
we then also have an equivariant homeomorphism $\geo \Ga \to \Lat \to \La$. 
Thus, Theorem   \ref{thm:haiss} together with Lemma \ref{lem:ca} imply
that the action of $\Ga\acts\Om$ is cocompact.
Therefore the action of $\Ga\acts\hat\Om$ is cocompact as well. 
By combining this with Theorem  \ref{thm:pdwconv}, we obtain the main result of this paper:

\begin{thm}\label{thm:main-coco}
Let $\Ga< G$ be a $\taumod$-asymptotically embedded subgroup, 
and let $\Th\subset W$ be a $W_{\taumod}$-invariant balanced thickening.
Then the action
$$
\Ga\acts \ol{X}^{Fins} - \Th^{Fins}(\Lat)  
$$
is properly discontinuous and cocompact. 
The quotient 
$$
\bigl(\ol{X}^{Fins} - \Th^{Fins}(\Lat)\bigr)/\Ga
$$
has a natural structure as a compact real-analytic orbifold with corners. 
\end{thm}

\begin{rem}
\label{rem:pfcmp}
The starting point of our proof of Theorem~\ref{thm:main-coco}, 
namely  the usage of the bundles  $E$ and $E_{mod}$ in the proof of Theorem \ref{thm:extended_cocompactness}, 
is similar to the one in \cite[Prop.\ 8.10]{GW}. 
However, we avoid the use of Poincar\'e duality and do not need homological assumptions on the space $\Si$.
An essential ingredient in our proof is the map of triads (\ref{eq:triads}),
i.e.\ the existence of a continuous extension of the equivariant proper map $\tilde f: Y\to\Om$ to a map of compactifications. \end{rem}

\begin{rem}[Cocompactness on $\DF X$]
Intersecting the domain in the theorem with $\DF X$ yields that the corresponding actions 
$$
\Ga\acts \DF X - \ThF(\Lat)  
$$
are cocompact,
thus recovering some of the cocompactness results obtained in \cite{coco15}.
\end{rem}

\subsection{An example of a non-regular discrete subgroup}

We now consider a simple example of a {\em non-regular} discrete subgroup
and show that the action on the domain of proper discontinuity constructed earlier in Theorem~\ref{thm:pdwdsc}
is cocompact.

Let $G=PSL(3,\R)$, 
and let $\Ga\cong \Z^2$ be a discrete subgroup of transvections
preserving a maximal flat $F\subset X$. 
As in Theorem \ref{thm:pdwdsc}, 
we choose a {\em multi-thickening}
$\Th_{\bullet}$, i.e. a collection of $\Wt$-left invariant thickenings $\emptyset\neq\Th_{\taumod}\subsetneq W$ 
for all face types $\taumod\subseteq\simod$.
In addition, we require that 
$$
\Th_{\iota\taumod}= \Th^c_{\taumod}.
$$ 
for all $\taumod$.
In particular, $\Th_{\simod}$ is the unique balanced thickening. 
To simplify notation, we drop the face type index,
and write $\Th(\tau)$ instead of $\Th_{\taumod}(\tau)$ for simplices $\tau\in\Flagt$.

There are exactly two distinct multi-thickenings $\Th_\bullet$ 
which are swapped by the involution $\iota:\amod\to\amod$. 
We let $p_{mod}$
denote the vertex of $\simod$ of the type {\em point} in terms of the projective
incidence geometry associated with the group $G$, and let $l_{mod}$ denote the
other vertex of $\simod$ of the type {\em line}. 
We require
$\Th(p_{mod})\subset a_{mod}$ to consist of the two chambers containing
$p_{mod}$ as a vertex. In other words, for a point $p\in \R P^2$ the thickening 
$\ThF(p)$ of $p$ in the full flag manifold $\DF X=\Flags$ consists of all the
flags $(p, l)$.  
Accordingly, for each line $l$ in the projective plane, its thickening
$\ThF(l)\subset\DF X$ consists of all the flags $(p, l')$ where $p\in
\R P^2$ are points incident to the line $l$. Topologically speaking, 
$\ThF(p)\cong S^1$ while $\ThF(l)$ is the 2-torus, the trivial
circle bundle over $l$ whose fibers are the thickenings $\ThF(p)$,
$p\in l$. 

We let $p_i, l_i$ for $i=1, 2, 3$ denote the singular points in $\geo F$, where
$p_i\in \R P^2$ are the fixed points of $\Ga$  and $l_i\in (\R P^2)^\vee$ are
the fixed lines of $\Ga$, labeled so that 
$l_i$ is the line through $p_{i-1}$ and $p_{i+1}$ (where $i$ is taken mod 3). 

We obtain:
$$
\La^{+}_{p_{mod}} =\{p_1, p_2, p_3\}= \La^{-}_{l_{mod}} ,
 \quad \La^{-}_{p_{mod}}=\{l_1, l_2, l_3\}= \La^{+}_{l_{mod}}.  
$$
Furthermore, 
$$
\La^{\pm}_{\simod}= \{(p_i, l_j):i\neq j\}.  
$$
Thus,
$$
\La_{\bullet}= \bigcup_{\taumod} \La^{\pm}_{\taumod}= 
\{p_1, p_2, p_3\} \cup \{l_1, l_2, l_3\} \cup \La^{\pm}_{\simod}. 
$$
We have:
$$
\ThF(\La^{+}_{p_{mod}})= \ThF^c(\La^{-}_{l_{mod}}), \quad 
\ThF(\La^{+}_{l_{mod}})= \ThF^c(\La^{-}_{p_{mod}}),
$$
while
$$
\ThF(\La^{+}_{\simod})= \ThF^c(\La^{-}_{\simod}). 
$$
Therefore, the union
$$
\ThF(\La_{\bullet})= \bigcup_{\taumod} (\ThF^c(\La^-_{\taumod})
\cup \ThF(\La^+_{\taumod}))\subset =\ThF(\La^+_{l_{mod}})=\DF X
$$
is the set of all flags $(p,l)$ such that $p$ is incident to one of the lines
$l_1, l_2, l_3$. Topologically speaking, this set is the union of three trivial
circle bundles $\ThF(l_i)$ over the circles $l_i$, such that 
$$
\ThF(l_{i-1})\cap \ThF(l_{i+1})= \ThF(p_i), 
$$
where $i$ is taken modulo 3. 

The fact that the action 
\begin{equation}\label{eq:domain-for-Z^2}
\Ga\acts\OmF:= \DF X - \ThF(\La_{\bullet}) 
\end{equation}
is properly discontinuous can be seen as a special case of \cite[Proposition 6.21]{coco15} in our earlier work. 
It can be also be seen directly by observing that the action of $\Ga$ is properly discontinuous on $\Om_{p_{mod}}:=\R P^2 - (l_1 \cup l_2 \cup l_3)$, since \eqref{eq:domain-for-Z^2} is the preimage of 
$\Om_{p_{mod}}$ under the fibration $\DF X\to \R P^2$. Since $\Ga$ acts cocompactly on  
$\Om_{p_{mod}}$ (the quotient is the disjoint union of four 2-tori), the group $\Ga$ also acts cocompactly on the domain 
\eqref{eq:domain-for-Z^2}.  The quotient $\OmF/\Ga$ is a circle bundle over $\Om_{p_{mod}}/\Ga$. 

We now discuss the corresponding Finsler thickening 
$$
\Th^{Fins}(\La_{\bullet})=\bigcup_{\taumod}\bigl((\Th^c)^{Fins}(\Latm)\cup\Th^{Fins}(\Latp)\bigr) 
$$ 
cf. Theorem \ref{thm:pdwdsc}. For each point $p_i$, the Finsler thickening of $p_i$ is the closed 2-disk, which is
the closed stratum $\ol{X}_{p_i}$ naturally isomorphic to a compactified hyperbolic
plane whose ideal boundary is the circle $\ThF(p_i)$. 

For each $l_i$, its Finsler thickening $\Th^{Fins}(l_i)$ is the union of a solid torus,
whose boundary is the torus $\ThF(l_i)$, 
and the closed 2-disk $\ol{X}_{l_i}$, whose boundary circle $C_i$ is the set of flags $(p,l_i)$ for $p\in l_i$.
The circle $C_i\subset \ThF(l_i)$ is a section 
of the circle bundle $\ThF(l_i)\to l_i$. In particular, $\Th^{Fins}(l_i)$ is 
contractible. 

Furthermore, we have
$$
\Th^{Fins}(l_{i-1}) \cap \Th^{Fins}(l_{i+1})= \Th^{Fins}(p_i),
$$
and the triple intersection of the $\Th^{Fins}(l_i)$ is empty. 
The thickening $\Th^{Fins}(\La_{\bullet})$ equals the union
$$
\bigcup_{i=1}^3 \Th^{Fins}(l_i)
$$
which is homotopy-equivalent to the circle.
The inclusion
$$
\geo^{Fins}(F)\embed \Th^{Fins}(\La_{\bullet})
$$
is a homotopy equivalence. This inclusion is the restriction of the natural embedding
$$
\ol{F}^{Fins} \embed \ol{X}^{Fins}. 
$$
We are now in the position to apply Theorems \ref{thm:pdwdsc} and \ref{thm:extended_cocompactness},  
taking $\Ga\acts\ol{F}^{Fins}$ as the compactified model action,
and conclude:

\begin{prop}
The action 
$$\Ga \acts \ol X^{Fins}  - \Th^{Fins}(\La_{\bullet}) $$
is properly discontinuous and cocompact.
\end{prop}

\section{Characterizations of $\taumod$-Anosov subgroups}\label{sec:ecRCA}

In our earlier papers \cite{morse,mlem} 
we gave various characterizations of Anosov subgroups 
in terms of dynamics and coarse extrinsic geometry, 
see also our survey \cite{anosov}.
The most relevant characterizations for this paper are {\em asymptotically embedded},
see Definition~\ref{def:asyemb} above, and {\em URU}.

We assume from now on that $\taumod$ is $\iota$-invariant. 
A discrete subgroup $\Ga<G$ is called {\em $\taumod$-URU}
if it is uniformly $\taumod$-regular and undistorted, 
cf.\ \cite{mlem}.
We proved in \cite{mlem} that the $\taumod$-URU property is equivalent to $\taumod$-Anosov.
In this section, we will give further characterizations of the Anosov property
in terms of dynamics (S-cocompactness) 
and coarse extrinsic geometry (Finsler quasiconvexity and existence of retractions).

\subsection{Finsler quasiconvexity}
\label{sec:fqc}

In this section we introduce the notion of {\em Finsler quasiconvex subgroups} of $G$, which mimics the notion of quasiconvex subgroups of hyperbolic groups. Recall that a subgroup $\Ga$ of a word hyperbolic group $\Ga'$ is called {\em quasiconvex} if discrete geodesic segments in $\Ga'$ with endpoints in $\Ga$ are uniformly close to $\Ga$.   

Fix a type $\bar\theta\in\interior(\taumod)$. Recall that  
$d^{\bar\theta}$ is in general only  a {\em pseudometric} on $X$. 

\begin{definition}\label{defn:quasiconvex}
A discrete subgroup $\Ga< G$ is {\em $\taumod$-Finsler quasiconvex}
if for each $x\in X$ there is a constant $R<+\infty$ 
such that any two points in $\Ga x$ can be connected by a $d^{\bar\theta}$-geodesic segment
contained in the $R$-neighborhood $N_R(\Ga x)$ of $\Ga x$
with respect to $d^{\bar\theta}$. 
\end{definition}

\begin{prop}\label{prop:FC=URU}
A uniformly $\taumod$-regular subgroup $\Ga< G$ is $\taumod$-Finsler quasiconvex iff it is $\taumod$-URU. 
\end{prop}
\proof 
We first reduce the assertion to the case when the pseudo-metric $d^{\bar\theta}$ is a metric. 
We recall, cf.\ the end of section~\ref{sec:finsdist},
that $X$ splits as a product $X_1\times X_2$ such that $d^{\bar\theta}$ is degenerate precisely in the $X_2$-direction
and induces a metric on $X_1$.
In particular, $\bar\theta$ points in the $X_1$-di\-rec\-tion,
i.e.\ the visual boundary points of type $\bar\theta$ are contained in $\geo X_1$.
Then $d^{\bar\theta}$-balls split off $X_2$-factors,
i.e.\ they are products of $X_2$ with $d^{\bar\theta}$-balls in $X_1$.
The same applies to $\taumod$-Weyl cones and $\taumod$-diamonds.
A map $I\to X$ from an interval is a $d^{\bar\theta}$-geodesic
iff its projection to $X_1$ is a $d^{\bar\theta}$-geodesic.
Furthermore, a map into $X$ is a quasiisometric embedding with respect to $d^{\bar\theta}$
iff its $X_1$-component is.
We can therefore assume that $d^{\bar\theta}$ is a metric. 

Suppose now that $\Ga< G$ is $\taumod$-Finsler quasiconvex. 
The closed $R$-neighborhood $N_R(\Ga x)$ is path-connected and $\Ga$ acts cocompactly on it.
Therefore, $\Ga$ is finitely generated and the orbit map 
$$
o_x: \Ga \to N_R(\Ga x)$$
is a quasiisometric embedding, where we equip $N_R(\Ga x)$ with a path-metric induced by $d^{\bar\theta}$. Since the metrics $d^{Riem}$ and $d^{\bar\theta}$ on $X$ are equivalent, 
the definition of Finsler quasiconvexity implies that the inclusion map
$$
N_R(\Ga x)\to (X, d^{Riem}) 
$$
is a quasiisometric embedding. Therefore, $\Ga< G$ is undistorted. Since $\Ga$ was assumed to be 
uniformly $\taumod$-regular, it is $\taumod$-URU. 

The converse direction follows from our Morse Lemma \cite[Thm.\ 1.3]{mlem}
and the description of the geometry of $d^{\bar\theta}$-geodesics, cf.\ section~\ref{sec:geodesics}.
\qed

\subsection{S-cocompactness and retractions}
\label{sec:sccrt}

We call an open subset $\Om\subset\geo^{Fins}X$ {\em saturated}
if it is a union of small strata $X_{\nu}$.

We start with the following simple observation about Finsler convergence at infinity: 
If $(x_n)$ and $(y_n)$ are sequences in $X$ which are bounded distance apart
(i.e.\ $d(x_n, y_n)$ is uniformly bounded) 
and $x_n\to [b], y_n\to [b']\in \geo^{Fins}X$, then 
the limit points $[b]$ and $[b']$ lie in the same small stratum $X_{\nu}$,
see Lemma~\ref{lem:limsamstr}.
In particular, for each saturated open subset $\Om\subset\geo^{Fins}X$,
$$[b]\in \om \Leftrightarrow [b']\in \Om .$$
It follows that if $[b]\in \Om$, then  the entire accumulation set 
of the sequence of balls $B(x_n,R)$
$$
\Acc((B(x_n,R))) \subset \geo^{Fins}X 
$$ 
is a compact subset of $\Om$.

\begin{lem}\label{lem:uniformly finite}
Let $\Ga<G$ be a discrete subgroup.
Suppose that $\Om\subset\geo^{Fins}X$ is a $\Ga$-invariant saturated open subset
such that the action
$$\Ga\acts X\sqcup \Om$$
is properly discontinuous.
Then each compact subset $C\subset X\sqcup \Om$
satisfies the following {\em uniform finiteness property}: 
There exists a function $f_C(R)$ such that for each ball $B(x, R)\subset X$ it holds that 
$$
\card(\{\ga\in \Ga: \ga C \cap B(x,R) \ne \emptyset\})\le f_C(R). 
$$
\end{lem}
\proof
Suppose the contrary.
Then there is a sequence of balls $B(x_i,R)$ intersecting $C$
and a sequence $\ga_i\to\infty$ in $\Ga$ 
such that also the balls $B(\ga_ix_i,R)$ intersect $C$.
We may assume after extraction that $x_i\to\bar x$ and $\ga_ix_i\to\bar x'$ in $\ol X^{Fins}$.
By the observation preceeding the lemma, it holds that $\bar x,\bar x'\in X\sqcup\Om$.
Since these points are dynamically related with respect to the $\Ga$-action,
we obtain a contradiction with proper discontinuity.
\qed

\medskip
The lemma leads to the following definition. 
\begin{definition}\label{def:scc}
A discrete subgroup $\Ga< G$ is {\em S-cocompact} if there exists a $\Ga$-invariant saturated open subset 
$\Om\subset \geo^{Fins}X$ such that the action $$\Ga\acts X\sqcup\Om$$
is properly discontinuous and cocompact.
\end{definition}
Note that each S-cocompact subgroup is necessarily finitely generated
because it acts properly discontinuously and cocompactly on a connected manifold with boundary.
\begin{thm}\label{thm:FC->U}
Each S-cocompact subgroup $\Ga< G$ admits a $\Ga$-equivariant coarse Lipschitz retraction $r: X\to \Ga$. 
In particular,  $\Ga$ is undistorted in $G$. 
\end{thm}
\proof Let $\Om\subset\geo^{Fins}X$ be as in the definition. 
Let $C\subset X\sqcup\Om$ be a compact subset 
whose $\Ga$-orbit covers the entire $X\sqcup\Om$. 
We define the coarse retraction $r$ first by sending each point $x\in X$ 
to the subset 
$$r(x) := \{ \ga\in \Ga : x\in\ga C  \} \subset \Ga .$$
This subset is clearly finite because of the proper discontinuity of the $\Ga$-action,
and the assignment $x\mapsto r(x)$ is equivariant.
According to Lemma \ref{lem:uniformly finite},
the cardinality of the subset 
$$
\{\ga\in \Ga: \ga\in r(B(x,1))\} = \{ \ga\in\Ga : B(x,1)\cap\ga C\neq\emptyset \}
$$ 
is bounded by $f_C(1)$, independently of $x$.
It follows that $r$ is coarse Lipschitz.
\qed

\medskip
We now apply the previous theorem to the cocompact domains of proper discontinuity 
obtained earlier by removing Finsler thickenings of limit sets. 
The next result relates conicality and S-cocompactness:
\begin{thm}
Suppose that $\Ga< G$ is uniformly $\taumod$-regular and $\taumod$-antipodal. 
Then $\Ga$ is $\taumod$-Anosov iff it is S-cocompact.
\end{thm}
\proof We use that $\taumod$-Anosov is equivalent to $\taumod$-asymptotically embedded.
That $\taumod$-asymptotically embedded implies S-cocompact 
is our main result Theorem~\ref{thm:main-coco}.  
To prove the converse, note that each S-cocompact subgroup is undistorted in $G$ 
by Theorem \ref{thm:FC->U}. 
Hence, $\Ga$ is $\taumod$-URU, and therefore $\taumod$-Anosov.
\qed 

\medskip
We now can prove a converse to Theorem~\ref{thm:main-coco}:
\begin{cor}
Suppose that $\Ga<G$ is uniformly $\taumod$-regular 
and that $\Th\subset W$ is a $W_{\taumod}$-invariant balanced thickening. 
Then the following are equivalent:

(i) The properly discontinuous action (see Theorem~\ref{thm:pdwdsc})
$$ \Ga\acts \ol X^{Fins} - \Th^{Fins}(\Lat)$$
is cocompact.

(ii) $\Ga$ is S-cocompact.

(iii) $\Ga$ is $\taumod$-Anosov.
\end{cor}
\proof
(i)$\Ra$(ii) is obvious.

(ii)$\Ra$(iii):
$\Ga$ is S-cocompact,
hence $\taumod$-URU by Theorem~\ref{thm:FC->U},
and therefore $\taumod$-Anosov.

(iii)$\Ra$(i):
Since $\taumod$-Anosov is equivalent to $\taumod$-asymptotically embedded,
the implication is the content of Theorem~\ref{thm:main-coco}.
\qed 

\medskip
We are now ready to state the equivalence of 
a variety of conditions on discrete subgroups,
extending the list of equivalent conditions from \cite{morse,mlem}.
\begin{thm}
The following are equivalent for uniformly $\taumod$-regular subgroups $\Ga< G$:

1. $\Ga$ is a equivariant coarse  retract. 

2. $\Ga$ is a coarse retract. 

3.  $\Ga$ is undistorted in $G$, i.e. $\taumod$-URU. 

4. $\Ga$ is $\taumod$-Finsler quasiconvex.

5. $\Ga$ is $\taumod$-asymptotically embedded. 

6. $\Ga$ is S-cocompact. 

7. $\Ga$ is $\taumod$-Anosov. 

\end{thm}
\proof The implications 1$\Ra$2$\Ra$3 are immediate. 
The equivalence 3$\Leftrightarrow$4 is proven in Proposition 
\ref{prop:FC=URU}. 
The equivalence 3$\Leftrightarrow$5 is one of the main results of \cite{mlem}, see Corollary 1.6 of that paper. The equivalence 5$\Leftrightarrow$7 is established in \cite{morse}. 
The implication 5$\Ra$6 is the main result Theorem \ref{thm:main-coco} of this paper, while the implication 6$\Ra$1 is established in Theorem \ref{thm:FC->U}. \qed

\medskip
We note that this list of equivalences is nearly a perfect match to the list of equivalent definitions of convex cocompact subgroups of rank 1 Lie groups, except that convex-cocompactness is (by necessity) missing, see \cite{convcoco}.

\end{document}